\documentclass[11pt,reqno]{amsart}
\makeatletter
\@ifundefined{subjclassname@2020}{%
  \@namedef{subjclassname@2020}{\textup{2020} Mathematics Subject Classification}}{}
\makeatother
\usepackage[T1]{fontenc}
\usepackage{lmodern}
\usepackage{amsmath,amssymb,mathtools,mathrsfs}
\usepackage[letterpaper,textwidth=6.35in,textheight=9.25in,centering]{geometry}
\usepackage{microtype}
\usepackage{enumitem}
\usepackage{xcolor}
\usepackage{needspace}
\usepackage{cases}
\usepackage[hidelinks]{hyperref}
\hypersetup{pdftitle={Asymptotic behaviors of positive solutions to critical fractional p-Laplace equations with Hardy potential},pdfauthor={Jungang Li, Yansheng Shen and Zeming Wang},pdfkeywords={Fractional p-Laplacian, Hardy potential, critical exponent, asymptotic behavior, comparison principle}}
\providecommand{\uppi}{\pi}
\ifdefined\fint\else
  \def\Xint#1{\mathchoice
    {\XXint\displaystyle\textstyle{#1}}%
    {\XXint\textstyle\scriptstyle{#1}}%
    {\XXint\scriptstyle\scriptscriptstyle{#1}}%
    {\XXint\scriptscriptstyle\scriptscriptstyle{#1}}%
    \!\int}
  \def\XXint#1#2#3{{\setbox0=\hbox{$#1{#2#3}{\int}$}\vcenter{\hbox{$#2#3$}}\kern-.5\wd0}}
  \def\fint{\Xint-}
\fi

\usepackage[normalem]{ulem}
\newif\ifCommandCAfound
\IfFileExists{D:/OneDrive/Desktop/LaTeXFiles/Command/CommandCA.tex}{\CommandCAfoundtrue}{\CommandCAfoundfalse}
\makeatletter
\ifCommandCAfound
  \input{D:/OneDrive/Desktop/LaTeXFiles/Command/CommandCA.tex}
\else
  \newcommand{\mk@box}[1]{\mbox{#1}}
  \newcommand{\mk@armour}{%
    \let\mk@@cite\cite
    \let\mk@@ref\ref
    \let\mk@@eqref\eqref
    \renewcommand{\cite}[2][]{%
      \mk@box{\if\relax\detokenize{##1}\relax\mk@@cite{##2}\else\mk@@cite[##1]{##2}\fi}}%
    \renewcommand{\ref}[1]{\mk@box{\mk@@ref{##1}}}%
    \renewcommand{\eqref}[1]{\mk@box{\mk@@eqref{##1}}}}
  \newcommand{\mkblue}[1]{\textcolor{blue}{\begingroup\mk@armour\dashuline{#1}\endgroup}}
  
  \newcommand{\blkendblue}{\textcolor{blue}{\ensuremath{\blacktriangleleft}}}
  \newcommand{\blkblue}{\color{blue}\ensuremath{\blacktriangleright}\aftergroup\blkendblue}
  \newcommand{\mkred}[1]{\textcolor{red}{\begingroup\mk@armour\uline{#1}\endgroup}}
  \newcommand{\blkendred}{\textcolor{red}{\ensuremath{\blacktriangleleft}}}
  \newcommand{\blkred}{\color{red}\ensuremath{\blacktriangleright}\aftergroup\blkendred}
\fi
\makeatother

\numberwithin{equation}{section}
\newtheorem{theorem}{Theorem}[section]
\newtheorem{lemma}[theorem]{Lemma}
\newtheorem{proposition}[theorem]{Proposition}
\newtheorem{corollary}[theorem]{Corollary}
\theoremstyle{definition}
\newtheorem{definition}[theorem]{Definition}
\theoremstyle{remark}
\newtheorem{remark}[theorem]{Remark}
\theoremstyle{plain}
\setenumerate[1]{itemsep=0pt,parsep=0pt,label={\textnormal{(\roman*)}}}
\allowdisplaybreaks[2]
\newcommand{\mbN}{\mathbb{N}}

\newcommand{\mbR}{\mathbb{R}}

\newcommand{\mcH}{\mathcal{H}}

\DeclareMathOperator{\supp}{supp}

\DeclareMathOperator{\dist}{dist}

\newcommand{\ve}{\varepsilon}
\newcommand{\vp}{\varphi}
\newcommand{\ol}{\overline}
\newcommand*{\diop}{\mathop{}\mspace{-3mu}\mathrm{d}}

\newcommand{\diopx}{\diop x}
\newcommand{\diopy}{\diop y}

\newcommand{\norm}[1]{\lVert #1\rVert}

\newcommand{\wtilde}[1]{\widetilde{#1}}

\newcommand{\mbRN}{{\mathbb{R}^N}}
\newcommand{\mbS}{\mathbb{S}^{N-1}}
\newcommand{\mcD}{\mathcal{D}^{s,p}(\mathbb{R}^N)}
\newcommand{\mcI}{\mathcal{I}}

\newcommand{\pst}{p_s^\ast}
\def\ds{\displaystyle}
\DeclareMathOperator*{\esssup}{ess\,sup}
\DeclareMathOperator*{\essinf}{ess\,inf}
\DeclareMathOperator{\Tail}{Tail}

\newcommand{\diopt}{\diop t}
\newcommand{\diopsigma}{\diop \sigma}
\newcommand{\semi}[1]{[#1]_{s,p}}

\title[Critical fractional $p$-Laplace equations with Hardy potential]{Asymptotic behaviors of solutions to the critical fractional $p$-Laplace equations with Hardy potential}
\author{Jungang Li}
\address{School of Mathematical Sciences, University of Science and Technology of China, Hefei 230026, The People's Republic of China}
\email{jungangli@ustc.edu.cn}
\author{Yansheng Shen}
\address{School of Mathematical Sciences, Jiangsu University, Zhenjiang 212013, The People's Republic of China}
\email{ysshen@mail.bnu.edu.cn}
\author{Zeming Wang}
\address{School of Mathematics and Statistics, Beijing Institute of Technology, Beijing 102488, The People's Republic of China}
\email{zmwangmath@bit.edu.cn}
\date{}
\begin{document}
\begin{abstract}
In this paper we consider positive weak solutions $u\in\mcD$ to the critical fractional $p$-Laplace equation with Hardy potential in $\mathbb{R}^{N}$
\begin{equation*}
    (-\Delta_p)^s u-\mu\frac{u^{p-1}}{|x|^{ps}}=u^{\pst-1} \ \ \ \text{in\ } \mbRN,
\end{equation*}
where $N\geqslant 2$, $s\in(0,1)$, $p\in(1,\infty)$ satisfying $ps<N$, $\pst=\frac{Np}{N-ps}$ and $\mu\in(0,\mathcal{C}_{N,s,p})$, $\mathcal{C}_{N,s,p}$ denotes the sharp constant of the fractional Hardy inequality. By means of a suitable comparison principle for weak subsolutions and supersolutions of the equation $(-\Delta_p)^sf= \big(\frac{\mu}{|x|^{ps}}+a(x)\big)f^{p-1}$ on general open set $\Omega$, with nonzero exterior data and $a\in L^{\frac{N}{ps}}(\mbRN)$, combined with careful constructions of the weak subsolutions and supersolutions, we prove the following asymptotic behaviors of positive weak solution $u$ to the above considered problem
\begin{equation*}
  C^{-1}\leqslant |x|^{\mu_1}\left(1+|x|\right)^{\mu_2-\mu_1}u(x)\leqslant C
  \quad\text{for a.e. } x\in\mbRN,
\end{equation*}
where $\mu_1<\mu_2$ are the two roots of the indicial equation associated with the Hardy operator and $C\in[1,\infty)$ depends on $u$. Also, we show that the constant $C$ cannot be chosen to depend only on $N$, $s$, $p$ and $\mu$.

\end{abstract}
\maketitle
\begingroup
\raggedbottom
\tableofcontents
\clearpage
\endgroup

\section{Introduction and main results }

This work is devoted to the study of the following critical fractional $p$-Laplacian equation involving a Hardy potential
\begin{equation}\label{equ:EL}
  \begin{cases}
    (-\Delta_p)^s u-\mu\frac{u^{p-1}}{|x|^{ps}}=u^{\pst-1} &\text{in\ } \mbRN,\\
    u>0 &\text{in\ } \mbRN,\\
    u\in\mcD,
  \end{cases}
\end{equation}
where $N\geqslant 2$, $s\in(0,1)$, $p\in(1,\infty)$ satisfying $ps<N$, $\pst:=\frac{Np}{N-ps}$ denotes the critical fractional Sobolev exponent, $\mu\in(0,\mathcal{C}_{N,s,p})$ and $\mathcal{C}_{N,s,p}$ is the sharp constant of the fractional Hardy inequality which was given by
Frank and Seiringer in \cite[Theorem 1.1]{bib:FS2008}
\begin{equation}\label{equ:hardyinequality}
\mathcal{C}_{N,s,p}\int_{\mbRN}\frac{|u(x)|^p}{|x|^{ps}}\diopx\leqslant\int_{\mathbb{R}^{N}}\int_{\mathbb{R}^{N}}
\frac{|u(x)-u(y)|^{p}}{|x-y|^{N+ps}}dxdy,  \ \ \ \text{for} \ \text{any} \  u\in\mcD.
\end{equation}
The fractional $p$-Laplace operator $(-\Delta_p)^s$ is defined on smooth functions as
\begin{equation}\label{equ:fractionalplaplacian}
(-\Delta_p)^s u(x):=2\lim_{\ve\to 0^+}\int_{\{y\in\mbRN\colon |x-y|>\ve\}}\frac{J_p\left(u(x)-u(y)\right)}{|x-y|^{N+ps}}\diopy,
\end{equation}
where $J_p(t):=|t|^{p-2}t$, $t\in\mbR$.

The existence, radial symmetry and qualitative properties of positive weak solutions to critical elliptic problems involving Hardy potential are crucial issues and have been widely studied in the past few decades, since these problems naturally appear in the study of Hardy-Sobolev type inequalities.
For the critical quasilinear $p$-Laplace equation
\begin{equation}\label{equ:localEL}
  -\Delta_p u-\mu\frac{u^{p-1}}{|x|^{p}}=u^{p^\ast-1}\quad\text{in } \mbRN,
\end{equation}
with $1<p<N$, $0\leq\mu<(\frac{N-p}{p})^p$ and $p^\ast=\frac{Np}{N-p}$, the existence of positive weak solutions in $\mathcal{D}^{1,p}(\mathbb{R}^{N})$ was proved by Abdellaoui, Felli and Peral \cite{bib:AFP2006}, Filippucci, Pucci and Robert \cite{FPR}, and Musina \cite{Mus}, where among other results, the authors showed the existence of the minimizers of the Hardy-Sobolev quotient
\begin{align*}
S_{p}=\inf_{u\in\mathcal{D}^{1,p}(\mathbb{R}^{N})\backslash\{0\}}\frac{\int_{\mathbb{R}^{N}}|\nabla u|^{p}dx-\mu\int_{\mathbb{R}^{N}}\frac{|u|^{p}}{|x|^{p}}dx}
{\big(\int_{\mathbb{R}^{N}}|u|^{p^{\ast}}dx\big)^{\frac{p}{p^{\ast}}}}.
\end{align*}
When $\mu=0$, the classification of positive weak solutions of \eqref{equ:localEL} has been well investigated by Caffarelli, Gidas and Spruck \cite{CGS} for the case of $p=2$, by Sciunzi \cite{Sc} for $p>2$ and by V\'{e}tois \cite{V} for $1<p<2$, respectively, which in turn yields the explicit form of the extremal functions of the classical Sobolev inequalities, see \cite{bib:Aubin1976,bib:Talenti1976,bib:Lieb1983} for details. In recent years, such result was generalized by Ciraolo, Figalli and Roncoroni \cite{CFR}, where the authors obtained a complete classification
result for positive solutions to the critical anisotropic $p$-Laplacian equations. When $\mu\neq0$, Terracini \cite{Te} proved for the first time the classification result for positive solutions to \eqref{equ:localEL} in the case $p=2$. Precisely, she showed that any solution to this problem is radial and radially decreasing about the origin through the moving-plane technique and the use of the Kelvin transformation. Regarding general $1<p<N$ and $\mu\in(0,(\frac{N-p}{p})^p)$, Abdellaoui, Felli and Peral \cite{bib:AFP2006} proved the asymptotic behavior of radial solutions of \eqref{equ:localEL} via the ODE analysis and then Oliva, Sciunzi and Vaira \cite{bib:OSV2020} proved that every positive weak solution $u\in D^{1,p}(\mbRN)$ of \eqref{equ:localEL} is radially symmetric and radially decreasing about the origin, which yields uniqueness up to scaling and precise asymptotic limits for all positive energy solutions of \eqref{equ:localEL}. Besides, without any symmetry assumption, Xiang \cite{bib:Xiang2015} proved that every positive weak solution $u\in D^{1,p}(\mbRN)$ of \eqref{equ:localEL} is comparable to $|x|^{-\gamma_1}$ near the origin and to $|x|^{-\gamma_2}$ at infinity, where $\gamma_1<\gamma_2$ are the two roots of the algebraic equation $\gamma^{p-1}\left(N-p-(p-1)\gamma\right)=\mu$ which expresses that $|x|^{-\gamma}$ solves the homogeneous equation $-\Delta_p v=\mu|x|^{-p}v^{p-1}$ in $\mbRN\setminus\{0\}$. In a recent work \cite{bib:EMSV2024}, by adapting similar arguments as in \cite{bib:Xiang2015}, Esposito, Montoro, Sciunzi and Vuono established analogous two-sided estimates for positive energy solutions of the anisotropic doubly critical equation. For other related literature, we refer to e.g. \cite{bib:BianchiEgnell1991,bib:CKN1984,CW,CDL,D,bib:GhoussoubYuan2000} and the references therein.

Concerning the fractional case, the existence of positive weak solutions of problem \eqref{equ:EL} was proved by the second author in \cite{bib:Shen2024}, where among other results, the author showed the existence of minimizers of the fractional Hardy--Sobolev quotient
\begin{equation*}
S_{ps}=\inf_{u\in\mathcal{D}^{s,p}(\mathbb{R}^{N})\backslash\{0\}}  \frac{\int_{\mathbb{R}^{N}}\int_{\mathbb{R}^{N}}
\frac{|u(x)-u(y)|^{p}}{|x-y|^{N+ps}}dxdy-\mu\ds\int_{\mbRN}\frac{|u(x)|^p}{|x|^{ps}}\diopx}{\left(\ds\int_{\mbRN}|u(x)|^{\pst}\diopx\right)^{\frac{p}{\pst}}},
\end{equation*}
for $\mu\in(0,\mathcal{C}_{N,s,p})$. When $p=2$ and $\mu\in(0,\mathcal{C}_{N,s,2})$, by employing the moving-plane method, the technique of the ground state representation of Frank, Lieb and Seiringer \cite{bib:FLS2008} together with variational arguments, Dipierro, Montoro, Peral and Sciunzi \cite{bib:DMPS2016} proved the radial symmetry as well as the asymptotic behaviors of positive weak solutions to problem \eqref{equ:EL}. Afterwards, Mallick \cite{bib:Mallick2019} studied a fractional Hardy--Sobolev--Maz'ya inequality for $p=2$, in which the Hardy potential is singular on a subspace of codimension at least $2$. For Hardy coefficients below the sharp constant, he proved the existence of extremal functions of the fractional Hardy--Sobolev--Maz'ya inequality and obtained a global pointwise upper bound for positive weak solutions of the associated Euler--Lagrange equation.
On the other hand, for general $p>1$ and $\mu=0$, Brasco, Mosconi and Squassina \cite{bib:BMS2016} showed that the extremals of the fractional Sobolev inequality, multiplied by $|x|^{\frac{N-ps}{p-1}}$, have a finite nonzero limit at infinity, and later Marano and Mosconi \cite{bib:MM2019} established two-sided bounds with the same decay exponent for the extremals of the fractional Hardy--Sobolev inequality.

It is noteworthy to mention that while writing this paper, we became aware of \cite[Theorem 1.5]{bib:BM2026} where the authors investigated the asymptotic behavior of positive radially symmetric nonincreasing weak solution of problem \eqref{equ:EL} for $p\geqslant 2$.
By means of the radial monotonicity of the positive weak solutions together with the two ground-state substitutions, weighted
Caccioppoli estimates and finite Moser iterations, they obtained upper bounds with the exponents $\mu_1$ and $\mu_2$, and lower bounds with an arbitrarily small loss in the exponents. Their approach completely differs from ours, and for $p\neq 2$, the sharp asymptotic behavior of general positive weak solutions of \eqref{equ:EL}, both near the origin and at infinity, still remains open up to now. In this work we intend to solve this problem.

Taking inspiration from \cite{bib:EMSV2024,bib:Xiang2015}, we shall employ the comparison principle together with the sub- and supersolutions methods to establish the asymptotic behavior of positive weak solutions for problem \eqref{equ:EL}. For this, we use a crucial result in the determination
of subsolutions and supersolutions. To be more precise, for $r\in(0,\frac{N}{p-1})$ and $x\neq0$, the function $|x|^{-r}$ solves the homogeneous fractional Hardy equation
\begin{equation}\label{equ:purepower}
  (-\Delta_p)^s\left(|\cdot|^{-r}\right)(x)=\Lambda(r)|x|^{-ps-r(p-1)},
\end{equation}
where
\begin{equation}\label{equ:Lambdaformula}
  \Lambda(r):=2\int_1^\infty K(\sigma)\left(\sigma^{r}-1\right)^{p-1}\left(\sigma^{N-1-r(p-1)}-\sigma^{ps-1}\right)\diopsigma, \ \ \
  K(\sigma)=\int_{\mbS}\frac{\diop\mcH^{N-1}(\varpi)}{|\sigma e-\varpi|^{N+ps}},
\end{equation}
$e\in\mbS$ is any fixed unit vector, and $\mcH^{N-1}$ is the $(N-1)$-dimensional Hausdorff measure on the unit sphere $\mbS$. Such a result has been given by Abdellaoui and Bentifour in \cite{bib:AB2017} (see also our Lemma \ref{lem:indexfunction} for details), where the authors showed in \cite[Lemma 3.1]{bib:AB2017} that for $0<\mu<\mathcal{C}_{N,s,p}$, there exist $\mu_{1}$ and $\mu_{2}$ such that
\begin{equation}\label{equ:tworoots}
  0<\mu_1<r_H<\mu_2<\frac{N-ps}{p-1}, \ \ \ \text{where} \ \  r_H=\frac{N-ps}{p}.
\end{equation}
and $\Lambda(\mu_{1})=\Lambda(\mu_{2})=\mu$. As a consequence, for every $\mu\in(0,\mathcal{C}_{N,s,p})$, the equation \eqref{equ:purepower} has exactly two roots $\mu_1$ and $\mu_2$ in $(0,\frac{N-ps}{p-1})$.

With the aid of the above results, we are able to establish the asymptotic behaviors of positive weak solutions of \eqref{equ:EL} for all $p\in(1,\infty)$. Let
\begin{equation*}
  m_0(R):=\int_{\{|x|<R\}}u(x)^{\pst}\diopx
  \qquad\text{and}\qquad
  m_\infty(R):=\int_{\{|x|>R\}}u(x)^{\pst}\diopx
\end{equation*}
be the critical masses of $u$ near the origin and near infinity, respectively. Our main result is the following

\begin{theorem}\label{thm:twosided}
  Let $N\geqslant 2$, $s\in(0,1)$, $p\in(1,\infty)$ with $ps<N$, and $\mu\in(0,\mathcal{C}_{N,s,p})$. Assume that $\mu_1<\mu_2$ are the two roots in \eqref{equ:tworoots}, and $u$ is a positive weak solution of \eqref{equ:EL}. Then there exists a constant $C\in[1,\infty)$, depending on $u$, such that the following hold:
  \begin{enumerate}
    \item\label{item:twosidedpointwise} For almost every $x\in\mbRN$,
    \begin{equation}\label{equ:twosidedglobal}
      C^{-1}\leqslant |x|^{\mu_1}\left(1+|x|\right)^{\mu_2-\mu_1}u(x)\leqslant C.
    \end{equation}
    \item\label{item:twosidedmass} For any $R\in(0,1]$,
    \begin{equation}\label{equ:massorigin}
      C^{-1}R^{N-\mu_1\pst}\leqslant m_0(R)\leqslant C R^{N-\mu_1\pst},
    \end{equation}
    and for any $R\in[1,\infty)$,
    \begin{equation}\label{equ:massinfinity}
      C^{-1}R^{-(\mu_2\pst-N)}\leqslant m_\infty(R)\leqslant C R^{-(\mu_2\pst-N)}.
    \end{equation}
  \end{enumerate}
\end{theorem}

\begin{remark}\label{rem:twosided}
Observe that near the origin the factor $(1+|x|)^{\mu_2-\mu_1}$ in \eqref{equ:twosidedglobal} is comparable to $1$, and at infinity it is comparable to $|x|^{\mu_2-\mu_1}$. Thus Theorem \ref{thm:twosided}\ref{item:twosidedpointwise} contains, with $C$ replaced by $2^{\mu_2-\mu_1}C$, the estimates
\begin{gather*}
  C^{-1}|x|^{-\mu_1}\leqslant u(x)\leqslant C|x|^{-\mu_1}\quad\text{for a.e. } 0<|x|<1,\\
  C^{-1}|x|^{-\mu_2}\leqslant u(x)\leqslant C|x|^{-\mu_2}\quad\text{for a.e. } |x|>1.
\end{gather*}
\end{remark}

The constant $C$ in Theorem \ref{thm:twosided} depends on the solution $u$. Indeed, for $t\in(0,\infty)$ let us set
\begin{equation*}
  u^{(t)}(x):=t^{r_H}u(tx),\qquad x\in\mbRN.
\end{equation*}
\begin{proposition}\label{prop:scaling}
  Let $N$, $s$, $p$, $\mu$, $\mu_1$, and $\mu_2$ be as in Theorem \ref{thm:twosided}, and let $u$ be a positive weak solution of \eqref{equ:EL}. Then, for any $t\in(0,\infty)$, $u^{(t)}$ is a positive weak solution of \eqref{equ:EL} with $\semi{u^{(t)}}=\semi{u}$. Moreover, for any $\rho\in(0,\infty)$,
  \begin{equation}\label{equ:scalingorigin}
    \lim_{t\to\infty}\esssup_{\{0<|x|<\rho\}}|x|^{\mu_1}u^{(t)}(x)=\infty
    \qquad\text{and}\qquad
    \lim_{t\to 0^+}\essinf_{\{0<|x|<\rho\}}|x|^{\mu_1}u^{(t)}(x)=0,
  \end{equation}
  and
  \begin{equation}\label{equ:scalinginfinity}
    \lim_{t\to 0^+}\esssup_{\{|x|>\rho\}}|x|^{\mu_2}u^{(t)}(x)=\infty
    \qquad\text{and}\qquad
    \lim_{t\to \infty}\essinf_{\{|x|>\rho\}}|x|^{\mu_2}u^{(t)}(x)=0.
  \end{equation}
\end{proposition}
From Proposition \ref{prop:scaling}, we see that neither the upper nor the lower bound in Theorem \ref{thm:twosided} holds, at either end, with a constant depending only on $N$, $s$, $p$, $\mu$, and $\semi{u}$.
\vskip0.2cm

Our main theorem generalized the results in \cite{bib:EMSV2024,bib:Xiang2015} to the fractional case and complement those in \cite{bib:DMPS2016} and \cite{bib:BMS2016,bib:MM2019} to the case when $1<p<\infty$ and $\mu\in(0,\mathcal{C}_{N,s,p})$. Moreover, compared with the main result of \cite[Theorem 1.5]{bib:BM2026} which holds only for the positive radially symmetric nonincreasing weak solution and $p\geq2$, our result holds for every $p\in(1,\infty)$
with $ps<N$ and for every positive weak solution, without any symmetry or monotonicity assumption, and its lower bounds involve no loss in the exponents. In order to prove the main theorem, several fundamental difficulties need to be overcome:
\begin{enumerate}
 \item\label{item:diffode} \emph{No ordinary differential reduction and no local comparison.} For $s=1$, the $p$-Laplacian of a radial function is an ordinary differential expression, which underlies the ODE analysis in \cite{bib:AFP2006}. Local comparison with a barrier in a ball or in an exterior domain uses the ordering on the boundary sphere, together with a suitable condition at infinity in the exterior case; the values of the barrier outside the comparison domain play no role.
     For our energy-space comparisons, due to the nonlocal characterization of $(-\Delta_p)^s$, the values of the barriers on the whole complement of the domain are also needed, not only its boundary.
  \item\label{item:diffgroundstate} \emph{No exact quadratic ground state identity for $p\neq 2$.} The method of \cite{bib:DMPS2016} relies on the exact weighted quadratic identity of \cite{bib:FLS2008}, which is specific to $p=2$. Frank and Seiringer \cite{bib:FS2008} also developed nonlinear ground state representations, which give remainder terms for $p\geqslant 2$, but these do not provide the same exact reduction for $p\neq 2$. We therefore work directly with the nonlinear operator $(-\Delta_p)^s-\mu|x|^{-ps}J_p$ and compare weak sub- and supersolutions, with a Hardy coefficient which is not small.
  \item\label{item:diffcritical} \emph{The critical term.} The nonlinearity $u^{\pst-1}$ has the same scaling as the Hardy term. To treat it as a perturbation one needs a pointwise bound for $u^{\pst-p}$ which is strictly better than $|x|^{-ps}$ at both ends, and this bound has to be derived for an arbitrary positive weak solution, without any symmetry assumption. When cutoff functions are used as test functions, $(-\Delta_p)^s$ produces nonlocal error terms, including the nonlocal tail of $u$ \cite{bib:DKP2014}.
  \item\label{item:diffpurepower} \emph{Construction of sub- and supersolutions.} The main difficulty of this work lies in the construction of sub- and
      supersolutions near the origin and at infinity, which differs substantially from the existing works. In
      the local case \cite{bib:EMSV2024,bib:Xiang2015}, pure powers and two-power profiles serve directly as sub- and supersolutions, which are not applicable in our setting. As we pointed out in \ref{item:diffode}, the barriers has to be defined on all of $\mathbb{R}^{N}$ and belong to $\mathcal{D}^{s,p}(\mathbb{R}^{N})$. Furthermore, the sub- and supersolutions should satisfy the requirements of the comparison principle, therefore
      new barriers are needed both for the lower bounds and for the upper bounds.
\end{enumerate}

To overcome these difficulties, some new ideas and techniques are introduced in our work, especially in the proof of the comparison principle and the constructions of the subsolutions and supersolutions near the origin and at infinity. Let us sketch the main steps of the proof for Theorem \ref{thm:twosided}.
First of all, by using the convexity of the Gagliardo energy along the interpolation of $p$-th powers and a lattice inequality for $\semi{\cdot}^p$, we prove a suitable comparison principle for the weak sub- and supersolutions of the following equation on general open set $\Omega$ of $\mbRN$
\begin{equation*}
  (-\Delta_p)^s f=\Big(\frac{\mu}{|x|^{ps}}+a(x)\Big)f^{p-1}\quad\text{in } \Omega,
\end{equation*}
with nonzero exterior data, where $\mu\in(0,\mathcal{C}_{N,s,p})$ and $a$ is a nonnegative function in $L^{\frac{N}{ps}}(\mbRN)$, see Lemma \ref{lem:comparison}. By applying this result with $\Omega$ being a ball or the complement of a ball, we are able to prove that if the subsolution lies below the supersolution in a neighborhood of $\mathbb{R}^{N}\backslash\Omega$, and strictly below on a set of positive measure outside $\Omega$, then it lies below in all of $\mbRN$, which helps us overcome the difficulties mentioned in \ref{item:diffode} and \ref{item:diffgroundstate}.
\vskip 0.1cm

Next we show that the critical term is a perturbation of lower order, which is difficulty \ref{item:diffcritical}. In the local problem, the preliminary estimates in \cite{bib:EMSV2024} are obtained by testing with cutoff multiples and truncated powers of $u$, followed by Moser iteration. Compared with the
local problem, as we mentioned in \ref{item:diffcritical}, the nonlocal error terms causing by $(-\Delta_p)^s$ and the cutoff test functions increase the difficulty of the proof. To face this fact, we first employ a two-point inequality and the De Giorgi type iteration to establish a local boundedness estimate for positive weak solution $u$ of \eqref{equ:EL} under small critical mass (see Proposition \ref{prop:localbound}), which involves the nonlocal tail of $u$. Then using an explicit integral of the kernel (see \eqref{equ:kernelconstant}) combined with the gap $\mu<\mathcal{C}_{N,s,p}$ in the Hardy inequality, we prove a preliminary decay of the critical mass
\begin{equation*}
  m_0(R)\leqslant CR^{\theta}\ \ \text{for } R\in(0,1],
  \qquad
  m_\infty(R)\leqslant CR^{-\theta}\ \ \text{for } R\in[1,\infty)
\end{equation*}
see Proposition \ref{prop:weighted}. Utilizing these two results, we obtain the pointwise bounds $u(x)\leqslant C|x|^{-r_H+\delta}$ near the origin and $u(x)\leqslant C|x|^{-r_H-\delta}$ at infinity for some $\delta>0$ (Lemma \ref{lem:coarse}).
\vskip 0.1cm

In the third step, we construct the barriers.

\textbf{The subsolutions:} First of all, it is easily seen that the positive weak solution $u$ of \eqref{equ:EL} is a supersolution of the corresponding homogeneous fractional Hardy equation. Unlike the local case that the pure power $|x|^{-\gamma_1}$(respectively, $|x|^{-\gamma_2}$) solves the homogeneous Hardy equation and is compared with the positive weak solution of \eqref{equ:localEL} in a small ball (respectively, outside a ball) to obtain the lower bound near the origin(respectively at infinity, see \cite{bib:EMSV2024,bib:Xiang2015}), in our setting a pure power has to be cut off to guarantee that the function belongs to $\mathcal{D}^{s,p}(\mathbb{R}^{N})$. Furthermore, by careful analysis, we are able to deduce that the single cutoff of the power $|x|^{-\mu_i}$($i=1,2$) are actually the strict supersolutions of the homogeneous Hardy equation, which is not applicable for taking as the subsolutions. To overcome this difficulty, we consider two barriers $\Psi_0(x)=(|x|^{-\mu_1}+1)\chi(x)$ near the origin and $\Psi_\infty(x)=\chi_1(x)|x|^{-\mu_2}+M\zeta(x)$ at infinity, where $\chi$ and $\chi_1$ are cutoff functions and $\zeta$ is a bump function, see Lemmas \ref{lem:lowerbarrierorigin} and \ref{lem:lowerbarrierinfinity}. As we will see, when near the origin, by virtue of \eqref{equ:purepower} and refined estimate of the nonlocal cutoff error, we are able to obtain $(-\Delta_p)^s\Psi_0(x)\leq(-\Delta_p)^s(|x|^{-\mu_{1}})(x)+\mathcal{C}|x|^{-(p-2)\mu_{1}}=\mu|x|^{-ps-\mu_{1}(p-1)}+
\mathcal{C}|x|^{-(p-2)\mu_{1}}<\mu|x|^{-ps}\Psi_0(x)^{p-1}$, this yields that $\Psi_0$ is a strict subsolution. On the other hand, when at infinity, by suitable selections of the constant $M$ and bump function $\zeta$ and using similar arguments mentioned above, we infer that as $|x|\rightarrow\infty$, $|x|^{N+ps}\big((-\Delta_p)^s\Psi_\infty(x)-\mu|x|^{-ps}\Psi_\infty(x)^{p-1}\big)$ converges to a negative limit which is given by the integral of $|y|^{-(p-1)\mu_{2}}-\Psi_\infty^{p-1}$ over a small ball, from which we get $\Psi_\infty$ is a desired strict subsolution of the homogeneous fractional Hardy equation outside a large ball.

\textbf{The supersolutions:} For the local case, since the comparison principle requires only the values of supersolutions within the domain and at its boundary(see \cite{bib:EMSV2024,bib:Xiang2015}), the authors in \cite{bib:EMSV2024,bib:Xiang2015} take the following types of
two-power profiles
\begin{equation*}
  |x|^{-\mu_1}-|x|^{-\vartheta}\ \ (\vartheta<\mu_1)\qquad\text{and}\qquad |x|^{-\mu_2}-|x|^{-\vartheta}\ \ (\vartheta>\mu_2),
\end{equation*}
as supersolutions directly near the origin and at infinity. Owing to the action of the $p$-Laplace and anisotropic $p$-Laplace operators on such two-power profiles can be calculated explicitly, the properties of the supersolutions are easy to obtain. In contrast, the situations in our setting are much more involved, since we need the information of the supersolutions in the entire space due to the nonlocal fractional $p$-Laplacian. To address this issue, we first cut off the above two-power profiles with suitable nonnegative smooth functions to ensure that the cut-off profiles belong to $\mathcal{D}^{s,p}(\mathbb{R}^{N})$. In order to make the supersolutions are defined positive in all of $\mbRN$ and comparable to $u$ outside the domain, we perform the glue method and construct the supersolutions by gluing the cut-off profiles with $u$ suitably. More specifically, near the origin, we construct the supersolutions as $A\chi_\ast(|x|^{-\mu_1}-|x|^{-\vartheta})+(1-\chi_\ast)u$, and $\chi_2u+A(1-\chi_2)(|x|^{-\mu_2}-|x|^{-\vartheta})$ at infinity, where $A\in[1,\infty)$ and $\chi_\ast,\chi_2$ are cutoff functions, see Propositions \ref{prop:upperorigin} and \ref{prop:upperinfinity} for details.
By estimating the nonlocal errors between the supersolutions and cut-off profiles, combined with the pointwise bounds of $u$ (Lemma \ref{lem:coarse}) and complicated and refined computations for $(-\Delta_p)^s$ on these two-power profiles through an inversion representation (see Lemmas \ref{lem:radialpairing} and \ref{lem:twopower}), we are able to show that these constructions for the supersolutions are required and satisfy the conditions of the comparison principle, which allows us to obtain the upper bounds of $u$ near the origin and at infinity and then prove Theorem \ref{thm:twosided}.
\vskip 0.1cm

The paper is organized as follows. Section \ref{sec:preliminaries} collects the properties of the function space, the index function $\Lambda$, weak sub- and supersolutions, and the comparison principle. Section \ref{sec:local} contains the local boundedness estimate and the weak Harnack inequality. In Section \ref{sec:mass} we prove the preliminary decay of the critical mass and the resulting pointwise bounds. The barriers are constructed in Section \ref{sec:barriers}, and in the last section, we prove the lower and upper bounds near the origin and at infinity and finish the proofs of Theorem \ref{thm:twosided} and Proposition \ref{prop:scaling}.
\vskip 0.1cm

\noindent\textbf{Notation.} For any $x\in\mbRN$ and $r\in(0,\infty)$, we write $B(x,r):=\{y\in\mbRN\colon |y-x|<r\}$ and $B_r:=B(0,r)$. For any measurable set $E\subset\mbRN$ we denote by $\mathbf{1}_E$ its characteristic function and by $|E|$ its Lebesgue measure. If $0<|E|<\infty$ and $f$ is integrable on $E$, we write $\fint_E f:=|E|^{-1}\int_E f$. We denote by $\omega_{N-1}:=\mcH^{N-1}(\mbS)$ the surface measure of the unit sphere. For any $f\colon\mbRN\to\mbR$, $z\in\mbRN$, and $R\in(0,\infty)$, the \emph{nonlocal tail} of $f$ is
\begin{equation*}
  \Tail(f;z,R):=\left(R^{ps}\int_{\mbRN\setminus B(z,R)}\frac{|f(y)|^{p-1}}{|y-z|^{N+ps}}\diopy\right)^{\frac{1}{p-1}}.
\end{equation*}
We use $C$ to denote a positive constant which may vary from line to line; its dependence on the parameters is indicated whenever it matters. For $a\in\mbR$ we write $a_+:=\max\{a,0\}$.

\section{Preliminaries}\label{sec:preliminaries}
In this section we collect the facts on the space $\mcD$, the index function $\Lambda$, and weak sub- and supersolutions that will be used later, and we prove a comparison principle.

\subsection{The space \texorpdfstring{$\mcD$}{D\textasciicircum{s,p}(R\textasciicircum{N})}}
Let $0<s<1<p<\infty$, the homogeneous fractional Sobolev space $\mcD$ is defined as
\begin{equation}\label{equ:spaceD}
  \mcD=\left\{f\in L^{\pst}(\mbRN)\colon \semi{f}<\infty\right\},
\end{equation}
endowed with the norm $[\cdot]_{s,p}$, where $[\cdot]_{s,p}$ denotes the Gagliardo seminorm given by
\begin{equation*}
  \semi{f}:=\left(\int_{\mbRN}\int_{\mbRN}\frac{|f(x)-f(y)|^p}{|x-y|^{N+ps}}\diopx\diopy\right)^{\frac{1}{p}}.
\end{equation*}
In light of the fractional Sobolev inequality, there exists a positive constant $S_\ast=S_\ast(N,s,p)$ such that, for any $f\in \mcD$,
\begin{equation}\label{equ:sobolev}
  S_\ast\norm{f}_{L^{\pst}(\mbRN)}^p\leqslant \semi{f}^p.
\end{equation}
Consequently, $\mcD$ is continuously embedded in $L^{\pst}(\mbRN)$, see \cite{bib:BMS2016,bib:BGCV2021,bib:DPV2012} for more details.

Let now $G\colon\mbR\to\mbR$ be a Lipschitz function satisfying $G(0)=0$, then for $f\in\mcD$, it holds $G\circ f\in\mcD$. More generally, let $f_1,f_2\in\mcD$ and let $G\colon\mbR^2\to\mbR$ satisfy $G(0,0)=0$ and
\begin{equation*}
  |G(a_1,a_2)-G(b_1,b_2)|\leqslant L(|a_1-b_1|+|a_2-b_2|).
\end{equation*}
Then $G(f_1,f_2)\in\mcD$.

We say that $u$ is a \emph{positive weak solution} of \eqref{equ:EL} if $u\in\mcD$, $u>0$ almost everywhere in $\mbRN$, and for any $\vp\in\mcD$,
\begin{equation}\label{equ:weakEL}
\int_{\mbR^{N}}\int_{\mbR^{N}}\frac{J_p\left(u(x)-u(y)\right)\left(\vp(x)-\vp(y)\right)}{|x-y|^{N+ps}}\diopx\diopy
  =\mu\int_{\mbRN}\frac{u^{p-1}\vp}{|x|^{ps}}\diopx+\int_{\mbRN}u^{\pst-1}\vp\diopx.
\end{equation}
All three integrals in \eqref{equ:weakEL} converge absolutely, by the H\"older inequality, \eqref{equ:hardyinequality}, and \eqref{equ:sobolev}.

For a positive weak solution $u$ of \eqref{equ:EL}, the nonlocal tail of $u$ is finite. Indeed, for any $z\in\mbRN$ and $R\in(0,\infty)$, the H\"older inequality with the exponents $\frac{\pst}{p-1}$ and $\frac{\pst}{\pst-p+1}$ gives
\begin{equation}\label{equ:tailfinite}
\begin{split}
  \Tail(u;z,R)^{p-1}
  &\leqslant R^{ps}\norm{u}_{L^{\pst}(\mbRN)}^{p-1}\left(\int_{\mbRN\setminus B(z,R)}|y-z|^{-\frac{(N+ps)\pst}{\pst-p+1}}\diopy\right)^{\frac{\pst-p+1}{\pst}}\\
  &=C(N,s,p)\,R^{-(p-1)r_H}\norm{u}_{L^{\pst}(\mbRN)}^{p-1},
\end{split}
\end{equation}
where we used $\frac{(N+ps)\pst}{\pst-p+1}>N$.

The following lemma gives a simple criterion for membership in $\mcD$.
\begin{lemma}\label{lem:powerenergy}
  \begin{enumerate}
    \item\label{item:powerenergyone} Let $a_0<r_H<a_\infty$ and let $F\in C^1(\mbRN\setminus\{0\})$ satisfy
    \begin{equation*}
      |F(x)|+|x||\nabla F(x)|\leqslant
      \begin{cases}
        C_F|x|^{-a_0} &\text{if\ } 0<|x|\leqslant 1,\\
        C_F|x|^{-a_\infty} &\text{if\ } |x|>1
      \end{cases}
    \end{equation*}
    for some constant $C_F$. Then $F\in\mcD$.
    \item\label{item:powerenergytwo} Let $f\in\mcD$ and let $\psi\in C^\infty(\mbRN)$ with $0\leqslant\psi\leqslant 1$, such that $\nabla\psi$ has compact support. Then $\psi f\in\mcD$.
  \end{enumerate}
\end{lemma}
\begin{proof}
  \ref{item:powerenergyone} Let $h(x):=|x|^{-a_0}$ if $0<|x|\leqslant 1$ and $h(x):=|x|^{-a_\infty}$ if $|x|>1$. Since $a_0\pst<N<a_\infty\pst$, we have $F\in L^{\pst}(\mbRN)$. If $|x-y|<\frac{|x|}{2}$, then every point $z$ of the segment joining $x$ and $y$ satisfies $\frac{|x|}{2}\leqslant|z|\leqslant\frac{3|x|}{2}$. For such $z$ we have $h(z)\leqslant Ch(x)$, with $C$ depending only on $a_0$ and $a_\infty$. Indeed, $h$ is continuous and is a power of $|z|$ on each of the sets $\{|z|\leqslant 1\}$ and $\{|z|>1\}$. Moreover, $\frac{|z|}{|x|}\in[\frac{1}{2},\frac{3}{2}]$. The assumption gives $|\nabla F(z)|\leqslant C_F|z|^{-1}h(z)\leqslant C|x|^{-1}h(x)$, and the mean value theorem yields
  \begin{equation*}
    |F(x)-F(y)|\leqslant C|x-y||x|^{-1}h(x)\quad\text{if\ } |x-y|<\frac{|x|}{2}.
  \end{equation*}
  Therefore
  \begin{equation*}
    \int_{\{|x-y|<\frac{|x|}{2}\}}\frac{|F(x)-F(y)|^p}{|x-y|^{N+ps}}\diopy
    \leqslant C|x|^{-p}h(x)^p\int_{\{|x-y|<\frac{|x|}{2}\}}|x-y|^{p-N-ps}\diopy
    =C h(x)^p|x|^{-ps},
  \end{equation*}
  and
  \begin{equation*}
    \int_{\{|x-y|\geqslant\frac{|x|}{2}\}}\frac{|F(x)|^p}{|x-y|^{N+ps}}\diopy\leqslant C|F(x)|^p|x|^{-ps}\leqslant Ch(x)^p|x|^{-ps}.
  \end{equation*}
  Moreover, if $|x-y|\geqslant\frac{|x|}{2}$, then $|x-y|\geqslant\frac{|y|}{3}$, and thus
  \begin{equation*}
    \int_{\mbRN}\int_{\{x\colon|x-y|\geqslant\frac{|x|}{2}\}}\frac{|F(y)|^p}{|x-y|^{N+ps}}\diopx\diopy
    \leqslant C\int_{\mbRN}h(y)^p|y|^{-ps}\diopy.
  \end{equation*}
  Since $|F(x)-F(y)|^p\leqslant 2^{p-1}(|F(x)|^p+|F(y)|^p)$, we conclude that $\semi{F}^p\leqslant C\int_{\mbRN}h(x)^p|x|^{-ps}\diopx$. In polar coordinates,
  \begin{equation*}
    \int_{\mbRN}h(x)^p|x|^{-ps}\diopx=\omega_{N-1}\left(\int_0^1t^{N-1-p(a_0+s)}\diopt+\int_1^\infty t^{N-1-p(a_\infty+s)}\diopt\right),
  \end{equation*}
  and both integrals are finite because $p(a_0+s)<N<p(a_\infty+s)$, which follows from $a_0<r_H<a_\infty$ and $p(r_H+s)=N$. By \eqref{equ:spaceD}, $F\in\mcD$.

  \ref{item:powerenergytwo} Clearly $\psi f\in L^{\pst}(\mbRN)$. Since $|\psi(x)f(x)-\psi(y)f(y)|\leqslant |f(x)-f(y)|+|f(y)||\psi(x)-\psi(y)|$, it suffices to show that
  \begin{equation*}
    g(y):=\int_{\mbRN}\frac{|\psi(x)-\psi(y)|^p}{|x-y|^{N+ps}}\diopx
  \end{equation*}
  belongs to $L^{\frac{N}{ps}}(\mbRN)$, because then $\int_{\mbRN}|f|^pg\leqslant\norm{f}_{L^{\pst}(\mbRN)}^p\norm{g}_{L^{\frac{N}{ps}}(\mbRN)}$ by the H\"older inequality. Let $\supp\nabla\psi\subset B_{\rho_\psi}$ for some $\rho_\psi\in(0,\infty)$. Then $g$ is bounded, since $|\psi(x)-\psi(y)|\leqslant\min\{\norm{\nabla\psi}_{L^\infty}|x-y|,1\}$. If $|y|\geqslant 2\rho_\psi$ and $\psi(x)\neq\psi(y)$, then, because $\psi$ is constant on the connected set $\mbRN\setminus B_{\rho_\psi}$, we have $|x|<\rho_\psi$, and hence $|x-y|\geqslant\frac{|y|}{2}$; thus $g(y)\leqslant C|B_{\rho_\psi}||y|^{-N-ps}$. Therefore $g\in L^{\frac{N}{ps}}(\mbRN)$, which completes the proof of Lemma \ref{lem:powerenergy}.
\end{proof}

\subsection{The index function}

Fix $e\in\mbS$. The function $K$ in \eqref{equ:Lambdaformula} does not depend on the choice of $e$, and, by the rotation invariance of $\mcH^{N-1}$, $K(\sigma)=\int_{\mbS}|e-\sigma\varpi|^{-N-ps}\diop\mcH^{N-1}(\varpi)$. Moreover, we have the following crucial result.
\begin{proposition}\label{prop:Kbound}
Let $K$ be given in \eqref{equ:Lambdaformula}, then there exists $C=C(N,s,p)$ such that
\begin{equation}\label{equ:Kbound}
  K(\sigma)\leqslant C(\sigma-1)^{-1-ps}\quad\text{for any\ } \sigma\in(1,2], \ \ \
  K(\sigma)\leqslant C\sigma^{-N-ps}\quad\text{for any\ } \sigma\in(2,\infty).
\end{equation}
\end{proposition}
\begin{proof}
For $\varpi\in\mbS$, we have
\begin{equation*}
  |\sigma e-\varpi|^2=\sigma^2-2\sigma e\cdot\varpi+1=(\sigma-1)^2+\sigma|e-\varpi|^2,
\end{equation*}
therefore $|\sigma e-\varpi|\geqslant\sigma-1$. If $\sigma\in(2,\infty)$, then $\sigma-1>\frac{\sigma}{2}$, this yields $K(\sigma)\leqslant 2^{N+ps}\omega_{N-1}\sigma^{-N-ps}$,

Let now $\sigma\in(1,2]$. Then $|\sigma e-\varpi|^2\geqslant(\sigma-1)^2+|e-\varpi|^2$. Since $|e-\varpi|^2=2-2e\cdot\varpi$, this lower bound depends on $\varpi$ only through $e\cdot\varpi$. We write $e\cdot\varpi=\cos\beta$ with $\beta\in[0,\pi]$. For integrands that depend on $\varpi$ only through $e\cdot\varpi$, the surface measure is $\diop\mcH^{N-1}(\varpi)=\omega_{N-2}\sin^{N-2}\beta\diop\beta$. Here $\omega_{N-2}$ denotes the surface measure of the unit sphere of $\mbR^{N-1}$, so that $\omega_0=2$. Moreover, $|e-\varpi|=2\sin\frac{\beta}{2}\geqslant\frac{2\beta}{\uppi}$ and $\sin\beta\leqslant\beta$ for $\beta\in[0,\uppi]$. Therefore
\begin{align*}
  K(\sigma)&\leqslant\omega_{N-2}\int_0^{\uppi}\left((\sigma-1)^2+\frac{4\beta^2}{\uppi^2}\right)^{-\frac{N+ps}{2}}\beta^{N-2}\diop\beta\\
  &=\omega_{N-2}\left(\frac{\uppi}{2}\right)^{N-1}(\sigma-1)^{-1-ps}\int_0^{\frac{2}{\sigma-1}}\frac{t^{N-2}}{(1+t^2)^{\frac{N+ps}{2}}}\diopt,
\end{align*}
where we substituted $\beta=\frac{\uppi}{2}(\sigma-1)t$. The last integral is at most $\int_0^\infty t^{N-2}(1+t^2)^{-\frac{N+ps}{2}}\diopt$ and finite, because $N-2-(N+ps)=-2-ps<-1$. This proves \eqref{equ:Kbound}.
\end{proof}

The following lemma expresses the fractional $p$-Laplacian of a radial function at a point of $\mbS$ by pairing each point $y$ with $|y|>1$ with its inversion $\frac{y}{|y|^2}$.
\begin{lemma}\label{lem:radialpairing}
  Let $f\colon(0,\infty)\to\mbR$ be continuous, twice continuously differentiable in a neighborhood of $1$, and such that
  \begin{equation}\label{equ:radialintegrability}
    \int_0^1|f(t)|^{p-1}t^{N-1}\diopt+\int_1^\infty|f(t)|^{p-1}t^{-1-ps}\diopt<\infty.
  \end{equation}
  If $p\in(1,2)$, assume moreover that $f'(1)\neq 0$. Let $F(x):=f(|x|)$ for $x\neq 0$. Then, for any $e\in\mbS$, the limit in \eqref{equ:fractionalplaplacian} defining $(-\Delta_p)^sF(e)$ exists and
  \begin{equation}\label{equ:radialpairing}
    (-\Delta_p)^sF(e)=2\int_1^\infty K(\sigma)\left[\sigma^{N-1}J_p\left(f(1)-f(\sigma)\right)+\sigma^{ps-1}J_p\left(f(1)-f\left(\frac{1}{\sigma}\right)\right)\right]\diopsigma,
  \end{equation}
  where the integral on the right-hand side converges absolutely.
\end{lemma}
\begin{proof}
  For $\ve\in(0,\frac{1}{4})$, let
  \begin{equation*}
    I_\ve:=2\int_{\{|e-y|>\ve\}}\frac{J_p\left(F(e)-F(y)\right)}{|e-y|^{N+ps}}\diopy,
  \end{equation*}
  which converges absolutely by \eqref{equ:radialintegrability}. We split the domain of integration into $\{|y|>1\}$ and $\{|y|<1\}$. In the second part we change variables $y=\frac{y^\ast}{|y^\ast|^2}$ with $|y^\ast|>1$. Then $\diopy=|y^\ast|^{-2N}\diop y^\ast$ and $|e-y|=\frac{|e-y^\ast|}{|y^\ast|}$, so that
  \begin{equation*}
    \int_{\{|y|<1,\,|e-y|>\ve\}}\frac{J_p\left(F(e)-F(y)\right)}{|e-y|^{N+ps}}\diopy
    =\int_{\{|y^\ast|>1,\,|e-y^\ast|>\ve|y^\ast|\}}\frac{|y^\ast|^{ps-N}J_p\left(f(1)-f(|y^\ast|^{-1})\right)}{|e-y^\ast|^{N+ps}}\diop y^\ast.
  \end{equation*}
  The two domains of integration $\{|y|>1,\,|e-y|>\ve\}$ and $\{|y|>1,\,|e-y|>\ve|y|\}$ differ by the set
  \begin{equation*}
    D_\ve:=\left\{y\in\mbRN\colon|y|>1,\ \ve<|e-y|\leqslant\ve|y|\right\}.
  \end{equation*}
  On $D_\ve$ we have $|y|\leqslant|e|+|e-y|\leqslant 1+\ve|y|$, hence $|y|\leqslant\frac{1}{1-\ve}\leqslant 1+2\ve$, and therefore $D_\ve\subset\{\ve<|e-y|\leqslant\ve(1+2\ve)\}$. This spherical shell has measure at most $C\ve^{N-1}\cdot 2\ve^2=C\ve^{N+1}$. Moreover, on $D_\ve$ we have $0<|y|-1\leqslant|e-y|\leqslant 2\ve$, and, since $f$ is Lipschitz near $1$, $|f(1)-f(|y|^{-1})|\leqslant C(1-|y|^{-1})\leqslant C(|y|-1)\leqslant C\ve$. Therefore the integral over $D_\ve$ of the last integrand is bounded by $C\ve^{N+1}\ve^{p-1}\ve^{-N-ps}=C\ve^{p-ps}$, which tends to $0$ as $\ve\to 0^+$. Consequently,
  \begin{equation*}
    I_\ve=2\int_{\{|y|>1,\,|e-y|>\ve\}}P(y)\diopy+o(1)\quad\text{as } \ve\to 0^+,
  \end{equation*}
  where
  \begin{equation*}
    P(y):=\frac{J_p\left(f(1)-f(|y|)\right)+|y|^{ps-N}J_p\left(f(1)-f(|y|^{-1})\right)}{|e-y|^{N+ps}}.
  \end{equation*}

  We now show that $P\in L^1(\{|y|>1\})$. Let $t:=|y|-1$. By the Taylor formula, $a:=f(1)-f(1+t)=-f'(1)t+O(t^2)$ and $b:=f(1)-f(\frac{1}{1+t})=f'(1)t+O(t^2)$ as $t\to 0^+$, so that $|a+b|\leqslant Ct^2$ and $|a|+|b|\leqslant Ct$. If $p\geqslant 2$, then
  \begin{equation*}
    |J_p(a)+J_p(b)|=|J_p(a)-J_p(-b)|\leqslant (p-1)\left(|a|+|b|\right)^{p-2}|a+b|\leqslant Ct^p.
  \end{equation*}
  If $p\in(1,2)$, then, since $f'(1)\neq 0$, there exists $t_0\in(0,\frac{1}{2})$ such that, for $t\in(0,t_0)$, the numbers $a$ and $-b$ have the same sign and $|a|,|b|\geqslant\frac{|f'(1)|t}{2}$; by the mean value theorem, $|J_p(a)-J_p(-b)|\leqslant(p-1)(\frac{|f'(1)|t}{2})^{p-2}|a+b|\leqslant Ct^p$. In both cases, since $||y|^{ps-N}-1|\leqslant Ct$ and $|J_p(b)|\leqslant Ct^{p-1}$, we obtain, for some $t_0\in(0,\frac{1}{2})$,
  \begin{equation*}
    |P(y)|\leqslant Ct^{p}|e-y|^{-N-ps}\leqslant C|e-y|^{p-N-ps}\quad\text{if\ } 1<|y|<1+t_0,
  \end{equation*}
  where we used $t=|y|-|e|\leqslant|e-y|$. The right-hand side is integrable near $e$ because $p-ps>0$. On $\{|y|\geqslant 1+t_0\}$ we have $|e-y|\geqslant t_0$ and $|e-y|\geqslant\frac{|y|}{2}$ if $|y|\geqslant 2$, whence
  \begin{equation*}
    |P(y)|\leqslant C\left[1+|f(|y|)|^{p-1}+|y|^{ps-N}\left(1+|f(|y|^{-1})|^{p-1}\right)\right]\left(1+|y|\right)^{-N-ps}.
  \end{equation*}
  In polar coordinates,
  \begin{gather*}
    \int_{\{|y|\geqslant 2\}}\frac{|f(|y|)|^{p-1}}{|y|^{N+ps}}\diopy=\omega_{N-1}\int_2^\infty\frac{|f(t)|^{p-1}}{t^{1+ps}}\diopt,\\
    \int_{\{|y|\geqslant 2\}}\frac{|f(|y|^{-1})|^{p-1}}{|y|^{2N}}\diopy=\omega_{N-1}\int_0^{\frac{1}{2}}|f(t)|^{p-1}t^{N-1}\diopt,
  \end{gather*}
  which are finite by \eqref{equ:radialintegrability}. Hence $P\in L^1(\{|y|>1\})$, and the dominated convergence theorem yields $\lim_{\ve\to 0^+}I_\ve=2\int_{\{|y|>1\}}P(y)\diopy$. Writing $y=\sigma\varpi$ with $\sigma>1$ and $\varpi\in\mbS$, we have $\diopy=\sigma^{N-1}\diopsigma\diop\mcH^{N-1}(\varpi)$ and $\int_{\mbS}|e-\sigma\varpi|^{-N-ps}\diop\mcH^{N-1}(\varpi)=K(\sigma)$, which gives \eqref{equ:radialpairing}. This finishes the proof of Lemma \ref{lem:radialpairing}.
\end{proof}

\begin{lemma}\label{lem:indexfunction}
  Let $r\in(0,\frac{N}{p-1})$. Then \eqref{equ:purepower} holds for any $x\neq 0$, with $\Lambda(r)$ given by \eqref{equ:Lambdaformula}, and the integral in \eqref{equ:Lambdaformula} converges absolutely. Moreover, the following hold.
  \begin{enumerate}
    \item\label{item:indexsign} $\Lambda(r)>0$ if $r\in(0,\frac{N-ps}{p-1})$, $\Lambda(\frac{N-ps}{p-1})=0$, and $\Lambda(r)<0$ if $r\in(\frac{N-ps}{p-1},\frac{N}{p-1})$.
    \item\label{item:indexmax} $\Lambda(r_H)=\mathcal{C}_{N,s,p}$.
    \item\label{item:indexmonotone} $\Lambda$ is strictly increasing on $(0,r_H]$, strictly decreasing on $[r_H,\frac{N-ps}{p-1}]$, and $\lim_{r\to 0^+}\Lambda(r)=0$.
  \end{enumerate}
  In particular, for any $\mu\in(0,\mathcal{C}_{N,s,p})$, the equation $\Lambda(r)=\mu$ has exactly two roots $\mu_1$ and $\mu_2$ in $(0,\frac{N}{p-1})$, and they satisfy \eqref{equ:tworoots}.
\end{lemma}
\begin{proof}
  Let $f(t):=t^{-r}$. Then $f'(1)=-r\neq 0$, and \eqref{equ:radialintegrability} holds because $r(p-1)<N$. By Lemma \ref{lem:radialpairing}, $(-\Delta_p)^s(|\cdot|^{-r})(e)$ is given by \eqref{equ:radialpairing}. Since $f(1)-f(\sigma)=1-\sigma^{-r}$ and $f(1)-f(\frac{1}{\sigma})=-(\sigma^r-1)$ for $\sigma>1$, and since $\sigma^{N-1}(1-\sigma^{-r})^{p-1}=\sigma^{N-1-r(p-1)}(\sigma^r-1)^{p-1}$, the right-hand side of \eqref{equ:radialpairing} equals the right-hand side of the first identity in \eqref{equ:Lambdaformula}. For $x\neq 0$, the change of variables $y=|x|y'$ in \eqref{equ:fractionalplaplacian}, together with the rotation invariance, gives $(-\Delta_p)^s(|\cdot|^{-r})(x)=|x|^{-ps-r(p-1)}(-\Delta_p)^s(|\cdot|^{-r})(\frac{x}{|x|})$, which proves \eqref{equ:purepower}.

  Assertion \ref{item:indexsign} follows directly from \eqref{equ:Lambdaformula}, since, for $\sigma>1$, the factor $\sigma^{N-1-r(p-1)}-\sigma^{ps-1}$ is positive, zero, or negative according as $r(p-1)$ is smaller than, equal to, or larger than $N-ps$.

  To prove \ref{item:indexmax}, note that $N-1-r_H(p-1)=r_H+ps-1$, so that $\Lambda(r_H)=2\int_1^\infty K(\sigma)\sigma^{ps-1}(\sigma^{r_H}-1)^p\diopsigma$. The change of variables $\sigma=\frac{1}{t}$, together with $K(\frac{1}{t})=t^{N+ps}\int_{\mbS}|e-t\varpi|^{-N-ps}\diop\mcH^{N-1}(\varpi)$ and $p r_H=N-ps$, yields
  \begin{equation*}
    \Lambda(r_H)=2\int_0^1 t^{ps-1}\left(1-t^{r_H}\right)^p\left(\int_{\mbS}\frac{\diop\mcH^{N-1}(\varpi)}{|e-t\varpi|^{N+ps}}\right)\diopt,
  \end{equation*}
  which is the value of the sharp constant in \eqref{equ:hardyinequality} given in \cite[Theorem 1.1]{bib:FS2008}.

  The monotonicity in \ref{item:indexmonotone} is contained in \cite[Lemma 3.1]{bib:AB2017} and its proof, where the function $\frac{1}{2}\Lambda$ is considered; the limit in \ref{item:indexmonotone} is proved below.

  We next show that $\Lambda$ is continuous on $(0,\frac{N}{p-1})$ and that $\lim_{r\to 0^+}\Lambda(r)=0$. Let $r_b\in(0,\frac{N}{p-1})$ and $r\in(0,r_b]$. For $\sigma\in(1,2]$, the mean value theorem gives $(\sigma^r-1)^{p-1}\leqslant(r_b2^{r_b}(\sigma-1))^{p-1}$ and $|\sigma^{N-1-r(p-1)}-\sigma^{ps-1}|\leqslant C(\sigma-1)$, where $C$ depends only on $N$, $s$, $p$, and $r_b$. For $\sigma\in(2,\infty)$ we have $(\sigma^r-1)^{p-1}|\sigma^{N-1-r(p-1)}-\sigma^{ps-1}|\leqslant\sigma^{N-1}+\sigma^{r_b(p-1)+ps-1}$. Hence, by \eqref{equ:Kbound}, the integrand in \eqref{equ:Lambdaformula} is bounded, uniformly in $r\in(0,r_b]$, by $C(\sigma-1)^{p-1-ps}$ for $\sigma\in(1,2]$ and by $C(\sigma^{-1-ps}+\sigma^{r_b(p-1)-N-1})$ for $\sigma\in(2,\infty)$. Both bounds are integrable, because $p-ps>0$, $ps>0$, and $r_b(p-1)<N$. For each $\sigma>1$ the integrand depends continuously on $r$ and tends to $0$ as $r\to 0^+$. Therefore the dominated convergence theorem shows that $\Lambda$ is continuous on $(0,r_b]$ and that $\Lambda(r)\to 0$ as $r\to 0^+$. Since $r_b\in(0,\frac{N}{p-1})$ is arbitrary, the claim follows.

  Finally, by \ref{item:indexsign}, $\Lambda\leqslant 0$ on $[\frac{N-ps}{p-1},\frac{N}{p-1})$, and, by \ref{item:indexmax}, \ref{item:indexmonotone}, the continuity of $\Lambda$, and the intermediate value theorem, $\Lambda$ takes each value in $(0,\mathcal{C}_{N,s,p})$ exactly once on $(0,r_H)$ and exactly once on $(r_H,\frac{N-ps}{p-1})$. This completes the proof of Lemma \ref{lem:indexfunction}.
\end{proof}

Since $\mu_1<r_H<\mu_2$ and $r_H\pst=N$, we have
\begin{equation}\label{equ:exponentfacts}
  \mu_1\pst<N<\mu_2\pst,\qquad p(\mu_1+s)<N<p(\mu_2+s),\qquad \mu_2(p-1)<N-ps,
\end{equation}
where the last inequality follows from \eqref{equ:tworoots}.

\subsection{Weak sub- and supersolutions}
We call $V$ an \emph{admissible potential} if
\begin{equation}\label{equ:admissiblepotential}
  V(x)=\frac{\mu}{|x|^{ps}}+a(x),\qquad x\in\mbRN\setminus\{0\},
\end{equation}
for some nonnegative function $a\in L^{\frac{N}{ps}}(\mbRN)$. For such $V$ and any $f,\vp\in\mcD$, by the H\"older inequality, \eqref{equ:hardyinequality}, and \eqref{equ:sobolev},
\begin{equation}\label{equ:potentialbound}
  \int_{\mbRN}V|f|^{p-1}|\vp|
  \leqslant\frac{\mu}{\mathcal{C}_{N,s,p}}\semi{f}^{p-1}\semi{\vp}+\norm{a}_{L^{\frac{N}{ps}}(\mbRN)}\norm{f}_{L^{\pst}(\mbRN)}^{p-1}\norm{\vp}_{L^{\pst}(\mbRN)}
  \leqslant C\semi{f}^{p-1}\semi{\vp}.
\end{equation}
For $f,\vp\in\mcD$ we write
\begin{equation*}
  \langle(-\Delta_p)^sf,\vp\rangle:=\int_{\mbRN}\int_{\mbRN}\frac{J_p\left(f(x)-f(y)\right)\left(\vp(x)-\vp(y)\right)}{|x-y|^{N+ps}}\diopx\diopy,
\end{equation*}
so that $|\langle(-\Delta_p)^sf,\vp\rangle|\leqslant\semi{f}^{p-1}\semi{\vp}$. For an open set $\Omega\subset\mbRN$ and $\delta\in(0,\infty)$ let $\Omega_\delta:=\{x\in\Omega\colon\dist(x,\mbRN\setminus\Omega)>\delta\}$ (with $\Omega_\delta:=\mbRN$ if $\Omega=\mbRN$), and let
\begin{equation*}
  X_+(\Omega):=\left\{\vp\in\mcD\colon \vp\geqslant 0\text{ a.e. in } \mbRN,\ \vp=0 \text{ a.e. in } \mbRN\setminus\Omega_\delta\text{ for some } \delta\in(0,\infty)\right\}.
\end{equation*}

\begin{definition}\label{def:subsuper}
  Let $\Omega\subset\mbRN$ be open, let $V$ be an admissible potential, and let $f\in\mcD$ be nonnegative. We say that $f$ is a \emph{weak supersolution} (respectively, \emph{weak subsolution}) of $(-\Delta_p)^sf=Vf^{p-1}$ in $\Omega$ if, for any $\vp\in X_+(\Omega)$,
  \begin{equation*}
    \langle(-\Delta_p)^sf,\vp\rangle\geqslant\int_{\Omega}Vf^{p-1}\vp
    \qquad\left(\text{respectively, } \langle(-\Delta_p)^sf,\vp\rangle\leqslant\int_{\Omega}Vf^{p-1}\vp\right).
  \end{equation*}
\end{definition}

Let $u$ be a positive weak solution of \eqref{equ:EL}. Since $(\pst-p)\frac{N}{ps}=\pst$, the function $u^{\pst-p}$ belongs to $L^{\frac{N}{ps}}(\mbRN)$, and hence
\begin{equation}\label{equ:Vu}
  V_u(x):=\frac{\mu}{|x|^{ps}}+u(x)^{\pst-p}
\end{equation}
is an admissible potential. By \eqref{equ:weakEL}, $\langle(-\Delta_p)^su,\vp\rangle=\int_{\mbRN}V_uu^{p-1}\vp$ for any $\vp\in\mcD$; in particular, for any open set $\Omega$, $u$ is both a weak supersolution and a weak subsolution of $(-\Delta_p)^sf=V_uf^{p-1}$ in $\Omega$. Moreover, since $u^{\pst-1}\geqslant 0$, $u$ is a weak supersolution of $(-\Delta_p)^sf=\mu|x|^{-ps}f^{p-1}$ in any open set, and $\langle(-\Delta_p)^su,\vp\rangle\geqslant 0$ for any nonnegative $\vp\in\mcD$.

\begin{lemma}\label{lem:approximation}
  Let $\Omega\subset\mbRN$ be open and $\vp\in X_+(\Omega)$. Then there exist nonnegative functions $\vp_k\in C_\mathrm{c}^\infty(\Omega\setminus\{0\})$, $k\in\mbN$, such that $\semi{\vp_k-\vp}\to 0$ as $k\to\infty$.
\end{lemma}
\begin{proof}
  Let $\delta\in(0,\infty)$ be such that $\vp=0$ a.e. in $\mbRN\setminus\Omega_\delta$. Fix $\chi\in C_\mathrm{c}^\infty(B_2)$ with $0\leqslant\chi\leqslant 1$ and $\chi=1$ on $B_1$, and let $\chi_L:=\chi(\frac{\cdot}{L})$ for $L\in[1,\infty)$.

  \emph{Step 1.} We show that $\semi{\vp-\chi_L\vp}\to 0$ as $L\to\infty$. We have
  \begin{equation*}
    \semi{(1-\chi_L)\vp}^p\leqslant 2^{p-1}\int_{\mbRN}\int_{\mbRN}\frac{(1-\chi_L(x))^p|\vp(x)-\vp(y)|^p}{|x-y|^{N+ps}}\diopx\diopy+2^{p-1}\int_{\mbRN}|\vp(y)|^pg_L(y)\diopy,
  \end{equation*}
  where $g_L(y):=\int_{\mbRN}|\chi_L(x)-\chi_L(y)|^p|x-y|^{-N-ps}\diopx$. The first term tends to $0$ by the dominated convergence theorem. As in the proof of Lemma \ref{lem:powerenergy}\ref{item:powerenergytwo}, $g_1\in L^\infty(\mbRN)\cap L^{\frac{N}{ps}}(\mbRN)$, and $g_L(y)=L^{-ps}g_1(\frac{y}{L})$. Hence $\norm{g_L}_{L^{\frac{N}{ps}}(\mbRN)}=\norm{g_1}_{L^{\frac{N}{ps}}(\mbRN)}$ and $\norm{g_L}_{L^\infty(\mbRN)}=L^{-ps}\norm{g_1}_{L^\infty(\mbRN)}$. For any $M\in(0,\infty)$, by the H\"older inequality,
  \begin{equation*}
    \int_{\mbRN}|\vp|^pg_L\leqslant L^{-ps}\norm{g_1}_{L^\infty(\mbRN)}\int_{B_M}|\vp|^p+\norm{g_1}_{L^{\frac{N}{ps}}(\mbRN)}\norm{\vp}_{L^{\pst}(\mbRN\setminus B_M)}^p.
  \end{equation*}
  Letting first $L\to\infty$ and then $M\to\infty$, we find that the second term also tends to $0$.

  \emph{Step 2.} Let $\phi:=\chi_L\vp$ for a fixed $L$, and let $\phi_\ve:=\phi\ast\eta_\ve$, where $\eta_\ve$ is a standard nonnegative mollifier supported in $B_\ve$ with $\ve\in(0,\frac{\delta}{2})$. Then $\phi_\ve\in C_\mathrm{c}^\infty(\Omega)$ is nonnegative, and $\semi{\phi_\ve-\phi}\to 0$ as $\ve\to 0^+$. Indeed, set $\Phi(x,y):=(\phi(x)-\phi(y))|x-y|^{-\frac{N+ps}{p}}$. Then $\semi{\phi}=\norm{\Phi}_{L^p(\mbR^{2N})}$, and
\begin{equation*}
  \phi_\ve(x)-\phi_\ve(y)=\int_{B_\ve}\eta_\ve(h)(\phi(x-h)-\phi(y-h))\diop h.
\end{equation*}
By the Minkowski inequality,
\begin{equation*}
  \semi{\phi_\ve-\phi}\leqslant\sup_{|h|<\ve}\norm{\Phi(\cdot-h,\cdot-h)-\Phi}_{L^p(\mbR^{2N})}.
\end{equation*}
The right-hand side tends to $0$ by the continuity of translations in $L^p(\mbR^{2N})$.

  \emph{Step 3.} If $0\notin\Omega$, Steps 1 and 2 complete the proof. If $0\in\Omega$, let $\psi\in C_\mathrm{c}^\infty(\Omega)$ be nonnegative, and let $\psi_\ve:=(1-\chi(\frac{\cdot}{\ve}))\psi$ for $\ve\in(0,\infty)$. Then $\psi_\ve\in C_\mathrm{c}^\infty(\Omega\setminus\{0\})$ is nonnegative, and
  \begin{equation*}
    \semi{\psi-\psi_\ve}^p=\semi{\chi(\frac{\cdot}{\ve})\psi}^p
    \leqslant 2^{p-1}\norm{\psi}_{L^\infty(\mbRN)}^p\semi{\chi(\frac{\cdot}{\ve})}^p+2^{p-1}\int_{B_{2\ve}}\int_{\mbRN}\frac{|\psi(x)-\psi(y)|^p}{|x-y|^{N+ps}}\diopx\diopy,
  \end{equation*}
  where $\semi{\chi(\frac{\cdot}{\ve})}^p=\ve^{N-ps}\semi{\chi}^p$ and the inner integral is bounded uniformly in $y$. Since $ps<N$, both terms tend to $0$ as $\ve\to 0^+$. To conclude, fix $k\in\mbN$. By Step 1, choose $L$ with $\semi{\vp-\chi_L\vp}<\frac{1}{3k}$. By Step 2, choose $\ve_0>0$ such that the nonnegative function $\psi:=\phi_{\ve_0}\in C_\mathrm{c}^\infty(\Omega)$ satisfies $\semi{\chi_L\vp-\psi}<\frac{1}{3k}$. If $0\in\Omega$, choose $\ve$ with $\semi{\psi-\psi_\ve}<\frac{1}{3k}$ and set $\vp_k:=\psi_\ve$. Otherwise, set $\vp_k:=\psi$. In the latter case, $\psi$ equals $0$ near $0$ because $\supp\psi\subset\Omega$. Then $\vp_k\in C_\mathrm{c}^\infty(\Omega\setminus\{0\})$ is nonnegative and $\semi{\vp-\vp_k}<\frac{1}{k}$. This completes the proof of Lemma \ref{lem:approximation}.
\end{proof}

By Lemma \ref{lem:approximation}, \eqref{equ:potentialbound}, and the bound $|\langle(-\Delta_p)^sf,\vp\rangle|\leqslant\semi{f}^{p-1}\semi{\vp}$, in Definition \ref{def:subsuper} it suffices to take the test functions $\vp$ in the set of nonnegative functions in $C_\mathrm{c}^\infty(\Omega\setminus\{0\})$. In the same way, applying Lemma \ref{lem:approximation} with $\Omega=\mbRN$ to $\vp_+$ and $\vp_-$, \eqref{equ:weakEL} holds for any $\vp\in\mcD$ as soon as it holds for any $\vp\in C_\mathrm{c}^\infty(\mbRN\setminus\{0\})$.

The next lemma converts pointwise inequalities for functions which are smooth away from the origin into weak inequalities.
\begin{lemma}\label{lem:pointwiseweak}
  Let $f\in\mcD$ and let $D\subset\mbRN\setminus\{0\}$ be open. Assume that $f\in C^2(D)$ and, if $p\in(1,2)$, that $\nabla f(x)\neq 0$ for any $x\in D$. Then the limit in \eqref{equ:fractionalplaplacian} defining $(-\Delta_p)^sf(x)$ exists for any $x\in D$, the function $(-\Delta_p)^sf$ is continuous in $D$, and, for any $\vp\in C_\mathrm{c}^\infty(D)$,
  \begin{equation}\label{equ:pointwiseweak}
    \langle(-\Delta_p)^sf,\vp\rangle=\int_{D}(-\Delta_p)^sf(x)\,\vp(x)\diopx.
  \end{equation}
  Consequently, let $\Omega\subset\mbRN$ be open, let $V$ be an admissible potential, let $f\in\mcD$ be nonnegative, and assume that the above assumptions hold with $D=\Omega\setminus\{0\}$. If $(-\Delta_p)^sf(x)\geqslant V(x)f(x)^{p-1}$ (respectively, $(-\Delta_p)^sf(x)\leqslant V(x)f(x)^{p-1}$) for almost every $x\in\Omega$, then $f$ is a weak supersolution (respectively, a weak subsolution) of $(-\Delta_p)^sf=Vf^{p-1}$ in $\Omega$.
\end{lemma}
\begin{proof}
  Let $E\subset D$ be compact, let $d\in(0,\frac{1}{2}\dist(E,\mbRN\setminus D))$, and let $E_{2d}:=\{x\colon\dist(x,E)<2d\}$. Then $\ol{E_{2d}}$ is a compact subset of $D$, so that $M:=\max_{\ol{E_{2d}}}\left(|\nabla f|+|\nabla^2f|\right)<\infty$ and, if $p\in(1,2)$, $m:=\min_{\ol{E_{2d}}}|\nabla f|>0$. For $x\in E$ and $\ve\in(0,d)$ let
  \begin{equation*}
    T_\ve(x):=2\int_{\{|h|>\ve\}}\frac{J_p\left(f(x)-f(x+h)\right)}{|h|^{N+ps}}\diop h.
  \end{equation*}
  Since $|J_p(f(x)-f(x+h))|\leqslant 2^{p}(|f(x)|^{p-1}+|f(x+h)|^{p-1})$ and $|f|^{p-1}\in L^{\frac{\pst}{p-1}}(\mbRN)$, $T_\ve$ is finite, and it is continuous on $E$: for $x,x'\in E$,
  \begin{equation*}
    |T_\ve(x)-T_\ve(x')|\leqslant 2\int_{\{|h|>\ve\}}\left|J_p\left(f(x)-f(x+h)\right)-J_p\left(f(x')-f(x'+h)\right)\right||h|^{-N-ps}\diop h,
  \end{equation*}
  which tends to $0$ as $x'\to x$ by the continuity of $f$ on $E$, the continuity of translations in $L^{\pst}(\mbRN)$, the H\"older inequality, and the bounds $|J_p(a)-J_p(b)|\leqslant 2^{2-p}|a-b|^{p-1}$ for $p\in(1,2)$ and $|J_p(a)-J_p(b)|\leqslant(p-1)(|a|+|b|)^{p-2}|a-b|$ for $p\geqslant 2$. For $0<\ve<\ve'<d$, by the symmetry $h\mapsto-h$,
  \begin{equation*}
    T_\ve(x)-T_{\ve'}(x)=\int_{\{\ve<|h|\leqslant\ve'\}}\frac{J_p\left(f(x)-f(x+h)\right)+J_p\left(f(x)-f(x-h)\right)}{|h|^{N+ps}}\diop h.
  \end{equation*}
  Let $A:=\nabla f(x)\cdot h$, $a_\pm:=f(x\pm h)-f(x)$. By the Taylor formula, $|a_\pm\mp A|\leqslant\frac{M}{2}|h|^2$, so that $|a_++a_-|\leqslant M|h|^2$ and the numerator above equals $-(J_p(a_+)-J_p(-a_-))$. If $p\geqslant 2$, then $|J_p(a_+)-J_p(-a_-)|\leqslant(p-1)(|a_+|+|a_-|)^{p-2}|a_++a_-|\leqslant C|h|^p$. If $p\in(1,2)$ and $|A|\geqslant 2M|h|^2$, then $a_+$ and $-a_-$ have the same sign as $A$ and absolute values at least $\frac{|A|}{2}$, so that $|J_p(a_+)-J_p(-a_-)|\leqslant(p-1)(\frac{|A|}{2})^{p-2}M|h|^2$; if $p\in(1,2)$ and $|A|<2M|h|^2$, then $|J_p(a_+)-J_p(-a_-)|\leqslant 2^{2-p}|a_++a_-|^{p-1}\leqslant C|h|^{2(p-1)}$. Writing $h=r\varpi$ with $\varpi\in\mbS$ and $\nu_x:=\frac{\nabla f(x)}{|\nabla f(x)|}$, we have $|A|=|\nabla f(x)|r|\nu_x\cdot\varpi|$. Hence $|A|<2M|h|^2$ implies $|\nu_x\cdot\varpi|<\frac{2Mr}{|\nabla f(x)|}\leqslant\frac{2Mr}{m}$. The set $\{\varpi\in\mbS\colon|\nu_x\cdot\varpi|<\frac{2Mr}{m}\}$ has measure at most $2\uppi\omega_{N-2}M\frac{r}{m}$. Indeed, we write $\nu_x\cdot\varpi=\cos\beta$ with $\beta\in[0,\uppi]$. As in the proof of \eqref{equ:Kbound}, the measure of this set equals $\omega_{N-2}$ times the integral of $\sin^{N-2}\beta$ over $\{\beta\in[0,\uppi]\colon|\cos\beta|<\frac{2Mr}{m}\}$. For $t\in(0,1)$, the set $\{\beta\in[0,\uppi]\colon|\cos\beta|<t\}$ is an interval of length $2\arcsin t$, and $2\arcsin t\leqslant\uppi t$ by the convexity of $\arcsin$ on $[0,1]$. For $t\geqslant 1$, the length of this set is at most $\uppi\leqslant\uppi t$. Since $\sin^{N-2}\beta\leqslant 1$, the claim follows. Also, $\int_{\mbS}|\nu_x\cdot\varpi|^{p-2}\diop\mcH^{N-1}(\varpi)<\infty$ because $p-2>-1$. Hence, for $p\in(1,2)$,
  \begin{equation*}
    |T_\ve(x)-T_{\ve'}(x)|\leqslant C\int_\ve^{\ve'}\left(r^{p-ps-1}+r^{2(p-1)}\,r\,r^{-ps-1}\right)\diop r
    =C\int_\ve^{\ve'}\left(r^{p-ps-1}+r^{2p-ps-2}\right)\diop r,
  \end{equation*}
  and, for $p\geqslant 2$, $|T_\ve(x)-T_{\ve'}(x)|\leqslant C\int_\ve^{\ve'}r^{p-ps-1}\diop r$, where $C$ does not depend on $x\in E$. Since $p-ps>0$ and $2p-ps-1>0$, $\{T_\ve\}$ converges uniformly on $E$ as $\ve\to 0^+$. This proves the existence of $(-\Delta_p)^sf(x)$ and its continuity in $D$.

  Let now $\vp\in C_\mathrm{c}^\infty(D)$ and $E:=\supp\vp$. By the dominated convergence theorem and then by the symmetry of the domain $\{|x-y|>\ve\}$,
  \begin{equation*}
    \langle(-\Delta_p)^sf,\vp\rangle
    =\lim_{\ve\to 0^+}\iint_{\{|x-y|>\ve\}}\frac{J_p\left(f(x)-f(y)\right)\left(\vp(x)-\vp(y)\right)}{|x-y|^{N+ps}}\diopx\diopy
    =\lim_{\ve\to 0^+}\int_{E}T_\ve(x)\vp(x)\diopx,
  \end{equation*}
  where the splitting is legitimate because $\iint_{\{|x-y|>\ve\}}|f(x)-f(y)|^{p-1}|\vp(x)||x-y|^{-N-ps}\diopx\diopy<\infty$. The uniform convergence of $T_\ve$ on $E$ gives \eqref{equ:pointwiseweak}.

  For the last assertion, \eqref{equ:pointwiseweak} yields the required inequality for any nonnegative $\vp\in C_\mathrm{c}^\infty(\Omega\setminus\{0\})$, and Lemma \ref{lem:approximation} together with \eqref{equ:potentialbound} extends it to any $\vp\in X_+(\Omega)$. This finishes the proof of Lemma \ref{lem:pointwiseweak}.
\end{proof}

\subsection{A comparison principle}
We now prove a comparison principle for weak sub- and supersolutions with nonzero exterior data.
\begin{lemma}\label{lem:comparison}
  Let $\Omega\subset\mbRN$ be open, let $V$ be an admissible potential as in \eqref{equ:admissiblepotential}, and let $f_1,f_2\in\mcD$ satisfy $f_1\geqslant 0$ and $f_2>0$ a.e. in $\mbRN$. Assume that
  \begin{enumerate}
    \item\label{item:comparisonsub} $f_1$ is a weak subsolution and $f_2$ is a weak supersolution of $(-\Delta_p)^sf=Vf^{p-1}$ in $\Omega$;
    \item\label{item:comparisonoutside} there exists $\delta\in(0,\infty)$ such that $f_1\leqslant f_2$ a.e. in $\mbRN\setminus\Omega_\delta$;
    \item\label{item:comparisonstrict} the set $\{x\in\mbRN\setminus\Omega\colon f_1(x)<f_2(x)\}$ has positive measure.
  \end{enumerate}
  Then $f_1\leqslant f_2$ a.e. in $\mbRN$.
\end{lemma}
\begin{proof}
  Let $m:=\min\{f_1,f_2\}$ and $M:=\max\{f_1,f_2\}$. Both belong to $\mcD$, and, by \eqref{equ:hardyinequality} and \eqref{equ:sobolev}, the quantity
  \begin{equation*}
    E(f):=\semi{f}^p-\int_{\Omega}V|f|^p
  \end{equation*}
  is finite for $f\in\{f_1,f_2,m,M\}$. By \ref{item:comparisonoutside}, the set $G:=\{f_1>f_2\}$ is contained in $\Omega_\delta$ up to a null set. We also use the following two elementary facts. First, for any $a,b\in\mbR$,
  \begin{equation}\label{equ:convexityp}
    |b|^p\geqslant|a|^p+pJ_p(a)(b-a).
  \end{equation}
  Second, let $g_1,g_2$ be nonnegative and let $\tau\in[0,1]$. Define $\sigma_\tau:=((1-\tau)g_1^p+\tau g_2^p)^{\frac{1}{p}}$. Apply the reverse triangle inequality for the $\ell^p$-norm on $\mbR^2$ to the vectors
\begin{equation*}
  ((1-\tau)^{\frac{1}{p}}g_1(x),\tau^{\frac{1}{p}}g_2(x))
  \quad\text{and}\quad
  ((1-\tau)^{\frac{1}{p}}g_1(y),\tau^{\frac{1}{p}}g_2(y)).
\end{equation*}
This inequality gives
  \begin{equation}\label{equ:hiddenconvexity}
    |\sigma_\tau(x)-\sigma_\tau(y)|^p\leqslant(1-\tau)|g_1(x)-g_1(y)|^p+\tau|g_2(x)-g_2(y)|^p,
  \end{equation}
  and in particular $\semi{\sigma_\tau}^p\leqslant(1-\tau)\semi{g_1}^p+\tau\semi{g_2}^p$.

  \emph{Step 1: $E(M)\geqslant E(f_2)$.} For $\tau\in(0,1]$ let $\sigma_\tau:=((1-\tau)f_2^p+\tau M^p)^{\frac{1}{p}}$ and $\vp_\tau:=\sigma_\tau-f_2$. Since $M\geqslant f_2$, we have $f_2\leqslant\sigma_\tau\leqslant M$, so that $\vp_\tau\geqslant 0$; and $\vp_\tau=0$ wherever $M=f_2$, that is, outside $G$. By \eqref{equ:hiddenconvexity} with $g_1=f_2$ and $g_2=M$, $\semi{\sigma_\tau}<\infty$, and $\sigma_\tau\leqslant M\in L^{\pst}(\mbRN)$; hence $\sigma_\tau\in\mcD$ by \eqref{equ:spaceD}, and so $\vp_\tau\in\mcD$. As $G\subset\Omega_\delta$ up to a null set, $\vp_\tau\in X_+(\Omega)$. We apply \eqref{equ:convexityp} with $a=f_2(x)-f_2(y)$ and $b=\sigma_\tau(x)-\sigma_\tau(y)$, so that $b-a=\vp_\tau(x)-\vp_\tau(y)$; dividing by $p|x-y|^{N+ps}$ and integrating over $\mbR^{2N}$ we obtain the first inequality below, and \ref{item:comparisonsub} gives the second one:
  \begin{equation*}
    \frac{1}{p}\semi{\sigma_\tau}^p-\frac{1}{p}\semi{f_2}^p\geqslant\langle(-\Delta_p)^sf_2,\vp_\tau\rangle\geqslant\int_{\Omega}Vf_2^{p-1}\vp_\tau.
  \end{equation*}
  Combined with $\semi{\sigma_\tau}^p\leqslant(1-\tau)\semi{f_2}^p+\tau\semi{M}^p$, this gives
  \begin{equation*}
    \semi{M}^p-\semi{f_2}^p\geqslant p\int_{\Omega}Vf_2^{p-1}\frac{\sigma_\tau-f_2}{\tau}.
  \end{equation*}
  By the concavity of $t\mapsto t^{\frac{1}{p}}$, $0\leqslant f_2^{p-1}\frac{\sigma_\tau-f_2}{\tau}\leqslant\frac{1}{p}(M^p-f_2^p)$, and the middle term tends to $\frac{1}{p}(M^p-f_2^p)$ as $\tau\to 0^+$, because $f_2>0$ a.e. Since $V(M^p-f_2^p)\in L^1(\Omega)$, the dominated convergence theorem yields $\semi{M}^p-\semi{f_2}^p\geqslant\int_{\Omega}V(M^p-f_2^p)$, that is, $E(M)\geqslant E(f_2)$.

  \emph{Step 2: $E(m)\geqslant E(f_1)$.} For $\tau\in(0,\frac{1}{2}]$ let $\sigma_\tau:=((1-\tau)f_1^p+\tau m^p)^{\frac{1}{p}}$ and $\psi_\tau:=f_1-\sigma_\tau$. Since $m\leqslant f_1$, we have $m\leqslant\sigma_\tau\leqslant f_1$, so that $\psi_\tau\geqslant 0$, and $\psi_\tau=0$ wherever $m=f_1$, that is, outside $G$. As in Step 1, $\sigma_\tau\in\mcD$ by \eqref{equ:hiddenconvexity}, and hence $\psi_\tau\in X_+(\Omega)$. By \eqref{equ:convexityp} with $a=f_1(x)-f_1(y)$ and $b=\sigma_\tau(x)-\sigma_\tau(y)$, and by \ref{item:comparisonsub},
  \begin{equation*}
    \frac{1}{p}\semi{\sigma_\tau}^p-\frac{1}{p}\semi{f_1}^p\geqslant-\langle(-\Delta_p)^sf_1,\psi_\tau\rangle\geqslant-\int_{\Omega}Vf_1^{p-1}\psi_\tau.
  \end{equation*}
  Combined with $\semi{\sigma_\tau}^p\leqslant(1-\tau)\semi{f_1}^p+\tau\semi{m}^p$, this gives
  \begin{equation*}
    \semi{m}^p-\semi{f_1}^p\geqslant-p\int_{\Omega}Vf_1^{p-1}\frac{f_1-\sigma_\tau}{\tau}.
  \end{equation*}
  On $G$ we have $f_1>0$, and, by the mean value theorem applied to $t\mapsto t^{\frac{1}{p}}$ on $[(1-\tau)f_1^p+\tau m^p,f_1^p]$ and the bound $(1-\tau)f_1^p+\tau m^p\geqslant\frac{1}{2}f_1^p$,
  \begin{equation*}
    0\leqslant f_1^{p-1}\frac{f_1-\sigma_\tau}{\tau}\leqslant\frac{2^{\frac{p-1}{p}}}{p}\left(f_1^p-m^p\right),
  \end{equation*}
  while $f_1^{p-1}\frac{f_1-\sigma_\tau}{\tau}\to\frac{1}{p}(f_1^p-m^p)$ as $\tau\to 0^+$. The dominated convergence theorem yields $\semi{m}^p-\semi{f_1}^p\geqslant-\int_{\Omega}V(f_1^p-m^p)$, that is, $E(m)\geqslant E(f_1)$.

  \emph{Step 3: conclusion.} Since $m^p+M^p=f_1^p+f_2^p$, Steps 1 and 2 give
  \begin{equation}\label{equ:latticereverse}
    \semi{m}^p+\semi{M}^p\geqslant\semi{f_1}^p+\semi{f_2}^p.
  \end{equation}
  On the other hand, for any $a_1,a_2,b_1,b_2\in\mbR$,
  \begin{equation}\label{equ:lattice}
    |\max\{a_1,b_1\}-\max\{a_2,b_2\}|^p+|\min\{a_1,b_1\}-\min\{a_2,b_2\}|^p\leqslant|a_1-a_2|^p+|b_1-b_2|^p,
  \end{equation}
  with strict inequality if $(a_1-b_1)(a_2-b_2)<0$. Indeed, if $a_1-b_1$ and $a_2-b_2$ do not have opposite signs, then both sides coincide. If $a_1>b_1$ and $a_2<b_2$, then the left-hand side equals $|a_1-b_2|^p+|b_1-a_2|^p$; the numbers $a_1-b_2$ and $b_1-a_2$ have the same sum as $a_1-a_2$ and $b_1-b_2$, and both lie strictly between $b_1-b_2$ and $a_1-a_2$, so the strict convexity of $t\mapsto|t|^p$ gives the strict inequality; the case $a_1<b_1$ and $a_2>b_2$ is symmetric. Applying \eqref{equ:lattice} with $a_i=f_1(\cdot)$ and $b_i=f_2(\cdot)$ at the points $x$ and $y$, and integrating, we find that $\semi{m}^p+\semi{M}^p\leqslant\semi{f_1}^p+\semi{f_2}^p$, with strict inequality if the set $\{(x,y)\colon(f_1(x)-f_2(x))(f_1(y)-f_2(y))<0\}$ has positive $2N$-dimensional measure. By \eqref{equ:latticereverse}, this set is null. If $G$ had positive measure, then, by \ref{item:comparisonstrict}, the set $G\times\{f_1<f_2\}$ would have positive measure and would be contained in the above set, a contradiction. Hence $|G|=0$, which completes the proof of Lemma \ref{lem:comparison}.
\end{proof}

\section{Local estimates away from the origin}\label{sec:local}

This section is devoted to proving that the positive weak solution $u$ of \eqref{equ:EL} is locally bounded away from the origin in regions of small critical mass. To this end, we establish the following elementary inequality.
\begin{lemma}\label{lem:caccioppolipointwise}
  There exists $C_p\in(0,\infty)$, depending only on $p$, such that, for any $a,b\in\mbR$, $T\in[0,\infty)$, and $t_1,t_2\in[0,1]$, with $A:=(a-T)_+$ and $B:=(b-T)_+$,
  \begin{equation}\label{equ:caccioppolipointwise}
    J_p(a-b)\left(At_1^p-Bt_2^p\right)\geqslant 2^{-p}|At_1-Bt_2|^p-C_p\max\{A,B\}^p|t_1-t_2|^p.
  \end{equation}
\end{lemma}
\begin{proof}
  Both sides of \eqref{equ:caccioppolipointwise} are invariant under the exchange of $(a,t_1)$ and $(b,t_2)$. Hence we may assume that $a\geqslant b$. If $a\leqslant T$, then $A=B=0$ and there is nothing to prove. If $b\leqslant T<a$, then $B=0$ and $a-b\geqslant A$, so that the left-hand side is at least $A^pt_1^p\geqslant 2^{-p}|At_1-Bt_2|^p$. Let now $b>T$, so that $a-b=A-B\geqslant 0$. Then
  \begin{equation*}
    J_p(a-b)\left(At_1^p-Bt_2^p\right)=(A-B)^pt_1^p+(A-B)^{p-1}B\left(t_1^p-t_2^p\right).
  \end{equation*}
  If $t_1\geqslant t_2$, the last term is nonnegative. If $t_1<t_2$, then, by the mean value theorem, $0<t_2^p-t_1^p\leqslant pt_2^{p-1}(t_2-t_1)\leqslant p2^{(p-2)_+}(t_1^{p-1}+(t_2-t_1)^{p-1})(t_2-t_1)$. By the Young inequality with the exponents $\frac{p}{p-1}$ and $p$,
  \begin{equation*}
    p2^{(p-2)_+}(A-B)^{p-1}t_1^{p-1}B(t_2-t_1)\leqslant\frac{1}{2}(A-B)^pt_1^p+C_pB^p(t_2-t_1)^p,
  \end{equation*}
  and $(A-B)^{p-1}B(t_2-t_1)^p\leqslant A^p(t_2-t_1)^p$. Therefore, in both cases,
  \begin{equation*}
    J_p(a-b)\left(At_1^p-Bt_2^p\right)\geqslant\frac{1}{2}(A-B)^pt_1^p-C_pA^p|t_1-t_2|^p.
  \end{equation*}
  Finally, $|At_1-Bt_2|\leqslant(A-B)t_1+B|t_1-t_2|$ implies $(A-B)^pt_1^p\geqslant 2^{1-p}|At_1-Bt_2|^p-B^p|t_1-t_2|^p$, which gives \eqref{equ:caccioppolipointwise}.
\end{proof}

The next proposition is a local boundedness estimate on balls away from the origin with small critical mass.
\begin{proposition}\label{prop:localbound}
  There exist $\ve_\ast\in(0,1)$, depending only on $N$, $s$, and $p$, and $C\in(0,\infty)$, depending only on $N$, $s$, $p$, and $\mu$, such that the following holds. Let $u$ be a positive weak solution of \eqref{equ:EL}, and let $z\in\mbRN$ and $R\in(0,\infty)$ satisfy $|z|\geqslant 5R$ and
  \begin{equation}\label{equ:smallmass}
    \int_{B(z,2R)}u(x)^{\pst}\diopx\leqslant\ve_\ast.
  \end{equation}
  Then
  \begin{equation}\label{equ:localbound}
    \esssup_{B(z,R)}u\leqslant C\left[R^{-r_H}\left(\int_{B(z,2R)}u(x)^{\pst}\diopx\right)^{\frac{1}{\pst}}+\Tail(u;z,2R)\right].
  \end{equation}
\end{proposition}
\begin{proof}
  \emph{Step 1: rescaling.} Let $\tilde u(x):=R^{r_H}u(z+Rx)$ for $x\in\mbRN$. Since $pr_H+ps=N$ and $r_H(\pst-p)=ps$, a change of variables in \eqref{equ:weakEL} shows that $\tilde u\in\mcD$ and, for any $\vp\in\mcD$,
  \begin{equation}\label{equ:rescaledweak}
    \langle(-\Delta_p)^s\tilde u,\vp\rangle=\int_{\mbRN}\left(\tilde V\tilde u^{p-1}+\tilde u^{\pst-1}\right)\vp,
    \qquad
    \tilde V(x):=\frac{\mu R^{ps}}{|z+Rx|^{ps}}.
  \end{equation}
  For $x\in B_2$ we have $|z+Rx|\geqslant|z|-2R\geqslant 3R$, and hence $0\leqslant\tilde V\leqslant\mu$ on $B_2$. Moreover, $\int_{B_2}\tilde u^{\pst}=\int_{B(z,2R)}u^{\pst}$ and $\Tail(\tilde u;0,2)=R^{r_H}\Tail(u;z,2R)$. Since $\esssup_{B(z,R)}u=R^{-r_H}\esssup_{B_1}\tilde u$, it suffices to show that, if $\int_{B_2}\tilde u^{\pst}\leqslant\ve_\ast$, then
  \begin{equation}\label{equ:localboundunit}
    \esssup_{B_1}\tilde u\leqslant C\left(\norm{\tilde u}_{L^{\pst}(B_2)}+\Tail(\tilde u;0,2)\right).
  \end{equation}

  \emph{Step 2: an energy inequality.} Let $T\in(0,\infty)$ be a constant to be chosen later, and, for $j\in\mbN\cup\{0\}$, let
  \begin{equation*}
    r_j:=1+2^{-j},\quad \tilde r_j:=\frac{r_j+r_{j+1}}{2},\quad T_j:=T\left(1-2^{-j-1}\right),\quad \tilde u_j:=(\tilde u-T_j)_+,\quad Y_j:=\int_{B_{r_j}}\tilde u_j^{\pst}.
  \end{equation*}
  Let $\psi_j\in C_\mathrm{c}^\infty(B_{\tilde r_j})$ satisfy $0\leqslant\psi_j\leqslant 1$, $\psi_j=1$ on $B_{r_{j+1}}$, and $|\nabla\psi_j|\leqslant 2^{j+4}$. The function $\vp:=\tilde u_{j+1}\psi_j^p$ is nonnegative and equals $0$ outside $B_{r_j}$. It belongs to $\mcD$ by the proof of Lemma \ref{lem:powerenergy}\ref{item:powerenergytwo}. That proof only requires $\psi_j^p$ to be bounded and Lipschitz with compactly supported gradient. Splitting $\mbR^{2N}$ into $B_{r_j}\times B_{r_j}$ and its complement, and using the symmetry of the kernel, we write $\langle(-\Delta_p)^s\tilde u,\vp\rangle=I_1+2I_2$, where
  \begin{gather*}
    I_1:=\int_{B_{r_j}}\int_{B_{r_j}}\frac{J_p\left(\tilde u(x)-\tilde u(y)\right)\left(\vp(x)-\vp(y)\right)}{|x-y|^{N+ps}}\diopx\diopy,\\
    I_2:=\int_{B_{r_j}}\int_{\mbRN\setminus B_{r_j}}\frac{J_p\left(\tilde u(x)-\tilde u(y)\right)\vp(x)}{|x-y|^{N+ps}}\diopy\diopx.
  \end{gather*}
  By Lemma \ref{lem:caccioppolipointwise} with $a=\tilde u(x)$, $b=\tilde u(y)$, $T=T_{j+1}$, $t_1=\psi_j(x)$, and $t_2=\psi_j(y)$,
  \begin{multline*}
    I_1\geqslant 2^{-p}\int_{B_{r_j}}\int_{B_{r_j}}\frac{|\tilde u_{j+1}\psi_j(x)-\tilde u_{j+1}\psi_j(y)|^p}{|x-y|^{N+ps}}\diopx\diopy\\
    -2C_p\int_{B_{r_j}}\tilde u_{j+1}(x)^p\int_{\mbRN}\frac{|\psi_j(x)-\psi_j(y)|^p}{|x-y|^{N+ps}}\diopy\diopx.
  \end{multline*}
  The inner integral in the last term is at most $\int_{\mbRN}\min\{2^{(j+4)p}|h|^p,1\}|h|^{-N-ps}\diop h\leqslant C2^{jps}$. If $\vp(x)>0$, then $\tilde u(x)>0$, and hence $J_p(\tilde u(x)-\tilde u(y))\geqslant-\tilde u(y)^{p-1}$; therefore
  \begin{equation*}
    I_2\geqslant-\int_{B_{r_j}}\vp(x)\int_{\mbRN\setminus B_{r_j}}\frac{\tilde u(y)^{p-1}}{|x-y|^{N+ps}}\diopy\diopx.
  \end{equation*}
  For $x\in B_{\tilde r_j}$ and $y\notin B_{r_j}$ we have $|x-y|\geqslant|y|-\tilde r_j\geqslant(1-\frac{\tilde r_j}{r_j})|y|\geqslant 2^{-j-3}|y|$, because $r_j-\tilde r_j=2^{-j-2}$ and $r_j\leqslant 2$. We split $\mbRN\setminus B_{r_j}$ into $B_2\setminus B_{r_j}$, where $|y|^{-N-ps}\leqslant 1$, and $\mbRN\setminus B_2$. By the H\"older inequality on $B_2$ and the definition of the tail,
  \begin{align*}
    \int_{\mbRN\setminus B_{r_j}}\frac{\tilde u(y)^{p-1}}{|x-y|^{N+ps}}\diopy
    &\leqslant 2^{(j+3)(N+ps)}\left(\int_{B_2\setminus B_1}\tilde u(y)^{p-1}\diopy+\int_{\mbRN\setminus B_2}\frac{\tilde u(y)^{p-1}}{|y|^{N+ps}}\diopy\right)\\
    &\leqslant C2^{j(N+ps)}\left(\norm{\tilde u}_{L^{\pst}(B_2)}^{p-1}+\Tail(\tilde u;0,2)^{p-1}\right)
    \leqslant C2^{j(N+ps)}T^{p-1},
  \end{align*}
  provided that
  \begin{equation}\label{equ:Tlarge1}
    T\geqslant\norm{\tilde u}_{L^{\pst}(B_2)}+\Tail(\tilde u;0,2).
  \end{equation}
  Moreover, since $\psi_j=0$ outside $B_{\tilde r_j}$ and $\dist(B_{\tilde r_j},\mbRN\setminus B_{r_j})=2^{-j-2}$,
  \begin{equation*}
    \semi{\tilde u_{j+1}\psi_j}^p\leqslant\int_{B_{r_j}}\int_{B_{r_j}}\frac{|\tilde u_{j+1}\psi_j(x)-\tilde u_{j+1}\psi_j(y)|^p}{|x-y|^{N+ps}}\diopx\diopy+C2^{jps}\int_{B_{r_j}}\tilde u_{j+1}^p.
  \end{equation*}
  We now combine these estimates. Since $0\leqslant\vp\leqslant\tilde u_{j+1}$, $\vp=0$ outside $B_{r_j}$, and $\tilde V\leqslant\mu$ on $B_2$, \eqref{equ:rescaledweak} gives
  \begin{equation*}
    I_1+2I_2\leqslant\mu\int_{B_{r_j}}\tilde u^{p-1}\tilde u_{j+1}+\int_{B_{r_j}}\tilde u^{\pst-1}\tilde u_{j+1}.
  \end{equation*}
  By the lower bounds for $I_1$ and $I_2$ and by \eqref{equ:Tlarge1},
  \begin{align*}
    2^{-p}\int_{B_{r_j}}\int_{B_{r_j}}\frac{|\tilde u_{j+1}\psi_j(x)-\tilde u_{j+1}\psi_j(y)|^p}{|x-y|^{N+ps}}\diopx\diopy
    \leqslant{}&I_1+C2^{jps}\int_{B_{r_j}}\tilde u_{j+1}^p\\
    \leqslant{}&I_1+2I_2+C2^{j(N+ps)}T^{p-1}\int_{B_{r_j}}\tilde u_{j+1}\\
    &+C2^{jps}\int_{B_{r_j}}\tilde u_{j+1}^p,
  \end{align*}
  where we used $\vp\leqslant\tilde u_{j+1}$ once more. Adding the estimate for $\semi{\tilde u_{j+1}\psi_j}^p$ and applying \eqref{equ:sobolev} to $\tilde u_{j+1}\psi_j$, we obtain
  \begin{equation}\label{equ:energyj}
    \begin{aligned}
      S_\ast\left(\int_{B_{r_j}}\left(\tilde u_{j+1}\psi_j\right)^{\pst}\right)^{\frac{p}{\pst}}
      \leqslant{}&C2^{jps}\int_{B_{r_j}}\tilde u_{j+1}^p+C2^{j(N+ps)}T^{p-1}\int_{B_{r_j}}\tilde u_{j+1}\\
      &+C\mu\int_{B_{r_j}}\tilde u^{p-1}\tilde u_{j+1}+C\int_{B_{r_j}}\tilde u^{\pst-1}\tilde u_{j+1},
    \end{aligned}
  \end{equation}
  where $C$ depends only on $N$, $s$, and $p$.

  \emph{Step 3: iteration.} On the set $\{\tilde u>T_{j+1}\}$ we have $\tilde u_j\geqslant T_{j+1}-T_j=2^{-j-2}T$ and $\tilde u=\tilde u_j+T_j\leqslant(1+2^{j+2})\tilde u_j\leqslant 2^{j+3}\tilde u_j$; moreover $\tilde u_{j+1}\leqslant \tilde u_j$. Hence
  \begin{align*}
    \int_{B_{r_j}}\tilde u_{j+1}^p&\leqslant\left(2^{-j-2}T\right)^{p-\pst}Y_j,
    &\int_{B_{r_j}}\tilde u_{j+1}&\leqslant\left(2^{-j-2}T\right)^{1-\pst}Y_j,\\
    \int_{B_{r_j}}\tilde u^{p-1}\tilde u_{j+1}&\leqslant 2^{(j+3)(p-1)}\left(2^{-j-2}T\right)^{p-\pst}Y_j,
    &\int_{B_{r_j}}\tilde u^{\pst-1}\tilde u_{j+1}&\leqslant 2^{(j+3)(\pst-1)}Y_j.
  \end{align*}
  Indeed, on the support of $\tilde u_{j+1}$ we have $\tilde u_{j+1}\leqslant\tilde u_j$ and $(2^{-j-2}T)^{-1}\tilde u_j\geqslant 1$, so that, for instance, $\tilde u_{j+1}^p\leqslant\tilde u_j^p\leqslant(2^{-j-2}T)^{p-\pst}\tilde u_j^{\pst}$ and $\tilde u^{p-1}\tilde u_{j+1}\leqslant(2^{j+3}\tilde u_j)^{p-1}\tilde u_j$; the other two bounds are obtained in the same way. Since $\psi_j=1$ on $B_{r_{j+1}}$, the left-hand side of \eqref{equ:energyj} is at least $S_\ast Y_{j+1}^{\frac{p}{\pst}}$. Inserting the four bounds into the right-hand side of \eqref{equ:energyj}, we find
  \begin{align*}
    S_\ast Y_{j+1}^{\frac{p}{\pst}}
    \leqslant{}&C\left[\left(2^{jps}+\mu2^{(j+3)(p-1)}\right)2^{(j+2)(\pst-p)}+2^{j(N+ps)}2^{(j+2)(\pst-1)}\right]T^{p-\pst}Y_j+C2^{(j+3)(\pst-1)}Y_j\\
    \leqslant{}&Cb_1^j\left[1+(1+\mu)T^{p-\pst}\right]Y_j,
  \end{align*}
  with $b_1:=2^{N+ps+\pst}$; indeed, the exponents of $2^j$ in the four coefficients are $ps+\pst-p$, $\pst-1$, $N+ps+\pst-1$, and $\pst-1$, all smaller than $N+ps+\pst$. Raising both sides to the power $\frac{\pst}{p}=1+\gamma_0$, we see that there exist $C_0\in(0,\infty)$ and $b\in(1,\infty)$, depending only on $N$, $s$, and $p$, such that, for any $j\in\mbN\cup\{0\}$,
  \begin{equation*}
    Y_{j+1}\leqslant C_0b^j\left[1+(1+\mu)T^{p-\pst}\right]^{\frac{\pst}{p}}Y_j^{1+\gamma_0},\qquad \gamma_0:=\frac{\pst-p}{p}.
  \end{equation*}
  By induction on $j$, if
  \begin{equation}\label{equ:iterationstart}
    Y_0\leqslant\left(C_0\left[1+(1+\mu)T^{p-\pst}\right]^{\frac{\pst}{p}}\right)^{-\frac{1}{\gamma_0}}b^{-\frac{1}{\gamma_0^2}},
  \end{equation}
  then $Y_j\leqslant b^{-\frac{j}{\gamma_0}}Y_0$ for any $j\in\mbN\cup\{0\}$. Indeed, if this holds for some $j$, then
  \begin{equation*}
    Y_{j+1}\leqslant C_0\left[1+(1+\mu)T^{p-\pst}\right]^{\frac{\pst}{p}}b^{j}b^{-\frac{j(1+\gamma_0)}{\gamma_0}}Y_0^{\gamma_0}Y_0\leqslant b^{-\frac{j+1}{\gamma_0}}Y_0.
  \end{equation*}

  \emph{Step 4: choice of the constants.} We now choose $\ve_\ast$ and $T$. The threshold $\ve_\ast$ will depend only on $N$, $s$, and $p$, whereas $T$ may depend on $\tilde u$. The condition \eqref{equ:iterationstart} involves $T$ only through the factor $(1+\mu)T^{p-\pst}$. We take $T$ so large, in terms of $\norm{\tilde u}_{L^{\pst}(B_2)}$, that this factor is compensated by the smallness of $Y_0$.

  If $\norm{\tilde u}_{L^{\pst}(B_2)}=0$, then $\tilde u=0$ a.e. in $B_2$, and \eqref{equ:localboundunit} is trivial. Hence we assume that $\norm{\tilde u}_{L^{\pst}(B_2)}>0$. Let
  \begin{equation}\label{equ:epsstar}
    \ve_\ast:=\min\left\{\frac{1}{2},\ \left(2^{\frac{\pst}{p}}C_0\right)^{-\frac{1}{\gamma_0}}b^{-\frac{1}{\gamma_0^2}}\right\}.
  \end{equation}
  Since $C_0$, $b$, and $\gamma_0$ depend only on $N$, $s$, and $p$, so does $\ve_\ast$. The bound $\ve_\ast\leqslant\frac{1}{2}$ only guarantees that $\ve_\ast\in(0,1)$. We choose
  \begin{equation}\label{equ:Tchoice}
    T:=\norm{\tilde u}_{L^{\pst}(B_2)}+\Tail(\tilde u;0,2)+(1+\mu)^{\frac{1}{\pst-p}}\ve_\ast^{-\frac{1}{\pst}}\norm{\tilde u}_{L^{\pst}(B_2)}.
  \end{equation}
  Then $T>0$. Moreover, \eqref{equ:Tlarge1} holds, because the last term in \eqref{equ:Tchoice} is nonnegative. Hence Steps 2 and 3 apply with this $T$.

  Assume now that $\int_{B_2}\tilde u^{\pst}\leqslant\ve_\ast$. We show that \eqref{equ:iterationstart} holds. Since $\gamma_0=\frac{\pst-p}{p}$, we have $\frac{1}{\gamma_0}\cdot\frac{\pst}{p}=\frac{\pst}{\pst-p}$. Hence the right-hand side of \eqref{equ:iterationstart} equals
  \begin{equation*}
    C_0^{-\frac{1}{\gamma_0}}b^{-\frac{1}{\gamma_0^2}}\left[1+(1+\mu)T^{p-\pst}\right]^{-\frac{\pst}{\pst-p}}.
  \end{equation*}
  We use $1+(1+\mu)T^{p-\pst}\leqslant 2\max\{1,(1+\mu)T^{p-\pst}\}$ and $2^{\frac{\pst}{\pst-p}}=(2^{\frac{\pst}{p}})^{\frac{1}{\gamma_0}}$. By \eqref{equ:epsstar}, the right-hand side of \eqref{equ:iterationstart} is therefore at least
  \begin{equation*}
    \left(2^{\frac{\pst}{p}}C_0\right)^{-\frac{1}{\gamma_0}}b^{-\frac{1}{\gamma_0^2}}\max\left\{1,(1+\mu)T^{p-\pst}\right\}^{-\frac{\pst}{\pst-p}}
    \geqslant\ve_\ast\max\left\{1,(1+\mu)T^{p-\pst}\right\}^{-\frac{\pst}{\pst-p}}.
  \end{equation*}
  Consequently, \eqref{equ:iterationstart} follows from the two inequalities
  \begin{equation}\label{equ:twoconditions}
    Y_0\leqslant\ve_\ast
    \qquad\text{and}\qquad
    \left((1+\mu)T^{p-\pst}\right)^{\frac{\pst}{\pst-p}}Y_0\leqslant\ve_\ast.
  \end{equation}
  Since $r_0=2$ and $T_0=\frac{T}{2}$, we have $Y_0=\int_{B_2}(\tilde u-\frac{T}{2})_+^{\pst}\leqslant\int_{B_2}\tilde u^{\pst}$. Thus the first inequality in \eqref{equ:twoconditions} follows from the smallness assumption. For the second one, note that $((1+\mu)T^{p-\pst})^{\frac{\pst}{\pst-p}}=(1+\mu)^{\frac{\pst}{\pst-p}}T^{-\pst}$. By \eqref{equ:Tchoice}, $T\geqslant(1+\mu)^{\frac{1}{\pst-p}}\ve_\ast^{-\frac{1}{\pst}}\norm{\tilde u}_{L^{\pst}(B_2)}$. Raising this inequality to the power $\pst$, we obtain
  \begin{equation*}
    (1+\mu)^{\frac{\pst}{\pst-p}}T^{-\pst}\leqslant\ve_\ast\left(\int_{B_2}\tilde u^{\pst}\right)^{-1}.
  \end{equation*}
  Multiplying by $Y_0\leqslant\int_{B_2}\tilde u^{\pst}$, we obtain the second inequality in \eqref{equ:twoconditions}. This second inequality does not use the smallness assumption. The smallness of the critical mass enters only through the first inequality, and the threshold $\ve_\ast$ does not depend on $T$.

  Hence \eqref{equ:iterationstart} holds, and Step 3 gives $Y_j\leqslant b^{-\frac{j}{\gamma_0}}Y_0$ for any $j\in\mbN\cup\{0\}$. Since $b>1$, $Y_j\to 0$ as $j\to\infty$. For every $j$ we have $T_j\leqslant T$, so that $(\tilde u-T)_+\leqslant\tilde u_j$; moreover, $B_1\subset B_{r_j}$. Therefore
  \begin{equation*}
    \int_{B_1}(\tilde u-T)_+^{\pst}\leqslant Y_j\quad\text{for any } j\in\mbN\cup\{0\}.
  \end{equation*}
  Letting $j\to\infty$, we obtain $\int_{B_1}(\tilde u-T)_+^{\pst}=0$, that is, $\tilde u\leqslant T$ a.e. in $B_1$. Finally, by \eqref{equ:Tchoice},
  \begin{equation*}
    T\leqslant\left(1+(1+\mu)^{\frac{1}{\pst-p}}\ve_\ast^{-\frac{1}{\pst}}\right)\left(\norm{\tilde u}_{L^{\pst}(B_2)}+\Tail(\tilde u;0,2)\right).
  \end{equation*}
  This proves \eqref{equ:localboundunit} with $C:=1+(1+\mu)^{\frac{1}{\pst-p}}\ve_\ast^{-\frac{1}{\pst}}$, which depends only on $N$, $s$, $p$, and $\mu$. Hence Proposition \ref{prop:localbound} is proved.
\end{proof}

We next recall the weak Harnack inequality of Di Castro, Kuusi and Palatucci \cite[Theorem 1.2]{bib:DKP2014}, in the form needed here. It is stated in \cite{bib:DKP2014} for functions in $W^{s,p}(\mbRN)$. We reduce the following lemma to that case by a cutoff.
\begin{lemma}\label{lem:weakharnack}
  Let $f\in\mcD$ be nonnegative in $\mbRN$, let $z\in\mbRN$ and $R\in(0,\infty)$, and assume that $\langle(-\Delta_p)^sf,\vp\rangle\geqslant 0$ for any nonnegative $\vp\in C_\mathrm{c}^\infty(B(z,2R))$. Then there exists $c\in(0,1]$, depending only on $N$, $s$, and $p$, such that
  \begin{equation*}
    \essinf_{B(z,R)}f\geqslant c\left(\fint_{B(z,R)}f^{\frac{p-1}{2}}\right)^{\frac{2}{p-1}}.
  \end{equation*}
\end{lemma}
\begin{proof}
  Fix $\chi_z\in C_\mathrm{c}^\infty(B(z,4R))$ with $0\leqslant\chi_z\leqslant 1$ and $\chi_z=1$ on $B(z,3R)$, and let $\tilde f:=\chi_zf$. By Lemma \ref{lem:powerenergy}\ref{item:powerenergytwo}, $\tilde f\in\mcD$. Since $\tilde f\in L^{\pst}(\mbRN)$ has compact support, $\tilde f\in L^p(\mbRN)$, and hence $\tilde f\in W^{s,p}(\mbRN)$.

  Let $\vp\in C_\mathrm{c}^\infty(B(z,2R))$ be nonnegative. Since $\tilde f=f$ in $B(z,3R)$ and $\vp=0$ outside $B(z,2R)$, the integrands of $\langle(-\Delta_p)^s\tilde f,\vp\rangle$ and $\langle(-\Delta_p)^sf,\vp\rangle$ differ only at pairs $(x,y)$ with one point in $B(z,2R)$ and the other outside $B(z,3R)$, and for such pairs $|x-y|\geqslant R$. By the symmetry of the integrands in $x$ and $y$,
  \begin{multline*}
    \langle(-\Delta_p)^s\tilde f,\vp\rangle-\langle(-\Delta_p)^sf,\vp\rangle\\
    =2\int_{B(z,2R)}\vp(x)\int_{\mbRN\setminus B(z,3R)}\frac{J_p\left(f(x)-\tilde f(y)\right)-J_p\left(f(x)-f(y)\right)}{|x-y|^{N+ps}}\diopy\diopx,
  \end{multline*}
  where the integrals converge absolutely, because $|x-y|\geqslant R$ and $f\in L^{\pst}(\mbRN)$. Since $0\leqslant\tilde f\leqslant f$ and $J_p$ is increasing, the right-hand side is nonnegative. Hence $\langle(-\Delta_p)^s\tilde f,\vp\rangle\geqslant 0$ for every nonnegative $\vp\in C_\mathrm{c}^\infty(B(z,2R))$. Since $|\langle(-\Delta_p)^s\tilde f,\eta\rangle|\leqslant\semi{\tilde f}^{p-1}\semi{\eta}$, and since every nonnegative function in the closure $W_0^{s,p}(B(z,2R))$ of $C_\mathrm{c}^\infty(B(z,2R))$ in $W^{s,p}(\mbRN)$ is the limit of nonnegative functions in $C_\mathrm{c}^\infty(B(z,2R))$, the same inequality holds for every nonnegative $\eta\in W_0^{s,p}(B(z,2R))$. Thus $\tilde f\in W^{s,p}(\mbRN)$ is a weak supersolution in $B(z,2R)$ in the sense of \cite{bib:DKP2014}. Moreover, $\tilde f$ is nonnegative in $\mbRN$, so that $\Tail(\tilde f_-;x,\rho)=0$ for any $x\in\mbRN$ and $\rho\in(0,\infty)$. Let $t:=\frac{p-1}{2}$, so that $0<t<\frac{N(p-1)}{N-ps}$, and, for $x\in\ol{B(z,R)}$, let
  \begin{equation*}
    M(x):=\left(\fint_{B(x,\frac{R}{8})}\tilde f^{t}\right)^{\frac{1}{t}}.
  \end{equation*}
  Let $x\in\ol{B(z,R)}$. Then $B(x,\frac{R}{2})\subset B(z,2R)$. The weak Harnack inequality \cite[Theorem 1.2]{bib:DKP2014}, applied to $\tilde f$ with the domain $B(z,2R)$, the ball $B(x,\frac{R}{2})$, the inner radius $\frac{R}{8}$, and the exponent $t$, gives $C_1\in[1,\infty)$, depending only on $N$, $s$, and $p$, such that
  \begin{equation*}
    M(x)\leqslant C_1\essinf_{B(x,\frac{R}{4})}\tilde f\qquad\text{for any }x\in\ol{B(z,R)}.
  \end{equation*}
  If $x,y\in\ol{B(z,R)}$ and $|x-y|\leqslant\frac{R}{8}$, then $B(y,\frac{R}{8})\subset B(x,\frac{R}{4})$, and hence
  \begin{equation*}
    M(x)\leqslant C_1\essinf_{B(y,\frac{R}{8})}\tilde f\leqslant C_1M(y).
  \end{equation*}
  Any two points of $\ol{B(z,R)}$ are joined by a segment in $\ol{B(z,R)}$ of length at most $2R$, and this segment is the union of $16$ segments of length at most $\frac{R}{8}$. Therefore $M(x)\leqslant C_1^{16}M(y)$ for any $x,y\in\ol{B(z,R)}$.

  By the compactness of $\ol{B_1}$, there exist $\kappa\in\mbN$, depending only on $N$, and $y_1,\dots,y_\kappa\in\ol{B_1}$ such that $\ol{B_1}\subset\bigcup_{i=1}^{\kappa}B(y_i,\frac{1}{8})$. Let $x_i:=z+Ry_i$ for $i\in\{1,\dots,\kappa\}$. Then $x_i\in\ol{B(z,R)}$ and $B(z,R)\subset\bigcup_{i=1}^{\kappa}B(x_i,\frac{R}{8})$; in particular, $\kappa\geqslant 8^N$. Hence
  \begin{equation*}
    \int_{B(z,R)}\tilde f^{t}\leqslant\sum_{i=1}^{\kappa}\int_{B(x_i,\frac{R}{8})}\tilde f^{t}=8^{-N}|B(z,R)|\sum_{i=1}^{\kappa}M(x_i)^{t}\leqslant \kappa 8^{-N}|B(z,R)|\,C_1^{16t}\min_{1\leqslant i\leqslant \kappa}M(x_i)^{t},
  \end{equation*}
  and, since the balls $B(x_i,\frac{R}{8})$ cover $B(z,R)$,
  \begin{equation*}
    \essinf_{B(z,R)}\tilde f\geqslant\min_{1\leqslant i\leqslant \kappa}\essinf_{B(x_i,\frac{R}{8})}\tilde f\geqslant\min_{1\leqslant i\leqslant \kappa}\essinf_{B(x_i,\frac{R}{4})}\tilde f\geqslant C_1^{-1}\min_{1\leqslant i\leqslant \kappa}M(x_i).
  \end{equation*}
  Combining the last two estimates, we obtain
  \begin{equation*}
    \left(\fint_{B(z,R)}\tilde f^{\frac{p-1}{2}}\right)^{\frac{2}{p-1}}\leqslant\left(\kappa 8^{-N}\right)^{\frac{1}{t}}C_1^{17}\essinf_{B(z,R)}\tilde f.
  \end{equation*}
  Since $\tilde f=f$ in $B(z,R)$, the assertion follows with $c:=(\kappa 8^{-N})^{-\frac{1}{t}}C_1^{-17}\in(0,1]$.
\end{proof}

\begin{corollary}\label{cor:annulus}
  Let $u$ be a positive weak solution of \eqref{equ:EL}. Then, for any $z\in\mbRN$ and $R\in(0,\infty)$, $\essinf_{B(z,R)}u>0$. Moreover, for any $0<\rho_1<\rho_2<\infty$,
  \begin{equation*}
    0<\essinf_{\{\rho_1<|x|<\rho_2\}}u\leqslant\esssup_{\{\rho_1<|x|<\rho_2\}}u<\infty.
  \end{equation*}
\end{corollary}
\begin{proof}
  As observed after \eqref{equ:Vu}, $\langle(-\Delta_p)^su,\vp\rangle\geqslant 0$ for any nonnegative $\vp\in\mcD$. Since $u>0$ a.e., Lemma \ref{lem:weakharnack} gives $\essinf_{B(z,R)}u>0$ for any ball, and the positive lower bound on the annulus follows by covering it with finitely many balls. For the upper bound, since $u^{\pst}\in L^1(\mbRN)$, the absolute continuity of the integral yields $\eta>0$ such that $\int_Eu^{\pst}\leqslant\ve_\ast$ for any measurable set $E$ with $|E|<\eta$. Choosing $R_0\in(0,\frac{\rho_1}{10})$ with $|B_{2R_0}|<\eta$, we obtain $\int_{B(z,2R_0)}u^{\pst}\leqslant\ve_\ast$ for any $z\in\mbRN$. Then $|z|\geqslant\rho_1\geqslant 5R_0$ for any $z$ with $\rho_1\leqslant|z|\leqslant\rho_2$, and Proposition \ref{prop:localbound}, together with \eqref{equ:tailfinite}, shows that $u$ is essentially bounded on $B(z,R_0)$ by a constant independent of $z$. Covering the annulus by finitely many such balls completes the proof of Corollary \ref{cor:annulus}.
\end{proof}

\section{A preliminary decay of the critical mass}\label{sec:mass}

In this section, we prove some extra integrability near the origin and at infinity on the positive weak solution $u$ of \eqref{equ:EL}, which allows to establish the required pointwise bounds.

For $\delta\in(0,s)$ let
\begin{equation}\label{equ:kernelconstant}
  \mcI(\delta):=\int_{\mbRN}\frac{\left(\mathrm{e}^{\delta|\ln|y||}-1\right)^p}{|e-y|^{N+ps}}\diopy,
\end{equation}
where $e\in\mbS$; by rotation invariance $\mcI(\delta)$ does not depend on $e$. The integrand is bounded by $C|y|^{-p\delta}$ near $0$, by $C|e-y|^{p-N-ps}$ near $e$, and by $C|y|^{p\delta-N-ps}$ near infinity. Since $p\delta<ps<N$ and $p-ps>0$, $\mcI(\delta)$ is finite. Moreover, for $0<\delta\leqslant\frac{s}{2}$, the convexity of $t\mapsto\mathrm{e}^{t}$ gives $\mathrm{e}^{\delta t}-1\leqslant\frac{2\delta}{s}(\mathrm{e}^{\frac{st}{2}}-1)$ for $t\geqslant 0$, whence
\begin{equation}\label{equ:kernelsmall}
  \mcI(\delta)\leqslant\left(\frac{2\delta}{s}\right)^p\mcI\left(\frac{s}{2}\right)\longrightarrow 0\quad\text{as } \delta\to 0^+.
\end{equation}

\begin{proposition}\label{prop:weighted}
  There exists $\delta_0\in(0,\frac{s}{2}]$, depending only on $N$, $s$, $p$, and $\mu$, such that, for any positive weak solution $u$ of \eqref{equ:EL} and any $\delta\in(0,\delta_0]$,
  \begin{equation}\label{equ:weighted}
    \int_{B_1}u(x)^{\pst}|x|^{-\pst\delta}\diopx+\int_{\mbRN\setminus B_1}u(x)^{\pst}|x|^{\pst\delta}\diopx<\infty.
  \end{equation}
  Consequently, with $\theta:=\pst\delta$, there exists $C\in(0,\infty)$, depending on $u$ and $\delta$, such that
  \begin{equation}\label{equ:H1}
    m_0(R)\leqslant CR^{\theta}\quad\text{for any } R\in(0,1],
    \qquad
    m_\infty(R)\leqslant CR^{-\theta}\quad\text{for any } R\in[1,\infty).
  \end{equation}
\end{proposition}
\begin{proof}
  \emph{Step 1: a two-point inequality.} Let $\ve_1:=\frac{1}{4}(1-\frac{\mu}{\mathcal{C}_{N,s,p}})\in(0,\frac{1}{4})$, and $z_0:=\frac{\ve_1}{p-1}$, and let
  \begin{equation*}
    c_1:=\frac{1-\ve_1}{(1+z_0)^{p-1}},\qquad C_1:=\left(\frac{p-1}{\ve_1}\right)^{p-1}+(1-\ve_1)z_0^{1-p}.
  \end{equation*}
  We claim that, for any $a,b\in[0,\infty)$ and $t_1,t_2\in(0,\infty)$,
  \begin{equation}\label{equ:twopoint}
    J_p(a-b)\left(at_1^p-bt_2^p\right)\geqslant c_1|at_1-bt_2|^p-C_1\left(a^p+b^p\right)|t_1-t_2|^p.
  \end{equation}
  Both sides are invariant under the exchange of $(a,t_1)$ and $(b,t_2)$, so we may assume $t_1\geqslant t_2$. Let $X:=t_1|a-b|$ and $Y:=b(t_1-t_2)$. Since $at_1^p-bt_2^p=(a-b)t_1^p+b(t_1^p-t_2^p)$ and $0\leqslant t_1^p-t_2^p\leqslant pt_1^{p-1}(t_1-t_2)$,
  \begin{equation*}
    J_p(a-b)\left(at_1^p-bt_2^p\right)\geqslant|a-b|^pt_1^p-|a-b|^{p-1}b\left(t_1^p-t_2^p\right)\geqslant X^p-pX^{p-1}Y.
  \end{equation*}
  By the Young inequality with the exponents $\frac{p}{p-1}$ and $p$, $pX^{p-1}Y\leqslant \ve_1X^p+(\frac{p-1}{\ve_1})^{p-1}Y^p$. Moreover, $|at_1-bt_2|=|(a-b)t_1+b(t_1-t_2)|\leqslant X+Y$, and the convexity of $t\mapsto t^p$, applied to $X+Y=\frac{1}{1+z_0}(1+z_0)X+\frac{z_0}{1+z_0}\frac{1+z_0}{z_0}Y$, gives
  \begin{equation*}
    (X+Y)^p\leqslant(1+z_0)^{p-1}X^p+(1+z_0^{-1})^{p-1}Y^p .
  \end{equation*}
  Since $(1+z_0)^{1-p}(1+z_0^{-1})^{p-1}=z_0^{1-p}$, we obtain
  \begin{equation*}
    J_p(a-b)\left(at_1^p-bt_2^p\right)\geqslant(1-\ve_1)X^p-\left(\frac{p-1}{\ve_1}\right)^{p-1}Y^p\geqslant c_1(X+Y)^p-C_1Y^p\geqslant c_1|at_1-bt_2|^p-C_1Y^p,
  \end{equation*}
  and $Y^p\leqslant b^p|t_1-t_2|^p\leqslant(a^p+b^p)|t_1-t_2|^p$, which proves \eqref{equ:twopoint}. Since $(1+z_0)^{p-1}\leqslant\mathrm{e}^{\ve_1}$ and $\mathrm{e}^{-\ve_1}\geqslant 1-\ve_1$, we have $c_1\geqslant(1-\ve_1)^2\geqslant 1-2\ve_1=\frac{1}{2}(1+\frac{\mu}{\mathcal{C}_{N,s,p}})$. By \eqref{equ:kernelsmall}, we can choose $\delta_0\in(0,\frac{s}{2}]$, depending only on $N$, $s$, $p$, and $\mu$, such that
  \begin{equation}\label{equ:dgap}
    d:=c_1-\frac{\mu}{\mathcal{C}_{N,s,p}}-\frac{2C_1\mcI(\delta)}{\mathcal{C}_{N,s,p}}\geqslant\frac{1}{4}\left(1-\frac{\mu}{\mathcal{C}_{N,s,p}}\right)\quad\text{for any } \delta\in(0,\delta_0].
  \end{equation}

  \emph{Step 2: the weight near the origin.} Fix $\delta\in(0,\delta_0]$ and $M\in[1,\infty)$, and let
  \begin{equation*}
    \xi(x):=\min\left\{M,\max\{1,|x|^{-\delta}\}\right\}.
  \end{equation*}
  Since $\ln$ is increasing, $\ln\xi(x)=\min\{\ln M,\max\{0,-\delta\ln|x|\}\}$. The maps $t\mapsto\max\{0,t\}$ and $t\mapsto\min\{\ln M,t\}$ on $\mbR$ are $1$-Lipschitz, and so is their composition. Hence $\ln\xi$ is obtained from $x\mapsto-\delta\ln|x|$ by composition with a $1$-Lipschitz function, so that
  \begin{equation*}
    |\ln\xi(x)-\ln\xi(y)|\leqslant\delta\left|\ln|x|-\ln|y|\right|{}=\delta\left|\ln\frac{|y|}{|x|}\right|.
  \end{equation*}
  Taking exponentials, we obtain
  \begin{equation*}
    \mathrm{e}^{-\delta|\ln\frac{|y|}{|x|}|}\leqslant\frac{\xi(y)}{\xi(x)}\leqslant\mathrm{e}^{\delta|\ln\frac{|y|}{|x|}|}.
  \end{equation*}
  If $\xi(y)\geqslant\xi(x)$, the second inequality gives $\xi(y)-\xi(x)\leqslant\xi(x)(\mathrm{e}^{\delta|\ln\frac{|y|}{|x|}|}-1)$. If $\xi(y)<\xi(x)$, the first inequality gives $\xi(x)-\xi(y)\leqslant\xi(x)(1-\mathrm{e}^{-\delta|\ln\frac{|y|}{|x|}|})$. Since $1-\mathrm{e}^{-a}\leqslant\mathrm{e}^{a}-1$ for $a\geqslant 0$, in both cases
  \begin{equation*}
    |\xi(x)-\xi(y)|\leqslant\xi(x)\left(\mathrm{e}^{\delta|\ln\frac{|y|}{|x|}|}-1\right).
  \end{equation*} By the change of variables $y=|x|Oy'$, with $O$ a rotation mapping $e$ to $\frac{x}{|x|}$, and by the definition \eqref{equ:kernelconstant} of $\mcI(\delta)$,
  \begin{equation}\label{equ:weightkernel}
    \int_{\mbRN}\frac{|\xi(x)-\xi(y)|^p}{|x-y|^{N+ps}}\diopy\leqslant\mcI(\delta)\,\xi(x)^p|x|^{-ps},\qquad x\neq 0.
  \end{equation}
  Let $h:=u\xi$ and $\vp:=u\xi^p$. Since $1\leqslant\xi\leqslant M$ and $|\xi^p(x)-\xi^p(y)|\leqslant pM^{p-1}|\xi(x)-\xi(y)|$, we have
  \begin{equation*}
    |\vp(x)-\vp(y)|\leqslant M^p|u(x)-u(y)|+pM^{p-1}u(y)|\xi(x)-\xi(y)|,
  \end{equation*}
  and, by \eqref{equ:weightkernel} with the roles of $x$ and $y$ exchanged and by \eqref{equ:hardyinequality},
  \begin{equation*}
    \int_{\mbRN}\int_{\mbRN}\frac{u(y)^p|\xi(x)-\xi(y)|^p}{|x-y|^{N+ps}}\diopx\diopy\leqslant\mcI(\delta)M^p\int_{\mbRN}\frac{u(y)^p}{|y|^{ps}}\diopy<\infty .
  \end{equation*}
  Hence $\semi{\vp}<\infty$, and $\vp\in\mcD$ by \eqref{equ:spaceD}; the same argument gives $h\in\mcD$. We test \eqref{equ:weakEL} with $\vp$. By \eqref{equ:twopoint} with $a=u(x)$, $b=u(y)$, $t_1=\xi(x)$, and $t_2=\xi(y)$, by the symmetry of the kernel, and by \eqref{equ:weightkernel},
  \begin{align*}
    \langle(-\Delta_p)^su,\vp\rangle
    &\geqslant c_1\semi{h}^p-C_1\int_{\mbRN}\int_{\mbRN}\frac{\left(u(x)^p+u(y)^p\right)|\xi(x)-\xi(y)|^p}{|x-y|^{N+ps}}\diopx\diopy\\
    &=c_1\semi{h}^p-2C_1\int_{\mbRN}u(x)^p\int_{\mbRN}\frac{|\xi(x)-\xi(y)|^p}{|x-y|^{N+ps}}\diopy\diopx\\
    &\geqslant c_1\semi{h}^p-2C_1\mcI(\delta)\int_{\mbRN}\frac{h^p}{|x|^{ps}}\diopx .
  \end{align*}
  The right-hand side of \eqref{equ:weakEL} with this $\vp$ equals $\mu\int h^p|x|^{-ps}+\int u^{\pst-p}h^p$. Therefore
  \begin{equation*}
    c_1\semi{h}^p-2C_1\mcI(\delta)\int_{\mbRN}\frac{h^p}{|x|^{ps}}\diopx\leqslant\mu\int_{\mbRN}\frac{h^p}{|x|^{ps}}\diopx+\int_{\mbRN}u^{\pst-p}h^p\diopx.
  \end{equation*}
  Since $\int h^p|x|^{-ps}\leqslant\mathcal{C}_{N,s,p}^{-1}\semi{h}^p$ by \eqref{equ:hardyinequality}, the left-hand side minus the first term on the right-hand side is at least $(c_1-\frac{\mu}{\mathcal{C}_{N,s,p}}-\frac{2C_1\mcI(\delta)}{\mathcal{C}_{N,s,p}})\semi{h}^p$. Hence, by \eqref{equ:dgap},
  \begin{equation}\label{equ:weightedenergy}
    d\,\semi{h}^p\leqslant\int_{\mbRN}u^{\pst-p}h^p\diopx.
  \end{equation}
  Choose $r_0\in(0,1)$, depending on $u$ and $\delta$, such that $m_0(r_0)^{\frac{\pst-p}{\pst}}\leqslant\frac{S_\ast d}{2}$. By the H\"older inequality and \eqref{equ:sobolev},
  \begin{equation*}
    \int_{B_{r_0}}u^{\pst-p}h^p\leqslant m_0(r_0)^{\frac{\pst-p}{\pst}}\norm{h}_{L^{\pst}(\mbRN)}^p\leqslant\frac{d}{2}\semi{h}^p,
  \end{equation*}
  while $\xi\leqslant r_0^{-\delta}$ outside $B_{r_0}$ gives $\int_{\mbRN\setminus B_{r_0}}u^{\pst-p}h^p\leqslant r_0^{-p\delta}\norm{u}_{L^{\pst}(\mbRN)}^{\pst}$. Hence, by \eqref{equ:weightedenergy} and \eqref{equ:sobolev},
  \begin{equation*}
    \int_{\mbRN}u^{\pst}\min\left\{M,\max\{1,|x|^{-\delta}\}\right\}^{\pst}\diopx
    =\norm{h}_{L^{\pst}(\mbRN)}^{\pst}\leqslant\left(\frac{2r_0^{-p\delta}}{S_\ast d}\norm{u}_{L^{\pst}(\mbRN)}^{\pst}\right)^{\frac{\pst}{p}},
  \end{equation*}
  and the right-hand side does not depend on $M$. Letting $M\to\infty$, the monotone convergence theorem gives the finiteness of the first integral in \eqref{equ:weighted}.

  \emph{Step 3: the weight at infinity.} We repeat Step 2 with $\xi(x):=\min\{M,\max\{1,|x|^{\delta}\}\}$, for which \eqref{equ:weightkernel} holds with the same proof. Choosing $L\in(1,\infty)$ such that $m_\infty(L)^{\frac{\pst-p}{\pst}}\leqslant\frac{S_\ast d}{2}$, and using $\xi\leqslant L^{\delta}$ in $B_L$, we obtain the finiteness of the second integral in \eqref{equ:weighted}.

  Finally, for $R\in(0,1]$ and $R'\in[1,\infty)$,
  \begin{equation*}
    m_0(R)\leqslant R^{\pst\delta}\int_{B_1}u^{\pst}|x|^{-\pst\delta}\diopx
    \qquad\text{and}\qquad
    m_\infty(R')\leqslant (R')^{-\pst\delta}\int_{\mbRN\setminus B_1}u^{\pst}|x|^{\pst\delta}\diopx.
  \end{equation*}
  This proves \eqref{equ:H1} and completes the proof of Proposition \ref{prop:weighted}.
\end{proof}

We now combine Propositions \ref{prop:localbound} and \ref{prop:weighted}. Since $r_H(p-1)=\frac{(N-ps)(p-1)}{p}<N$, the number
\begin{equation}\label{equ:deltaone}
  \delta_1:=\min\left\{\delta_0,\ \frac{r_H}{2},\ \frac{N-r_H(p-1)}{2(p-1)}\right\}
\end{equation}
is positive; it depends only on $N$, $s$, $p$, and $\mu$.

\begin{lemma}\label{lem:coarse}
  Let $u$ be a positive weak solution of \eqref{equ:EL} and let $\delta:=\delta_1$. Then there exist $\rho_1\in(0,1)$, $R_1\in(1,\infty)$, and $C\in(0,\infty)$, depending on $u$, such that
  \begin{equation}\label{equ:coarse}
    \begin{aligned}
      u(x)&\leqslant C|x|^{-r_H+\delta}\quad\text{for a.e. } 0<|x|<\rho_1,\\
      u(x)&\leqslant C|x|^{-r_H-\delta}\quad\text{for a.e. } |x|>R_1.
    \end{aligned}
  \end{equation}
  Consequently, since $r_H(\pst-p)=ps$,
  \begin{equation}\label{equ:coarsecritical}
    \begin{aligned}
      u(x)^{\pst-p}&\leqslant C^{\pst-p}|x|^{-ps+\delta(\pst-p)}\quad\text{for a.e. } 0<|x|<\rho_1,\\
      u(x)^{\pst-p}&\leqslant C^{\pst-p}|x|^{-ps-\delta(\pst-p)}\quad\text{for a.e. } |x|>R_1.
    \end{aligned}
  \end{equation}
\end{lemma}
\begin{proof}
  Let $\theta:=\pst\delta$, so that $\frac{\theta}{\pst}=\delta$ and \eqref{equ:H1} holds. We use $\frac{N}{\pst}=r_H$ repeatedly. Choose finitely many points $z_1,\ldots,z_{n}\in\{1\leqslant|z|\leqslant 2\}$, with $n$ depending only on $N$, such that the balls $B(z_i,\frac{1}{8})$ cover $\{1\leqslant|z|\leqslant 2\}$. For $r\in(0,\infty)$ and $i\in\{1,\ldots,n\}$ we apply Proposition \ref{prop:localbound} with $z=rz_i$ and $R=\frac{r}{8}$; note that $|rz_i|\geqslant r=8R$ and $B(rz_i,\frac{r}{4})\subset\{\frac{3r}{4}<|x|<\frac{9r}{4}\}$.

  \emph{Near the origin.} By \eqref{equ:H1}, there exists $r_1\in(0,\frac{1}{8})$ such that $m_0(3r)\leqslant\ve_\ast$ for $r\in(0,r_1)$, so \eqref{equ:smallmass} holds. For $r\in(0,r_1)$, the first term in \eqref{equ:localbound} is at most $Cr^{-r_H}m_0(3r)^{\frac{1}{\pst}}\leqslant Cr^{-r_H+\delta}$. For the tail, we split $\mbRN\setminus B(rz_i,\frac{r}{4})$ into $\{|y|<4r\}$, $\{4r\leqslant|y|<\frac{1}{2}\}$, and $\{|y|\geqslant\frac{1}{2}\}$. On the first set the kernel is at most $(\frac{r}{4})^{-N-ps}$, and the H\"older inequality with \eqref{equ:H1} gives
  \begin{equation*}
    \int_{B_{4r}}u^{p-1}\leqslant|B_{4r}|^{1-\frac{p-1}{\pst}}m_0(4r)^{\frac{p-1}{\pst}}\leqslant Cr^{N-(p-1)(r_H-\delta)}.
  \end{equation*}
  On the second set we have $|y-rz_i|\geqslant|y|-2r\geqslant\frac{|y|}{2}$. For $k\geqslant 2$ with $2^{k+1}r\leqslant 1$, the H\"older inequality and \eqref{equ:H1} give, on the dyadic shell $\{2^kr\leqslant|y|<2^{k+1}r\}$,
  \begin{equation*}
    \begin{aligned}
      \int_{\{2^kr\leqslant|y|<2^{k+1}r\}}\frac{u(y)^{p-1}}{|y-rz_i|^{N+ps}}\diopy
      &\leqslant\left(2^{k-1}r\right)^{-N-ps}\left|B_{2^{k+1}r}\right|^{1-\frac{p-1}{\pst}}m_0\left(2^{k+1}r\right)^{\frac{p-1}{\pst}}\\
      &\leqslant C\left(2^kr\right)^{-ps-(p-1)(r_H-\delta)}.
    \end{aligned}
  \end{equation*}
  Summing over these $k$, we obtain
  \begin{equation*}
    \int_{\{4r\leqslant|y|<\frac{1}{2}\}}\frac{u(y)^{p-1}}{|y-rz_i|^{N+ps}}\diopy
    \leqslant C\sum_{k\geqslant 2}\left(2^kr\right)^{-ps-(p-1)(r_H-\delta)}
    \leqslant Cr^{-ps-(p-1)(r_H-\delta)},
  \end{equation*}
  because $ps+(p-1)(r_H-\delta)>0$ by \eqref{equ:deltaone}. On the third set we have $|y-rz_i|\geqslant\frac{|y|}{2}$, and the integral of $u^{p-1}|y|^{-N-ps}$ over $\{|y|\geqslant\frac{1}{2}\}$ is finite by the H\"older inequality, as in \eqref{equ:tailfinite}. Multiplying the three contributions by $(\frac{r}{4})^{ps}$, we obtain
  \begin{equation*}
    \Tail\left(u;rz_i,\frac{r}{4}\right)^{p-1}\leqslant C\left(r^{-N}r^{N-(p-1)(r_H-\delta)}+r^{-(p-1)(r_H-\delta)}+r^{ps}\right)\leqslant Cr^{-(p-1)(r_H-\delta)},
  \end{equation*}
  since $r_H-\delta>0$. Therefore
  \begin{equation*}
    \Tail\left(u;rz_i,\frac{r}{4}\right)\leqslant Cr^{-r_H+\delta}\quad\text{for any } r\in(0,r_1).
  \end{equation*}
  By \eqref{equ:localbound}, $u\leqslant Cr^{-r_H+\delta}$ a.e. in $B(rz_i,\frac{r}{8})$ for any $i$, and hence a.e. in $\{r\leqslant|x|\leqslant 2r\}$. This gives the first estimate in \eqref{equ:coarse} with $\rho_1:=r_1$.

  \emph{At infinity.} By \eqref{equ:H1}, there exists $r_2\in(8,\infty)$ such that $m_\infty(\frac{3r}{4})\leqslant\ve_\ast$ for $r\in(r_2,\infty)$. For such $r$, the first term in \eqref{equ:localbound} is at most $Cr^{-r_H}m_\infty(\frac{3r}{4})^{\frac{1}{\pst}}\leqslant Cr^{-r_H-\delta}$. For the tail, we split $\mbRN\setminus B(rz_i,\frac{r}{4})$ into $\{|y|<\frac{r}{2}\}$, $\{\frac{r}{2}\leqslant|y|<4r\}$, and $\{|y|\geqslant 4r\}$. On the first set, $|y-rz_i|\geqslant\frac{r}{2}$, and, by the H\"older inequality with \eqref{equ:H1} on the dyadic shells $\{2^k\leqslant|y|<2^{k+1}\}$,
  \begin{equation*}
    \int_{B_{\frac{r}{2}}}u^{p-1}\leqslant\int_{B_1}u^{p-1}+C\sum_{0\leqslant k\leqslant\log_2 r}2^{k[N-(p-1)(r_H+\delta)]}\leqslant Cr^{N-(p-1)(r_H+\delta)},
  \end{equation*}
  because $N-(p-1)(r_H+\delta)>0$ by \eqref{equ:deltaone}. On the second set the kernel is at most $(\frac{r}{4})^{-N-ps}$ and
  \begin{equation*}
    \int_{\{\frac{r}{2}\leqslant|y|<4r\}}u^{p-1}\leqslant Cr^{N(1-\frac{p-1}{\pst})}m_\infty\left(\frac{r}{2}\right)^{\frac{p-1}{\pst}}\leqslant Cr^{N-(p-1)(r_H+\delta)}.
  \end{equation*} On the third set, $|y-rz_i|\geqslant\frac{|y|}{2}$, and, on the shells $\{2^kr\leqslant|y|<2^{k+1}r\}$ with $k\geqslant 2$, the H\"older inequality and \eqref{equ:H1} give the bound $C(2^kr)^{-N-ps}(2^kr)^{N-(p-1)(r_H+\delta)}$; summing over $k$ yields the bound $Cr^{-ps-(p-1)(r_H+\delta)}$. Multiplying by $(\frac{r}{4})^{ps}$ we find $\Tail(u;rz_i,\frac{r}{4})\leqslant Cr^{-r_H-\delta}$, and, as before, $u\leqslant Cr^{-r_H-\delta}$ a.e. in $\{r\leqslant|x|\leqslant 2r\}$ for $r\in(r_2,\infty)$. This gives the second estimate in \eqref{equ:coarse} with $R_1:=r_2$, and completes the proof of Lemma \ref{lem:coarse}.
\end{proof}

\section{Barriers}\label{sec:barriers}
In this section we construct the barriers that will be compared with $u$.

\subsection{Two-power profiles}
Let $\mu_1$ and $\mu_2$ be the two roots given in \eqref{equ:tworoots}, for $i\in\{1,2\}$ and $\vartheta\in(0,\infty)$ we consider
\begin{equation*}
  g_{i,\vartheta}(x):=|x|^{-\mu_i}-|x|^{-\vartheta},\qquad x\neq 0.
\end{equation*}
We say that $(i,\vartheta)$ is \emph{admissible} if
\begin{equation}\label{equ:admissiblepair}
  \text{either} \ \  i=1\ \text{and}\ 0<\vartheta<\mu_1,  \ \
  \text{or} \ \ i=2,\ \vartheta>\mu_2,\ \vartheta+\mu_2(p-2)<N,\ \text{and}\ \vartheta(p-1)<N.
\end{equation}
By \eqref{equ:exponentfacts}, the set of $\vartheta$ satisfying the second alternative is a nonempty interval with left endpoint $\mu_2$. For an admissible pair $(i,\vartheta)$ we define
\begin{equation}\label{equ:Lambdatilde}
  \wtilde\Lambda_i(\vartheta):=2\int_1^\infty K(\sigma)\left(1-\sigma^{-\mu_i}\right)^{p-2}\left[\sigma^{N-1}\left(1-\sigma^{-\vartheta}\right)-\sigma^{ps-1+\mu_i(p-2)}\left(\sigma^{\vartheta}-1\right)\right]\diopsigma.
\end{equation}

\begin{lemma}\label{lem:Lambdatilde}
  Let $(i,\vartheta)$ be admissible. Then the integral in \eqref{equ:Lambdatilde} converges absolutely, and
  \begin{equation}\label{equ:Lambdatildeless}
    \wtilde\Lambda_i(\vartheta)<\mu.
  \end{equation}
\end{lemma}
\begin{proof}
  If we replace $\vartheta$ by $\mu_i$ in the bracket of \eqref{equ:Lambdatilde}, then, since $\sigma^{\mu_i(p-2)}(\sigma^{\mu_i}-1)=\sigma^{\mu_i(p-1)}(1-\sigma^{-\mu_i})$ and $\sigma^{N-1}(1-\sigma^{-\mu_i})^{p-1}=\sigma^{N-1-\mu_i(p-1)}(\sigma^{\mu_i}-1)^{p-1}$, the integrand becomes that of \eqref{equ:Lambdaformula} with $r=\mu_i$, whose integral equals $\Lambda(\mu_i)=\mu$. Subtracting, and using $\sigma^{-\mu_i}-\sigma^{-\vartheta}=\sigma^{-\mu_i-\vartheta}(\sigma^\vartheta-\sigma^{\mu_i})$, we obtain
  \begin{equation}\label{equ:Lambdatildedifference}
  \begin{split}
    \wtilde\Lambda_i(\vartheta)-\mu={}&2\int_1^\infty K(\sigma)\left(1-\sigma^{-\mu_i}\right)^{p-2}\left(\sigma^{\vartheta}-\sigma^{\mu_i}\right)\\
    &\times\left(\sigma^{N-1-\mu_i-\vartheta}-\sigma^{ps-1+\mu_i(p-2)}\right)\diopsigma,
  \end{split}
  \end{equation}
  provided that all integrals involved converge absolutely. Near $\sigma=1$, the bracket in \eqref{equ:Lambdatilde} is $O((\sigma-1)^2)$, because both of its terms equal $0$ at $\sigma=1$ and have the same derivative $\vartheta$ there, and the product of the last two factors in \eqref{equ:Lambdatildedifference} is also $O((\sigma-1)^2)$. Since $(1-\sigma^{-\mu_i})^{p-2}\leqslant C(\sigma-1)^{p-2}$ on $(1,2]$, the factors multiplying $K(\sigma)$ in both integrands are $O((\sigma-1)^{p})$. Thus, by \eqref{equ:Kbound}, both integrands are bounded in absolute value by $C(\sigma-1)^{p-1-ps}$ on $(1,2]$. This bound is integrable because $p-ps>0$. For $\sigma>2$, by \eqref{equ:Kbound}, both integrands are bounded by
  \begin{equation*}
    C\sigma^{-N-ps}\left(\sigma^{N-1}+\sigma^{ps-1+\mu_i(p-2)+\max\{\vartheta,\mu_i\}}\right)\leqslant C\left(\sigma^{-1-ps}+\sigma^{\max\{\vartheta,\mu_i\}+\mu_i(p-2)-N-1}\right),
  \end{equation*}
  which is integrable because $\max\{\vartheta,\mu_i\}+\mu_i(p-2)<N$: for $i=1$ this follows from $\mu_1(p-1)<N$, and for $i=2$ it is part of \eqref{equ:admissiblepair}.

  It remains to determine the sign of the integrand in \eqref{equ:Lambdatildedifference} for $\sigma>1$. The exponents of the last factor satisfy
  \begin{equation*}
    \left(N-1-\mu_i-\vartheta\right)-\left(ps-1+\mu_i(p-2)\right)=N-ps-\mu_i(p-1)-\vartheta.
  \end{equation*}
  If $i=1$, then $\vartheta<\mu_1$ gives $N-ps-\mu_1(p-1)-\vartheta>N-ps-p\mu_1>0$, since $\mu_1<r_H$; thus the last factor is positive, while $\sigma^{\vartheta}-\sigma^{\mu_1}<0$. If $i=2$, then $\vartheta>\mu_2$ gives $N-ps-\mu_2(p-1)-\vartheta<N-ps-p\mu_2<0$, since $\mu_2>r_H$; thus the last factor is negative, while $\sigma^{\vartheta}-\sigma^{\mu_2}>0$. In both cases the integrand in \eqref{equ:Lambdatildedifference} is negative for $\sigma>1$, which proves \eqref{equ:Lambdatildeless}.
\end{proof}

\begin{lemma}\label{lem:twopower}
  Let $(i,\vartheta)$ be admissible. Then there exist $\tau_1\in(0,\frac{1}{2}]$ and $c_1\in(0,\infty)$, depending only on $N$, $s$, $p$, $\mu$, $i$, and $\vartheta$, such that, for any $x\neq 0$ with $|x|^{\mu_i-\vartheta}\leqslant\tau_1$, the following hold:
  \begin{enumerate}
    \item\label{item:twopowershape} $\frac{1}{2}|x|^{-\mu_i}\leqslant g_{i,\vartheta}(x)\leqslant|x|^{-\mu_i}$ and $\nabla g_{i,\vartheta}(x)\neq 0$;
    \item\label{item:twopowerresidual} the limit in \eqref{equ:fractionalplaplacian} defining $(-\Delta_p)^sg_{i,\vartheta}(x)$ exists, and
    \begin{equation}\label{equ:twopowerresidual}
      (-\Delta_p)^sg_{i,\vartheta}(x)-\frac{\mu}{|x|^{ps}}g_{i,\vartheta}(x)^{p-1}\geqslant c_1|x|^{-ps-\mu_i(p-1)+\mu_i-\vartheta}.
    \end{equation}
  \end{enumerate}
\end{lemma}
\begin{proof}
  Write $g:=g_{i,\vartheta}$ and $\tau:=|x|^{\mu_i-\vartheta}$. For $\tau\in(0,1)$ let $f_\tau(t):=t^{-\mu_i}-\tau t^{-\vartheta}$ for $t>0$. Then $g(|x| y)=|x|^{-\mu_i}f_\tau(|y|)$ for $y\neq 0$, and $g(x)=|x|^{-\mu_i}(1-\tau)$; this proves the first part of \ref{item:twopowershape} if $\tau\leqslant\frac{1}{2}$. Moreover, $\nabla g(x)=\frac{x}{|x|}|x|^{-\mu_i-1}(-\mu_i+\tau\vartheta)\neq 0$ if $\tau<\frac{\mu_i}{\vartheta}$.

  The function $f_\tau$ satisfies \eqref{equ:radialintegrability}, because $\max\{\mu_i,\vartheta\}(p-1)<N$, and $f_\tau'(1)=-\mu_i+\tau\vartheta\neq 0$ if $\tau<\frac{\mu_i}{\vartheta}$. Let $F_\tau(y):=f_\tau(|y|)$, so that $g(y)=|x|^{-\mu_i}F_\tau(\frac{y}{|x|})$ for $y\neq 0$. The change of variables $y=|x|y'$ in \eqref{equ:fractionalplaplacian} shows that
  \begin{equation*}
    (-\Delta_p)^sg(x)=|x|^{-ps-\mu_i(p-1)}(-\Delta_p)^sF_\tau\left(\frac{x}{|x|}\right),
  \end{equation*}
  whenever the right-hand side exists. By Lemma \ref{lem:radialpairing} and the rotation invariance, $(-\Delta_p)^sF_\tau(\frac{x}{|x|})$ exists and is given by \eqref{equ:radialpairing} with $f=f_\tau$. Moreover $\frac{\mu}{|x|^{ps}}g(x)^{p-1}=\mu|x|^{-ps-\mu_i(p-1)}(1-\tau)^{p-1}$. Consequently,
  \begin{equation*}
    |x|^{ps+\mu_i(p-1)}\left[(-\Delta_p)^sg(x)-\frac{\mu}{|x|^{ps}}g(x)^{p-1}\right]=Q(\tau),
  \end{equation*}
  where, with the notation
  \begin{equation*}
    a(\sigma):=1-\sigma^{-\mu_i},\quad b(\sigma):=1-\sigma^{-\vartheta},\quad \hat a(\sigma):=\sigma^{\mu_i}-1,\quad \hat b(\sigma):=\sigma^{\vartheta}-1,\qquad\sigma>1,
  \end{equation*}
  so that $f_\tau(1)-f_\tau(\sigma)=(1-\sigma^{-\mu_i})-\tau(1-\sigma^{-\vartheta})=a-\tau b$ and $f_\tau(1)-f_\tau(\frac{1}{\sigma})=-(\sigma^{\mu_i}-1)+\tau(\sigma^{\vartheta}-1)=-(\hat a-\tau \hat b)$, and since $J_p$ is odd,
  \begin{equation*}
    Q(\tau):=2\int_1^\infty K(\sigma)\left[\sigma^{N-1}J_p\left(a-\tau b\right)-\sigma^{ps-1}J_p\left(\hat a-\tau \hat b\right)\right]\diopsigma-\mu(1-\tau)^{p-1}.
  \end{equation*}
  For $\tau=0$ the integral is the right-hand side of \eqref{equ:radialpairing} for $f(t)=t^{-\mu_i}$, and hence, by Lemma \ref{lem:indexfunction}, $Q(0)=\Lambda(\mu_i)-\mu=0$. We claim that
  \begin{equation}\label{equ:Qderivative}
    \lim_{\tau\to 0^+}\frac{Q(\tau)}{\tau}=(p-1)\left(\mu-\wtilde\Lambda_i(\vartheta)\right).
  \end{equation}
  By Lemma \ref{lem:Lambdatilde}, the right-hand side of \eqref{equ:Qderivative} is positive, so \eqref{equ:Qderivative} yields $\tau_1\in(0,\frac{1}{2}]$ with $\tau_1<\frac{\mu_i}{\vartheta}$ such that $Q(\tau)\geqslant\frac{p-1}{2}(\mu-\wtilde\Lambda_i(\vartheta))\tau$ for $\tau\in(0,\tau_1]$. Multiplying this inequality by $|x|^{-ps-\mu_i(p-1)}$ and using $|x|^{-ps-\mu_i(p-1)}\tau=|x|^{-ps-\mu_i(p-1)+\mu_i-\vartheta}$, we obtain \eqref{equ:twopowerresidual} with $c_1:=\frac{p-1}{2}(\mu-\wtilde\Lambda_i(\vartheta))$.

  It remains to prove \eqref{equ:Qderivative}. Since $\frac{\mu(1-(1-\tau)^{p-1})}{\tau}\to\mu(p-1)$, and since
  \begin{equation*}
    \sigma^{N-1}a^{p-2}b-\sigma^{ps-1}\hat a^{p-2}\hat b=\left(1-\sigma^{-\mu_i}\right)^{p-2}\left[\sigma^{N-1}\left(1-\sigma^{-\vartheta}\right)-\sigma^{ps-1+\mu_i(p-2)}\left(\sigma^{\vartheta}-1\right)\right],
  \end{equation*}
  it suffices to show that
  \begin{equation}\label{equ:Pidifference}
    \lim_{\tau\to 0^+}\int_1^\infty K(\sigma)\frac{\Pi_\tau(\sigma)-\Pi_0(\sigma)}{\tau}\diopsigma=-(p-1)\int_1^\infty K(\sigma)\left[\sigma^{N-1}a^{p-2}b-\sigma^{ps-1}\hat a^{p-2}\hat b\right]\diopsigma,
  \end{equation}
  where $\Pi_\tau:=\sigma^{N-1}J_p(a-\tau b)-\sigma^{ps-1}J_p(\hat a-\tau \hat b)$. The integrands converge pointwise, since $a,\hat a>0$ for $\sigma>1$. We prove \eqref{equ:Pidifference} by the dominated convergence theorem on $(1,2]$ and on $(2,\infty)$ separately.

  \emph{The range $\sigma\in(1,2]$.} There exist constants $0<c\leqslant C$, depending only on $\mu_i$ and $\vartheta$, such that $a,\hat a\in[c(\sigma-1),C(\sigma-1)]$ and $b,\hat b\leqslant C(\sigma-1)$ for $\sigma\in(1,2]$. Let $\tau_2\in(0,\frac{1}{2}]$ be small, depending only on $\mu_i$ and $\vartheta$; for $\tau\in(0,\tau_2]$ and $\lambda\in[0,1]$ the numbers $X_\lambda:=a-\lambda\tau b$ and $Y_\lambda:=\hat a-\lambda\tau \hat b$ lie in $[\frac{c}{2}(\sigma-1),C(\sigma-1)]$. Then
  \begin{equation*}
    \frac{\Pi_\tau-\Pi_0}{\tau}=-(p-1)\int_0^1\left(\sigma^{N-1}X_\lambda^{p-2}b-\sigma^{ps-1}Y_\lambda^{p-2}\hat b\right)\diop\lambda.
  \end{equation*}
  Since $\hat a-a=\sigma^{\mu_i}+\sigma^{-\mu_i}-2=O((\sigma-1)^2)$ and $\hat b-b=O((\sigma-1)^2)$, we have $|X_\lambda-Y_\lambda|\leqslant C(\sigma-1)^2$, and the mean value theorem gives $|X_\lambda^{p-2}-Y_\lambda^{p-2}|\leqslant C(\sigma-1)^{p-3}|X_\lambda-Y_\lambda|\leqslant C(\sigma-1)^{p-1}$. Hence
  \begin{align*}
    \left|\sigma^{N-1}X_\lambda^{p-2}b-\sigma^{ps-1}Y_\lambda^{p-2}\hat b\right|
    &\leqslant|\sigma^{N-1}-\sigma^{ps-1}|X_\lambda^{p-2}b+\sigma^{ps-1}\left(|X_\lambda^{p-2}-Y_\lambda^{p-2}|b+Y_\lambda^{p-2}|b-\hat b|\right)\\
    &\leqslant C(\sigma-1)^{p},
  \end{align*}
  and, by \eqref{equ:Kbound}, $K(\sigma)|\frac{\Pi_\tau-\Pi_0}{\tau}|\leqslant C(\sigma-1)^{p-1-ps}$, which is integrable on $(1,2]$.

  \emph{The range $\sigma\in(2,\infty)$.} Here $a\geqslant 1-2^{-\mu_i}$ and $0<b<1$; decreasing $\tau_2$ if necessary, $a-\lambda\tau b\geqslant\frac{1}{2}(1-2^{-\mu_i})$ for $\lambda\in[0,1]$ and $\tau\in(0,\tau_2]$, and hence $\sigma^{N-1}\frac{|J_p(a-\tau b)-J_p(a)|}{\tau}\leqslant C\sigma^{N-1}$. By \eqref{equ:Kbound}, this is dominated by $C\sigma^{-1-ps}$ after multiplication by $K$. For the second term we distinguish two cases.

  If $i=1$, then $\hat b\leqslant \hat a$ and $\hat a-\lambda\tau \hat b\geqslant\frac{\hat a}{2}$, so the mean value theorem gives
  \begin{equation*}
    \sigma^{ps-1}\frac{|J_p(\hat a-\tau \hat b)-J_p(\hat a)|}{\tau}\leqslant C\sigma^{ps-1}\hat a^{p-2}\hat b\leqslant C\sigma^{ps-1+\mu_1(p-2)+\vartheta},
  \end{equation*}
  and, after multiplication by $K$, this is dominated by $C\sigma^{\mu_1(p-2)+\vartheta-N-1}$, which is integrable because $\mu_1(p-2)+\vartheta<\mu_1(p-1)<N$.

  If $i=2$, let $G_\tau:=\{\sigma>2\colon\tau \hat b(\sigma)\leqslant\frac{\hat a(\sigma)}{2}\}$. On $G_\tau$, as in the case $i=1$, the integrand is dominated by $C\sigma^{\mu_2(p-2)+\vartheta-N-1}$, which is integrable by \eqref{equ:admissiblepair}, and $\mathbf{1}_{G_\tau}\to 1$ pointwise as $\tau\to 0^+$. On $(2,\infty)\setminus G_\tau$ we have $\hat a<2\tau \hat b$, hence $|J_p(\hat a-\tau \hat b)-J_p(\hat a)|\leqslant C(\tau \hat b)^{p-1}\leqslant C\tau^{p-1}\sigma^{\vartheta(p-1)}$, and $\sigma\geqslant c\tau^{-\frac{1}{\vartheta-\mu_2}}$ for some $c=c(\mu_2,\vartheta)>0$. Therefore
  \begin{multline*}
    \int_{(2,\infty)\setminus G_\tau}K(\sigma)\sigma^{ps-1}\frac{|J_p(\hat a-\tau \hat b)-J_p(\hat a)|}{\tau}\diopsigma\\
    \leqslant C\tau^{p-2}\int_{c\tau^{-\frac{1}{\vartheta-\mu_2}}}^\infty\sigma^{\vartheta(p-1)-N-1}\diopsigma
    =C\tau^{p-2+\frac{N-\vartheta(p-1)}{\vartheta-\mu_2}},
  \end{multline*}
  where we used $\vartheta(p-1)<N$. The last exponent is positive, because
  \begin{equation*}
    N-\vartheta(p-1)-(2-p)(\vartheta-\mu_2)=N-\vartheta-\mu_2(p-2)>0
  \end{equation*}
  by \eqref{equ:admissiblepair}. Hence this contribution tends to $0$ as $\tau\to 0^+$. Moreover, the limit integrand is integrable on $(2,\infty)$ and $\mathbf{1}_{(2,\infty)\setminus G_\tau}\to 0$ pointwise, so its integral over $(2,\infty)\setminus G_\tau$ also tends to $0$.

  Combining the two ranges, the dominated convergence theorem gives \eqref{equ:Pidifference}, and hence \eqref{equ:Qderivative}. This finishes the proof of Lemma \ref{lem:twopower}.
\end{proof}

\subsection{Barriers for the lower bounds}
Fix, once and for all, radial functions $\chi,\chi_1,\zeta\in C^\infty(\mbRN)$ such that $0\leqslant\chi\leqslant 1$, $\chi=1$ on $B_1$, $\supp\chi\subset B_2$; $0\leqslant\chi_1\leqslant 1$, $\chi_1=0$ on $B_{\frac{1}{4}}$, $\chi_1=1$ outside $B_2$; and $\zeta\geqslant 0$, $\zeta\not\equiv 0$, $\supp\zeta\subset\{\frac{1}{2}<|x|<1\}$. Let $U(x):=|x|^{-\mu_1}$.

\begin{lemma}\label{lem:lowerbarrierorigin}
  Let $\Psi_0:=(U+1)\chi$. Then $\Psi_0\in\mcD$, $\Psi_0\in C^\infty(B_1\setminus\{0\})$ with $\nabla\Psi_0\neq 0$ in $B_1\setminus\{0\}$, and there exists $\eta_0\in(0,\frac{1}{4}]$, depending only on $N$, $s$, $p$, and $\mu$, such that, for any $x\in B_{\eta_0}\setminus\{0\}$,
  \begin{equation}\label{equ:lowerbarrierorigin}
    (-\Delta_p)^s\Psi_0(x)<\frac{\mu}{|x|^{ps}}\Psi_0(x)^{p-1}.
  \end{equation}
\end{lemma}
\begin{proof}
  By Lemma \ref{lem:powerenergy}\ref{item:powerenergyone} and $\mu_1<r_H$, $\Psi_0\in\mcD$. The smoothness and $\nabla\Psi_0=\nabla U\neq 0$ in $B_1\setminus\{0\}$ are clear. Let $x\in B_{\frac{1}{4}}\setminus\{0\}$ with $U(x)\geqslant 4$. Since adding a constant does not change the integrand in \eqref{equ:fractionalplaplacian}, \eqref{equ:purepower} gives $(-\Delta_p)^s(U+1)(x)=\mu|x|^{-ps}U(x)^{p-1}$. As $\Psi_0=U+1$ on $B_1$,
  \begin{equation*}
    (-\Delta_p)^s\Psi_0(x)-\frac{\mu U(x)^{p-1}}{|x|^{ps}}=2\int_{\{|y|>1\}}\frac{J_p\left(U(x)+1-\Psi_0(y)\right)-J_p\left(U(x)-U(y)\right)}{|x-y|^{N+ps}}\diopy.
  \end{equation*}
  For $|y|>1$ both arguments of $J_p$ lie in $[U(x)-1,U(x)+1]\subset[\frac{U(x)}{2},2U(x)]$ and differ by $|U(y)+1-\Psi_0(y)|\leqslant 2$. By the mean value theorem the difference of the two values of $J_p$ is at most $C_pU(x)^{p-2}$, where $C_p$ depends only on $p$. Since $|x-y|\geqslant\frac{3|y|}{4}$ for $|y|>1$, we obtain $(-\Delta_p)^s\Psi_0(x)\leqslant\mu|x|^{-ps}U(x)^{p-1}+C_2U(x)^{p-2}$ with $C_2$ depending only on $N$, $s$, and $p$. On the other hand, by the mean value theorem, $(U+1)^{p-1}-U^{p-1}\geqslant(p-1)\min\{1,2^{p-2}\}U^{p-2}$ for $U\geqslant 1$. Therefore
  \begin{equation*}
    (-\Delta_p)^s\Psi_0(x)-\frac{\mu}{|x|^{ps}}\Psi_0(x)^{p-1}\leqslant U(x)^{p-2}\left[C_2-\mu(p-1)\min\{1,2^{p-2}\}|x|^{-ps}\right].
  \end{equation*}
  Choosing $\eta_0\in(0,\frac{1}{4}]$ such that $\eta_0^{-\mu_1}\geqslant 4$ and $\mu(p-1)\min\{1,2^{p-2}\}\eta_0^{-ps}>C_2$, we obtain \eqref{equ:lowerbarrierorigin} for $x\in B_{\eta_0}\setminus\{0\}$.
\end{proof}

\begin{lemma}\label{lem:lowerbarrierinfinity}
  There exist $M_0\in\mbN$ and $R_0\in[4,\infty)$, depending only on $N$, $s$, $p$, and $\mu$, such that the function $\Psi_\infty:=\chi_1|\cdot|^{-\mu_2}+M_0\zeta$ satisfies $\Psi_\infty\in\mcD\cap C^\infty(\mbRN)$, $\Psi_\infty=|\cdot|^{-\mu_2}$ and $\nabla\Psi_\infty\neq 0$ outside $\ol{B_2}$, and, for any $x\in\mbRN\setminus B_{R_0}$,
  \begin{equation}\label{equ:lowerbarrierinfinity}
    (-\Delta_p)^s\Psi_\infty(x)<\frac{\mu}{|x|^{ps}}\Psi_\infty(x)^{p-1}.
  \end{equation}
\end{lemma}
\begin{proof}
  For $M\in\mbN$ let $\Psi_M:=\chi_1|\cdot|^{-\mu_2}+M\zeta$. Since $\mu_2(p-1)<N$, the function $|y|^{-\mu_2(p-1)}$ is integrable on $B_2$, while $\int_{B_2}\Psi_M^{p-1}\geqslant M^{p-1}\int_{B_2}\zeta^{p-1}\to\infty$ as $M\to\infty$. Let $M_0$ be the smallest $M\in\mbN$ such that
  \begin{equation}\label{equ:Mzero}
    \int_{B_2}\left(|y|^{-\mu_2(p-1)}-\Psi_M(y)^{p-1}\right)\diopy<-1,
  \end{equation}
  and let $\Psi_\infty:=\Psi_{M_0}$. Then $\Psi_\infty$ is smooth, equals $0$ on $B_{\frac{1}{4}}$, and equals $|\cdot|^{-\mu_2}$ outside $B_2$; by Lemma \ref{lem:powerenergy}\ref{item:powerenergyone} and $\mu_2>r_H$, $\Psi_\infty\in\mcD$. Let $G(y):=|y|^{-\mu_2}$. For $|x|\geqslant 4$, \eqref{equ:purepower} gives $(-\Delta_p)^sG(x)=\mu|x|^{-ps}G(x)^{p-1}$, and $\Psi_\infty=G$ outside $B_2$, so that
  \begin{align*}
    &|x|^{N+ps}\left[(-\Delta_p)^s\Psi_\infty(x)-\frac{\mu}{|x|^{ps}}\Psi_\infty(x)^{p-1}\right]\\
    &\quad=2\int_{B_2}\frac{|x|^{N+ps}}{|x-y|^{N+ps}}\left[J_p\left(G(x)-\Psi_\infty(y)\right)-J_p\left(G(x)-G(y)\right)\right]\diopy.
  \end{align*}
  For $|x|\geqslant 4$ and $y\in B_2$, the kernel ratio is at most $2^{N+ps}$ and the absolute value of the bracket is at most $C(1+\Psi_\infty(y)^{p-1}+G(y)^{p-1})$, which is integrable on $B_2$. As $|x|\to\infty$, $G(x)\to 0$ and the kernel ratio tends to $1$ uniformly in $y\in B_2$. Moreover, for $y\in B_2$ the bracket tends to $J_p(-\Psi_\infty(y))-J_p(-G(y))=G(y)^{p-1}-\Psi_\infty(y)^{p-1}$ as $|x|\to\infty$. By the dominated convergence theorem, the right-hand side tends to $2\int_{B_2}(G^{p-1}-\Psi_\infty^{p-1})$ as $|x|\to\infty$, which is less than $-2$ by \eqref{equ:Mzero}. Since $\Psi_\infty$ and $G$ are radial, the left-hand side is a radial function of $x$. Hence there exists $R_0\in[4,\infty)$ such that the left-hand side is less than $-1$ for $|x|\geqslant R_0$, which gives \eqref{equ:lowerbarrierinfinity}.
\end{proof}

\section{Bounds near the origin and at infinity and proofs of the main results}\label{sec:bounds}

With the previous preparations, we are able to prove the lower and upper bounds for positive weak solution of \eqref{equ:EL} near the origin and at infinity, which yields Theorem \ref{thm:twosided}.

\subsection{Lower bounds}
\begin{proposition}\label{prop:lower}
  Let $\eta_0$ and $R_0$ be as in Lemmas \ref{lem:lowerbarrierorigin} and \ref{lem:lowerbarrierinfinity} and assume that $u\in\mathcal{D}^{s,p}(\mathbb{R}^{N})$ is a positive weak solution of \eqref{equ:EL}. Then there exists $c\in(0,\infty)$, depending on $u$, such that
  \begin{equation}\label{equ:lower}
    u(x)\geqslant c|x|^{-\mu_1}\quad\text{for a.e. } x\in B_{\eta_0}\setminus\{0\},
    \qquad
    u(x)\geqslant c|x|^{-\mu_2}\quad\text{for a.e. } x\in\mbRN\setminus B_{2R_0}.
  \end{equation}
\end{proposition}
\begin{proof}
  \emph{Near the origin.} Let $V:=\mu|x|^{-ps}$, which is an admissible potential. By Corollary \ref{cor:annulus}, $m_1:=\essinf_{B_2}u>0$. Let $c_0:=\frac{m_1}{2((\frac{\eta_0}{2})^{-\mu_1}+1)}$ and $f_1:=c_0\Psi_0$, $f_2:=u$. By Lemmas \ref{lem:lowerbarrierorigin} and \ref{lem:pointwiseweak}, and by the $(p-1)$-homogeneity of \eqref{equ:lowerbarrierorigin}, $f_1$ is a weak subsolution of $(-\Delta_p)^sf=Vf^{p-1}$ in $\Omega:=B_{\eta_0}$. As observed after \eqref{equ:Vu}, $f_2$ is a weak supersolution of the same equation in $\Omega$. Since $\Psi_0\leqslant(\frac{\eta_0}{2})^{-\mu_1}+1$ outside $B_{\frac{\eta_0}{2}}$ and $\Psi_0=0$ outside $B_2$, we have $f_1\leqslant\frac{m_1}{2}\leqslant\frac{u}{2}$ a.e. in $B_2\setminus B_{\frac{\eta_0}{2}}$ and $f_1=0<f_2$ outside $B_2$. Hence $f_1\leqslant f_2$ a.e. in $\mbRN\setminus\Omega_{\frac{\eta_0}{2}}=\mbRN\setminus B_{\frac{\eta_0}{2}}$, and $\{f_1<f_2\}\setminus\Omega$ has positive measure. By Lemma \ref{lem:comparison}, $u\geqslant c_0\Psi_0\geqslant c_0|x|^{-\mu_1}$ a.e. in $B_{\eta_0}$.

  \emph{At infinity.} Let $\Omega:=\mbRN\setminus\ol{B_{R_0}}$, $m_2:=\essinf_{B_{2R_0}}u>0$, $c_\infty:=\frac{m_2}{2\max\Psi_\infty}$, $f_1:=c_\infty\Psi_\infty$, and $f_2:=u$. By Lemmas \ref{lem:lowerbarrierinfinity} and \ref{lem:pointwiseweak}, $f_1$ is a weak subsolution of $(-\Delta_p)^sf=\mu|x|^{-ps}f^{p-1}$ in $\Omega$, and $f_2$ is a weak supersolution of the same equation in $\Omega$. Moreover, $f_1\leqslant\frac{u}{2}$ a.e. in $\ol{B_{2R_0}}=\mbRN\setminus\Omega_{R_0}$, and $f_1=0<f_2$ in $B_{\frac{1}{4}}\subset\mbRN\setminus\Omega$. By Lemma \ref{lem:comparison}, $u\geqslant c_\infty\Psi_\infty=c_\infty|x|^{-\mu_2}$ a.e. in $\mbRN\setminus B_{2R_0}$. Taking $c:=\min\{c_0,c_\infty\}$ completes the proof of Proposition \ref{prop:lower}.
\end{proof}

\subsection{Upper bound near the origin}
Let $\delta_1$ be as in \eqref{equ:deltaone}, and let
\begin{equation}\label{equ:nuorigin}
  \nu_0:=\frac{1}{2}\min\left\{\mu_1,\ ps,\ \delta_1(\pst-p)\right\},\qquad\vartheta_0:=\mu_1-\nu_0.
\end{equation}
Then $(1,\vartheta_0)$ is admissible in the sense of \eqref{equ:admissiblepair}, and $g_{1,\vartheta_0}(x)=|x|^{-\mu_1}(1-|x|^{\nu_0})$.

\begin{proposition}\label{prop:upperorigin}
Assume that $u\in\mathcal{D}^{s,p}(\mathbb{R}^{N})$ is a positive weak solution of \eqref{equ:EL}, then there exist $\rho\in(0,1)$ and $C\in(0,\infty)$, depending on $u$, such that
  \begin{equation}\label{equ:upperorigin}
    u(x)\leqslant C|x|^{-\mu_1}\quad\text{for a.e. } x\in B_\rho\setminus\{0\}.
  \end{equation}
\end{proposition}
\begin{proof}
  Write $g:=g_{1,\vartheta_0}$. Let $\tau_1$ and $c_1$ be given by Lemma \ref{lem:twopower} for the pair $(1,\vartheta_0)$, so that the conclusions of Lemma \ref{lem:twopower} hold for $0<|x|\leqslant\tau_1^{\frac{1}{\nu_0}}$. Let $\rho_1$ be as in Lemma \ref{lem:coarse}, and let
  \begin{equation*}
    \rho_\ast:=\frac{1}{4}\min\left\{\tau_1^{\frac{1}{\nu_0}},\ \rho_1\right\}.
  \end{equation*}
  Then $g>0$ in $B_{2\rho_\ast}\setminus\{0\}$. With $\chi$ as in Section \ref{sec:barriers}, let $\chi_\ast:=\chi(\frac{\cdot}{\rho_\ast})$, so that $\chi_\ast=1$ on $B_{\rho_\ast}$ and $\supp\chi_\ast\subset B_{2\rho_\ast}$. For $A\in[1,\infty)$ we define
  \begin{equation*}
    \ol u_A:=A\chi_\ast g+(1-\chi_\ast)u.
  \end{equation*}
  By Lemma \ref{lem:powerenergy}, $\ol u_A\in\mcD$, and $\ol u_A>0$ a.e. in $\mbRN$. Indeed, $\chi_\ast g$ satisfies the assumption of Lemma \ref{lem:powerenergy}\ref{item:powerenergyone} with $a_0=\mu_1<r_H$ and any $a_\infty>r_H$, since it vanishes outside $B_{2\rho_\ast}$; and $(1-\chi_\ast)u\in\mcD$ by Lemma \ref{lem:powerenergy}\ref{item:powerenergytwo}. Positivity follows from $g>0$ in $B_{2\rho_\ast}\setminus\{0\}$ and $u>0$ a.e.

  \emph{Step 1: $\ol u_A$ is a pointwise supersolution near the origin.} Let $x\in B_{\frac{\rho_\ast}{2}}\setminus\{0\}$. Since $\ol u_A=Ag$ in $B_{\rho_\ast}$ and $|x-y|\geqslant\frac{|y|}{2}$ for $|y|\geqslant\rho_\ast$, the limit defining $(-\Delta_p)^s\ol u_A(x)$ exists by Lemma \ref{lem:twopower}, and
  \begin{equation*}
    (-\Delta_p)^s\ol u_A(x)=A^{p-1}(-\Delta_p)^sg(x)+E(x),
  \end{equation*}
  where
  \begin{equation*}
    E(x):=2\int_{\{|y|\geqslant\rho_\ast\}}\frac{J_p\left(Ag(x)-\ol u_A(y)\right)-J_p\left(Ag(x)-Ag(y)\right)}{|x-y|^{N+ps}}\diopy.
  \end{equation*}
  For $|y|\geqslant\rho_\ast$ we have $|\ol u_A(y)-Ag(y)|=(1-\chi_\ast(y))|u(y)-Ag(y)|\leqslant A(u(y)+|g(y)|)$ and $|g(y)|\leqslant 2\rho_\ast^{-\mu_1}$. If $p\in(1,2)$, then $|J_p(t_1)-J_p(t_2)|\leqslant 2^{2-p}|t_1-t_2|^{p-1}$, and hence
  \begin{equation*}
    |E(x)|\leqslant CA^{p-1}\int_{\{|y|\geqslant\rho_\ast\}}\frac{u(y)^{p-1}+|g(y)|^{p-1}}{|y|^{N+ps}}\diopy\leqslant C_3A^{p-1},
  \end{equation*}
  where $C_3$ is finite by \eqref{equ:tailfinite} and depends on $u$. If $p\geqslant 2$, then $|J_p(t_1)-J_p(t_2)|\leqslant(p-1)(|t_1|+|t_2|)^{p-2}|t_1-t_2|$, and, since $g(x)\leqslant|x|^{-\mu_1}$,
  \begin{multline*}
    |E(x)|\leqslant CA^{p-1}\int_{\{|y|\geqslant\rho_\ast\}}\frac{\left(|x|^{-\mu_1}+u(y)+|g(y)|\right)^{p-2}\left(u(y)+|g(y)|\right)}{|y|^{N+ps}}\diopy\\
    \leqslant C_3A^{p-1}\left(1+|x|^{-\mu_1(p-2)}\right),
  \end{multline*}
  where we used that $\int_{\{|y|\geqslant\rho_\ast\}}(u+|g|)^q|y|^{-N-ps}\diopy<\infty$ for $q\in\{1,p-1\}$, by the H\"older inequality and $u\in L^{\pst}(\mbRN)$. Next, by Lemma \ref{lem:twopower}\ref{item:twopowerresidual},
  \begin{equation*}
    A^{p-1}\left[(-\Delta_p)^sg(x)-\frac{\mu}{|x|^{ps}}g(x)^{p-1}\right]\geqslant c_1A^{p-1}|x|^{-ps-\mu_1(p-1)+\nu_0},
  \end{equation*}
  and, by \eqref{equ:coarsecritical} and $g\leqslant|x|^{-\mu_1}$,
  \begin{equation*}
    u(x)^{\pst-p}\left(Ag(x)\right)^{p-1}\leqslant C_4A^{p-1}|x|^{-ps-\mu_1(p-1)+\delta_1(\pst-p)}\quad\text{for a.e. } x\in B_{\frac{\rho_\ast}{2}},
  \end{equation*}
  with $C_4$ depending on $u$. Collecting these estimates, for a.e. $x\in B_{\frac{\rho_\ast}{2}}\setminus\{0\}$,
  \begin{align*}
    &(-\Delta_p)^s\ol u_A(x)-V_u(x)\ol u_A(x)^{p-1}\\
    &\quad\geqslant A^{p-1}|x|^{-ps-\mu_1(p-1)}\left[c_1|x|^{\nu_0}-C_3|x|^{ps+\mu_1(p-1)}\left(1+|x|^{-\mu_1(p-2)_+}\right)-C_4|x|^{\delta_1(\pst-p)}\right].
  \end{align*}
  Since $\nu_0<ps$, $\nu_0<ps+\mu_1$, and $\nu_0<\delta_1(\pst-p)$ by \eqref{equ:nuorigin}, and since $ps+\mu_1(p-1)-\mu_1(p-2)_+\geqslant\min\{ps+\mu_1(p-1),ps+\mu_1\}$, there exists $\rho\in(0,\frac{\rho_\ast}{2})$, depending on $u$ but not on $A$, such that the bracket is positive for $0<|x|<\rho$. Since $\ol u_A=Ag$ is smooth with nonzero gradient in $B_{\rho_\ast}\setminus\{0\}$, Lemma \ref{lem:pointwiseweak} shows that $\ol u_A$ is a weak supersolution of $(-\Delta_p)^sf=V_uf^{p-1}$ in $B_\rho$ for any $A\in[1,\infty)$.

  \emph{Step 2: ordering outside $B_{\frac{\rho}{2}}$ and comparison.} By Corollary \ref{cor:annulus} and the continuity and positivity of $g$ in $B_{2\rho_\ast}\setminus\{0\}$,
  \begin{equation*}
    M_u:=\esssup_{\{\frac{\rho}{2}<|x|<2\rho_\ast\}}u<\infty\qquad\text{and}\qquad m_g:=\min_{\{\frac{\rho}{2}\leqslant|x|\leqslant 2\rho_\ast\}}g>0.
  \end{equation*} Let $A:=\max\{1,\frac{2M_u}{m_g}\}$. Then $Ag\geqslant 2u$ a.e. in $\{\frac{\rho}{2}<|x|<2\rho_\ast\}$, and hence $\ol u_A=\chi_\ast Ag+(1-\chi_\ast)u\geqslant u$ a.e. in this annulus, while $\ol u_A=u$ outside $B_{2\rho_\ast}$. Thus $u\leqslant\ol u_A$ a.e. in $\mbRN\setminus B_{\frac{\rho}{2}}$, and $\ol u_A=Ag\geqslant 2u>u$ a.e. in $\{\rho<|x|<\rho_\ast\}$. By Lemma \ref{lem:comparison} with $\Omega:=B_\rho$, $V:=V_u$, $f_1:=u$, and $f_2:=\ol u_A$, we conclude that $u\leqslant\ol u_A=Ag\leqslant A|x|^{-\mu_1}$ a.e. in $B_\rho$. This proves \eqref{equ:upperorigin} with $C:=A$.
\end{proof}

\subsection{Upper bound at infinity}
Let $\delta_1$ be as in \eqref{equ:deltaone}, and let
\begin{equation}\label{equ:nuinfinity}
  \nu_\infty:=\frac{1}{2p}\min\left\{N-\mu_2(p-1),\ \delta_1(\pst-p)\right\},\qquad\vartheta_\infty:=\mu_2+\nu_\infty.
\end{equation}
By \eqref{equ:exponentfacts}, $\nu_\infty>0$. Moreover, $\vartheta_\infty+\mu_2(p-2)=\mu_2(p-1)+\nu_\infty<N$ and $\vartheta_\infty(p-1)=\mu_2(p-1)+(p-1)\nu_\infty<N$, so that $(2,\vartheta_\infty)$ is admissible in the sense of \eqref{equ:admissiblepair}, and $g_{2,\vartheta_\infty}(x)=|x|^{-\mu_2}(1-|x|^{-\nu_\infty})$.

\begin{proposition}\label{prop:upperinfinity}
Assume that $u\in\mathcal{D}^{s,p}(\mathbb{R}^{N})$ is a positive weak solution of \eqref{equ:EL}, then there exist $R\in(1,\infty)$ and $C\in(0,\infty)$, depending on $u$, such that
  \begin{equation}\label{equ:upperinfinity}
    u(x)\leqslant C|x|^{-\mu_2}\quad\text{for a.e. } x\in\mbRN\setminus B_R.
  \end{equation}
\end{proposition}
\begin{proof}
  Write $g:=g_{2,\vartheta_\infty}$. Let $\tau_1$ and $c_1$ be given by Lemma \ref{lem:twopower} for the pair $(2,\vartheta_\infty)$, so that its conclusions hold for $|x|\geqslant\tau_1^{-\frac{1}{\nu_\infty}}$, and let $R_1$ be as in Lemma \ref{lem:coarse}. Let
  \begin{equation*}
    R_2:=\max\left\{\tau_1^{-\frac{1}{\nu_\infty}},\ R_1\right\}.
  \end{equation*}
  Then $g>0$ in $\mbRN\setminus\ol{B_1}$, and in particular outside $B_{R_2}$. Let $\chi_2:=\chi(\frac{\cdot}{R_2})$, so that $\chi_2=1$ on $B_{R_2}$ and $\supp\chi_2\subset B_{2R_2}$. For $A\in[1,\infty)$ we define
  \begin{equation*}
    \ol u_A:=\chi_2u+A(1-\chi_2)g.
  \end{equation*}
  By Lemma \ref{lem:powerenergy}, $\ol u_A\in\mcD$, and $\ol u_A>0$ a.e. in $\mbRN$. Indeed, $\chi_2u\in\mcD$ by Lemma \ref{lem:powerenergy}\ref{item:powerenergytwo}, and $(1-\chi_2)g$ satisfies the assumption of Lemma \ref{lem:powerenergy}\ref{item:powerenergyone} with any $a_0<r_H$ and $a_\infty=\mu_2>r_H$, since it vanishes in $B_{R_2}$ and $|g|+|x||\nabla g|\leqslant C|x|^{-\mu_2}$ for $|x|\geqslant R_2$. Moreover, $\ol u_A=u>0$ a.e.\ in $B_{R_2}$, and in $\mbRN\setminus B_{R_2}$ the function $\ol u_A$ is a convex combination of $u$ and $Ag$, where $g>0$ there; hence $\ol u_A>0$ a.e.\ in $\mbRN$.

  \emph{Step 1: $\ol u_A$ is a pointwise supersolution near infinity.} Let $|x|\geqslant 4R_2$. Since $\ol u_A=Ag$ outside $B_{2R_2}$ and $|x-y|\geqslant\frac{|x|}{2}$ for $y\in B_{2R_2}$, the limit defining $(-\Delta_p)^s\ol u_A(x)$ exists by Lemma \ref{lem:twopower}, and
  \begin{equation*}
    (-\Delta_p)^s\ol u_A(x)=A^{p-1}(-\Delta_p)^sg(x)+E(x),
  \end{equation*}
  where
  \begin{equation*}
    E(x):=2\int_{B_{2R_2}}\frac{J_p\left(Ag(x)-\ol u_A(y)\right)-J_p\left(Ag(x)-Ag(y)\right)}{|x-y|^{N+ps}}\diopy.
  \end{equation*}
  For $y\in B_{2R_2}$ we have $|\ol u_A(y)-Ag(y)|=\chi_2(y)|u(y)-Ag(y)|\leqslant A(u(y)+|g(y)|)$, and $|g(y)|\leqslant|y|^{-\mu_2}+|y|^{-\vartheta_\infty}$. Since $\vartheta_\infty(p-1)<N$, the function $|g|^{p-1}$ is integrable on $B_{2R_2}$, and so is $|g|$ if $p\geqslant 2$; moreover, $u$ and $u^{p-1}$ are integrable on $B_{2R_2}$. We also have $0<g(x)\leqslant 1$. Since $|x-y|\geqslant\frac{|x|}{2}$ for $y\in B_{2R_2}$, we have $|x-y|^{-N-ps}\leqslant 2^{N+ps}|x|^{-N-ps}$ there. If $p\in(1,2)$, then $|J_p(t_1)-J_p(t_2)|\leqslant 2^{2-p}|t_1-t_2|^{p-1}$, and hence
  \begin{equation*}
    |E(x)|\leqslant CA^{p-1}|x|^{-N-ps}\int_{B_{2R_2}}\left(u(y)+|g(y)|\right)^{p-1}\diopy.
  \end{equation*}
  If $p\geqslant 2$, then $|J_p(t_1)-J_p(t_2)|\leqslant(p-1)(|t_1|+|t_2|)^{p-2}|t_1-t_2|$, and both arguments of $J_p$ in $E(x)$ are bounded in modulus by $A(1+u(y)+|g(y)|)$, because $Ag(x)\leqslant A$ and $|\ol u_A(y)|\leqslant A(u(y)+|g(y)|)$; hence
  \begin{equation*}
    |E(x)|\leqslant CA^{p-1}|x|^{-N-ps}\int_{B_{2R_2}}\left(1+u(y)+|g(y)|\right)^{p-2}\left(u(y)+|g(y)|\right)\diopy .
  \end{equation*}
  In both cases the last integral is finite and depends on $u$, and we obtain
  \begin{equation*}
    |E(x)|\leqslant C_5A^{p-1}|x|^{-N-ps}\quad\text{for any } |x|\geqslant 4R_2,
  \end{equation*}
  with $C_5$ depending on $u$. By Lemma \ref{lem:twopower},
  \begin{equation*}
    A^{p-1}\left[(-\Delta_p)^sg(x)-\frac{\mu}{|x|^{ps}}g(x)^{p-1}\right]\geqslant c_1A^{p-1}|x|^{-ps-\mu_2(p-1)-\nu_\infty},
  \end{equation*}
  and, by \eqref{equ:coarsecritical} and $g\leqslant|x|^{-\mu_2}$,
  \begin{equation*}
    u(x)^{\pst-p}\left(Ag(x)\right)^{p-1}\leqslant C_6A^{p-1}|x|^{-ps-\mu_2(p-1)-\delta_1(\pst-p)}\quad\text{for a.e. } |x|\geqslant 4R_2,
  \end{equation*}
  with $C_6$ depending on $u$. Therefore, for a.e. $x$ with $|x|\geqslant 4R_2$,
  \begin{align*}
    &(-\Delta_p)^s\ol u_A(x)-V_u(x)\ol u_A(x)^{p-1}\\
    &\quad\geqslant A^{p-1}|x|^{-ps-\mu_2(p-1)}\left[c_1|x|^{-\nu_\infty}-C_5|x|^{-(N-\mu_2(p-1))}-C_6|x|^{-\delta_1(\pst-p)}\right].
  \end{align*}
  Since $\nu_\infty<N-\mu_2(p-1)$ and $\nu_\infty<\delta_1(\pst-p)$ by \eqref{equ:nuinfinity}, there exists $R_3\in[4R_2,\infty)$, depending on $u$ but not on $A$, such that the bracket is positive for $|x|>R_3$. As $\ol u_A=Ag$ is smooth with nonzero gradient outside $B_{2R_2}$, Lemma \ref{lem:pointwiseweak} shows that $\ol u_A$ is a weak supersolution of $(-\Delta_p)^sf=V_uf^{p-1}$ in $\Omega:=\mbRN\setminus\ol{B_{R_3}}$ for any $A\in[1,\infty)$.

  \emph{Step 2: ordering in $B_{2R_3}$ and comparison.} By Corollary \ref{cor:annulus} and the continuity and positivity of $g$ outside $\ol{B_1}$,
  \begin{equation*}
    M_u:=\esssup_{\{R_2<|x|<2R_3\}}u<\infty\qquad\text{and}\qquad m_g:=\min_{\{R_2\leqslant|x|\leqslant 2R_3\}}g>0.
  \end{equation*} Let $A:=\max\{1,\frac{2M_u}{m_g}\}$. Then $Ag\geqslant 2u$ a.e. in $\{R_2<|x|<2R_3\}$, so that $\ol u_A=\chi_2u+(1-\chi_2)Ag\geqslant u$ a.e. there, while $\ol u_A=u$ in $B_{R_2}$. Hence $u\leqslant\ol u_A$ a.e. in $\ol{B_{2R_3}}=\mbRN\setminus\Omega_{R_3}$, and $\ol u_A=Ag\geqslant 2u>u$ a.e. in $\{2R_2<|x|<R_3\}\subset\mbRN\setminus\Omega$. By Lemma \ref{lem:comparison} with $V:=V_u$, $f_1:=u$, and $f_2:=\ol u_A$, we conclude that $u\leqslant\ol u_A$ a.e. in $\mbRN$. Since $\ol u_A=Ag\leqslant A|x|^{-\mu_2}$ outside $B_{2R_2}$, this proves \eqref{equ:upperinfinity} with $R:=2R_2$ and $C:=A$.
\end{proof}

\subsection{\texorpdfstring{Proofs of Theorem \ref{thm:twosided} and Proposition \ref{prop:scaling}}{Proofs of Theorem 1.1 and Proposition 1.2}}\label{sec:proofs}

\begin{proof}[Proof of Theorem~\ref{thm:twosided}]
  Let $u$ be a positive weak solution of \eqref{equ:EL}. By \eqref{equ:lower} of Proposition \ref{prop:lower}, \eqref{equ:upperorigin} of Proposition \ref{prop:upperorigin}, and \eqref{equ:upperinfinity} of Proposition \ref{prop:upperinfinity}, there exist $\rho_0\in(0,1)$, $R_\infty\in(1,\infty)$, and $C_7\in[1,\infty)$, depending on $u$, such that
  \begin{align*}
    C_7^{-1}|x|^{-\mu_1}&\leqslant u(x)\leqslant C_7|x|^{-\mu_1}\quad\text{for a.e. } 0<|x|<\rho_0,\\
    C_7^{-1}|x|^{-\mu_2}&\leqslant u(x)\leqslant C_7|x|^{-\mu_2}\quad\text{for a.e. } |x|>R_\infty.
  \end{align*}
  For $0<|x|<\rho_0$ we have $1\leqslant(1+|x|)^{\mu_2-\mu_1}\leqslant 2^{\mu_2-\mu_1}$, and for $|x|>R_\infty$ we have $|x|^{\mu_2-\mu_1}\leqslant(1+|x|)^{\mu_2-\mu_1}\leqslant 2^{\mu_2-\mu_1}|x|^{\mu_2-\mu_1}$. Hence \eqref{equ:twosidedglobal} holds for a.e. $x$ with $|x|<\rho_0$ or $|x|>R_\infty$, with $C:=2^{\mu_2-\mu_1}C_7$. On the closed annulus $\{\rho_0\leqslant|x|\leqslant R_\infty\}$, the function $|x|^{\mu_1}(1+|x|)^{\mu_2-\mu_1}$ is bounded from above and below by positive constants, and, by Corollary \ref{cor:annulus}, so is $u$. Enlarging $C$, we obtain \eqref{equ:twosidedglobal} for a.e. $x\in\mbRN$. This proves Theorem \ref{thm:twosided}\ref{item:twosidedpointwise}.

  By \eqref{equ:twosidedglobal}, for a.e. $x\in B_1$ we have $2^{-(\mu_2-\mu_1)}C^{-1}|x|^{-\mu_1}\leqslant u(x)\leqslant C|x|^{-\mu_1}$, and, for a.e. $x\in\mbRN\setminus B_1$, $2^{-(\mu_2-\mu_1)}C^{-1}|x|^{-\mu_2}\leqslant u(x)\leqslant C|x|^{-\mu_2}$. Since $\mu_1\pst<N<\mu_2\pst$ by \eqref{equ:exponentfacts}, integration in polar coordinates gives, for $R\in(0,1]$,
  \begin{equation*}
    \int_{B_R}|x|^{-\mu_1\pst}\diopx=\frac{\omega_{N-1}}{N-\mu_1\pst}R^{N-\mu_1\pst},
  \end{equation*}
  and, for $R\in[1,\infty)$,
  \begin{equation*}
    \int_{\mbRN\setminus B_R}|x|^{-\mu_2\pst}\diopx=\frac{\omega_{N-1}}{\mu_2\pst-N}R^{-(\mu_2\pst-N)}.
  \end{equation*}
  Therefore \eqref{equ:massorigin} and \eqref{equ:massinfinity} hold after a further enlargement of $C$, which proves Theorem \ref{thm:twosided}\ref{item:twosidedmass}.
\end{proof}

\begin{proof}[Proof of Proposition~\ref{prop:scaling}]
  Let $t\in(0,\infty)$. Clearly $u^{(t)}>0$ a.e. The change of variables $x\mapsto tx$, $y\mapsto ty$ gives $\semi{u^{(t)}}^p=t^{pr_H+ps-N}\semi{u}^p=\semi{u}^p$, because $pr_H=N-ps$, and $\norm{u^{(t)}}_{L^{\pst}(\mbRN)}=\norm{u}_{L^{\pst}(\mbRN)}$, because $r_H\pst=N$; hence $u^{(t)}\in\mcD$ by \eqref{equ:spaceD}. Let $\vp\in\mcD$ and $\vp_t(x):=\vp(\frac{x}{t})$. The same change of variables and \eqref{equ:weakEL} tested with $\vp_t$ give
  \begin{align*}
    \langle(-\Delta_p)^su^{(t)},\vp\rangle
    &=t^{r_H(p-1)+ps-N}\langle(-\Delta_p)^su,\vp_t\rangle\\
    &=\mu\int_{\mbRN}\frac{u^{(t)}(x)^{p-1}\vp(x)}{|x|^{ps}}\diopx+t^{r_H(p-1)+ps-r_H(\pst-1)}\int_{\mbRN}u^{(t)}(x)^{\pst-1}\vp(x)\diopx.
  \end{align*}
  Since $r_H(\pst-p)=ps$, the last exponent is $0$. Hence $u^{(t)}$ is a positive weak solution of \eqref{equ:EL} with $\semi{u^{(t)}}=\semi{u}$.

  Let $C$ be the constant in Theorem \ref{thm:twosided} for $u$, and let $\rho\in(0,\infty)$. For a.e. $x\neq 0$,
  \begin{equation*}
    |x|^{\mu_1}u^{(t)}(x)=t^{r_H-\mu_1}|tx|^{\mu_1}u(tx)
    \qquad\text{and}\qquad
    |x|^{\mu_2}u^{(t)}(x)=t^{r_H-\mu_2}|tx|^{\mu_2}u(tx).
  \end{equation*}
  By \eqref{equ:twosidedglobal}, $2^{-(\mu_2-\mu_1)}C^{-1}\leqslant|y|^{\mu_1}u(y)\leqslant C$ for a.e. $y\in B_1$, and $2^{-(\mu_2-\mu_1)}C^{-1}\leqslant|y|^{\mu_2}u(y)\leqslant C$ for a.e. $y\in\mbRN\setminus B_1$. If $t\geqslant\frac{1}{\rho}$, then $\{tx\colon 0<|x|<\rho\}\supset B_1\setminus\{0\}$, and hence
  \begin{equation*}
    \esssup_{\{0<|x|<\rho\}}|x|^{\mu_1}u^{(t)}(x)\geqslant 2^{-(\mu_2-\mu_1)}C^{-1}t^{r_H-\mu_1}\longrightarrow\infty\quad\text{as } t\to\infty,
  \end{equation*}
  since $r_H-\mu_1>0$. If $t\leqslant\frac{1}{\rho}$, then $\{tx\colon 0<|x|<\rho\}\subset B_1$, and hence $\essinf_{\{0<|x|<\rho\}}|x|^{\mu_1}u^{(t)}(x)\leqslant Ct^{r_H-\mu_1}\to 0$ as $t\to 0^+$. This proves \eqref{equ:scalingorigin}. Similarly, suppose first that $t\leqslant\frac{1}{\rho}$. Then $\{tx\colon|x|>\rho\}\supset\mbRN\setminus B_1$, so that
\begin{equation*}
  \esssup_{\{|x|>\rho\}}|x|^{\mu_2}u^{(t)}(x)\geqslant 2^{-(\mu_2-\mu_1)}C^{-1}t^{r_H-\mu_2}\longrightarrow\infty
  \quad\text{as } t\to 0^+,
\end{equation*}
because $r_H-\mu_2<0$. If $t\geqslant\frac{1}{\rho}$, then $\{tx\colon|x|>\rho\}\subset\mbRN\setminus B_1$. Hence
\begin{equation*}
  \essinf_{\{|x|>\rho\}}|x|^{\mu_2}u^{(t)}(x)\leqslant Ct^{r_H-\mu_2}\longrightarrow 0
  \quad\text{as } t\to\infty.
\end{equation*} This proves \eqref{equ:scalinginfinity} and completes the proof of Proposition \ref{prop:scaling}.
\end{proof}

{\bf Acknowledgment.}
The first author is supported by the National Natural Science Foundation of China (Grant No. 12571127). The second author was supported by the National Natural Science Foundation of China (Grant No. 12526592). \\

{\bf Conflict of interest.}
The authors declare that they have no conflict of interest. \\

{\bf Declaration on the use of AI tools.}
AI-based tools were used for: (i) the organization and wording of the text; and (ii) the verification of mathematical arguments and computations.  The mathematical statements and their proofs were reviewed by
the authors, who take full responsibility for the content of the manuscript.

\end{document}